\documentclass[a4paper, final]{amsart}

\usepackage[T1]{fontenc}
\usepackage[utf8]{inputenc}
\usepackage{lmodern}
\usepackage[ngerman, english]{babel}

\usepackage[only,mapsfrom]{stmaryrd}
\usepackage{amsmath}
\usepackage{amsthm}
\usepackage{amssymb}
\usepackage{mathtools}
\usepackage{mathdots}
\usepackage{amsfonts}
\usepackage{centernot}
\usepackage{thmtools}
\usepackage{ytableau}
\usepackage{multirow}
\usepackage{multicol}
\usepackage[inline]{enumitem}

\usepackage{url}
\usepackage[numbers,sort]{natbib}

\usepackage{xcolor}
\usepackage{setspace}
\usepackage{graphicx}
\usepackage{tikz}
	\usetikzlibrary{cd}
\usepackage{faktor}
\usepackage{etoolbox}
\usepackage{makecell}

\usepackage{environ} %

\usepackage[colorlinks=true, linkcolor=black, citecolor=black, urlcolor=blue]{hyperref}
\usepackage{cleveref}

\crefname{thm}{theorem}{theorems}
\crefname{prop}{proposition}{propositions}
\crefname{lem}{lemma}{lemmas}
\crefname{cor}{corollary}{corollaries}
\crefname{defi}{definition}{definitions}
\crefname{exa}{example}{examples}
\crefname{notn}{notation}{notations}
\crefname{figure}{figure}{figures}
\crefname{section}{section}{sections}
\crefname{subsection}{subsection}{subsections}
\crefname{table}{table}{tables}
\crefname{equation}{}{}

\makeatletter
\g@addto@macro\@floatboxreset\centering
\makeatother

\newtheorem{thm}{Theorem}[section]
\newtheorem{cor}[thm]{Corollary}
\newtheorem{prop}[thm]{Proposition}
\newtheorem{lem}[thm]{Lemma}

\newtheorem*{cla*}{Claim}
\newtheorem{defi}[thm]{Definition}

\theoremstyle{definition}

\newtheorem{rem}[thm]{Remark}
\newtheorem*{rem*}{Remark}
\newtheorem{exa}[thm]{Example}

\newtheorem{notn}[thm]{Notation}

\newtheoremstyle{efronremark}{6pt}{6pt}{}{}{\itshape}{\quad}{ }{\thmnote{#3}}
\theoremstyle{efronremark}

\newtheoremstyle{longproofstyle}%
{3pt}%
{3pt}%
{\color{darkgray}}%
{}%
{\itshape\color{darkgray}}%
{.}%
{.5em}%
{}%

\theoremstyle{longproofstyle}

\numberwithin{equation}{section}

\ytableausetup{centertableaux, boxsize=1.2em}

\renewcommand{\bf}[1]{\textbf{#1} }

\newcommand{\bs}[1]{\boldsymbol{#1} }
\newcommand{\mc}[1]{\mathcal{#1} }

\newcommand{\N}{\mathbb N}

\newcommand{\Z}{\mathbb Z}

\newcommand{\K}{\mathbb K}

\newcommand{\set}[1]{\left\{ #1 \right\}}
\newcommand{\card}[1]{\left| #1 \right|}

\newcommand{\id}{\operatorname{id}}

\newcommand{\cl}{\operatorname{cl}}

\newcommand{\SG}{\mathfrak{S}} %

\newcommand{\fs}{s}
\newcommand{\myif}{\text {if }}

\newcommand{\He}{{H_n(0)}}
\newcommand{\Hec}{{H_W(0)}}
\newcommand{\Heq}{{H_W(q)}}

\newcommand{\ngen}{T}

\newcommand{\parts}[2]{{(#1_1, \dots, #1_{#2})}}

\newcommand{\fakt}[3]{\faktor{{#1}_{#2,#3}}{\approx_{#2}}}
\newcommand{\faktid}[2]{\faktor{{#1}_{#2}}{\approx}}

\newcommand{\wfakt}[2]{\fakt{W}{#1}{#2}}
\newcommand{\wfaktid}[1]{\faktid{W}{#1}}
\newcommand{\wmaxa}{\wfaktid{\max}}

\newcommand{\sfakt}[2]{\fakt{(\SG_n)}{#1}{#2}}
\newcommand{\sfaktid}[1]{\faktid{(\SG_n)}{#1}}
\newcommand{\smaxa}{\sfaktid{\max}}

\newcommand{\del}{\operatorname{del}}
\newcommand{\ins}{\operatorname{ins}}
\newcommand{\cst}{\operatorname{cst}}

\newcommand{\inv}{^{-1}}
\newcommand{\Inv}{\operatorname{Inv}}
\newcommand{\rnk}{\rho}

\newcommand{\partition}{P}
\newcommand{\iprod}{\mathbin{\odot}}

\newcommand{\typems}{C}
\newcommand{\hdash}{\vDash_e}
\newcommand{\ord}{\operatorname{ord}}
\newcommand{\na}{\nu} %

\DeclarePairedDelimiter\ceil{\lceil}{\rceil}
\DeclarePairedDelimiter\floor{\lfloor}{\rfloor}

\newcommand\cf{cf.\kern.3em}
\newcommand\eg{e.g.\kern.3em}
\newcommand\ie{i.e.\kern.3em}
\newcommand\resp{\text{resp.}\kern.3em}
\newcommand\etc{etc.\kern.3em}
\newcommand\etal{et al.\kern.3em}

\setlist{                %
	topsep = 0.5ex,
	itemsep = 0.5ex,
	parsep = 0ex
	}
	
\newlist{wideenumerate}{enumerate}{5}
\newlist{caseenum}{enumerate}{3}

\setlist[enumerate,1]{label = \textup{(}\arabic*\textup{)} }
\setlist[enumerate,2]{label = \textup{(}\roman*\textup{)} }

\setlist[wideenumerate]{label = \textup{(}\arabic*\textup{)}, wide, nosep}
\setlist[wideenumerate, 2]{label = \textup{(}\roman*\textup{)}}

\setlist[caseenum]{label = \textbf{Case \arabic*.}, wide, nosep}
\setlist[caseenum, 2]{label = \textbf{Case \arabic{caseenumi}.\arabic*.}, wide, nosep}
\setlist[caseenum, 3]{label = \textbf{Case \arabic{caseenumi}.\arabic{caseenumii}.\arabic*.}, wide, nosep}

\definecolor{g1}{gray}{0.85}
\definecolor{g2}{gray}{0.7}
\definecolor{g3}{gray}{0.55}
\definecolor{g4}{gray}{0.4}
\definecolor{g5}{gray}{0.25}
\definecolor{g6}{gray}{0.1}

\title{
	Combinatorics of Centers of 0-Hecke Algebras in Type A 
	}
\author{Sebastian König}
\address{
	Leibniz Universität Hannover \\	
	Institute of Algebra, Number Theory and Discrete Mathematics \\
	Welfengarten 1 \\
	30167 Hannover \\
	Germany
}
\email{sebastian.koenig@math.uni-hannover.de}

\hypersetup{%
	pdfauthor = {Sebastian König},
	pdftitle = {Combinatorics of Centers of 0-Hecke Algebras in Type A},
	pdfkeywords = {0-Hecke algebra, center, symmetric group, conjugacy class}
}

\begin{document}

\begin{abstract} 		 
	A basis of the center of the 0-Hecke algebra of an arbitrary finite Coxeter group was described by He in 2015.
	This basis is indexed by certain equivalence classes of the Coxeter group whose explicit description is rather complicated. 
	Even their number is not obvious.
	
	We consider case of the symmetric group $\SG_n$.
	Building on work of Geck, Kim and Pfeiffer, we obtain a complete set of representatives of the equivalence classes.
	This set is naturally parametrized by certain compositions of $n$ called \emph{maximal}.
	It follows that the maximal compositions of $n$ index the basis of the center of the 0-Hecke algebra of $\SG_n$.
	We then develop an explicit combinatorial description for the equivalence classes that are parametrized by the maximal compositions whose odd parts form a hook.
\end{abstract}

\maketitle

\section{Introduction}

Let $W$ be a finite Coxeter group.
The Iwahori-Hecke algebra $\Heq$ of $W$ is a  deformation of the group algebra of $W$ with nonzero parameter $q$.
Iwahori-Hecke algebras arise in the representation theory of finite groups of Lie type and Knot theory \cite{Geck2000a}.
Setting $q=0$ results in the $0$-Hecke algebra $\Hec$.
A first (and thorough) study of  $\Hec$ was carried out by Norton \cite{Norton1979}.
Its structure diverges considerably from the generic $q\neq 0$ case \cite{Carter1986}.
The 0-Hecke algebras appear in the modular representation theory of finite groups of Lie type   \cite{Norton1979,Carter1976}.
The Grothendieck ring of the finitely generated modules of the 0-Hecke algebras of the symmetric groups is isomorphic to the Hopf algebra of quasisymmetric functions \cite{Krob1997}. 
This article is concerned with the center $Z(\He)$ of the 0-Hecke algebra $\He$ of the symmetric group $\SG_n$.

Fayers mentions the description of $Z(\Hec)$ as an open problem in \cite{Fayers2005}.
Brichard gives a formula for the dimension of the center in type $A$ \cite{Brichard2008}.
Yang and Li obtain a lower bound for the dimension of $Z(\Hec)$ for irreducible $W$ in several types other than $A$ \cite{Yang2015}.
Moreover, they specify the dimension in type $I_2(n)$ for $n\geq 5$.
In \cite{He2015} He describes a basis of $Z(\Hec)$ in arbitrary type indexed by a set of equivalence classes $\wmaxa$ of $W$. 
These classes are rather subtle. 
In fact, \cite{He2015} contains no result on the number of these classes which is the dimension of  the center.
Motivated by the connection to the center of the 0-Hecke algbera $\He$, the present paper aims to shed some light on the quotient set $\smaxa$. 
To be precise, we parametrize it by compositions, obtain a complete set of representatives and develop a combinatorial description for certain elements of $\smaxa$.

Let $S$ be the set of Coxeter generators of $W$ and $\ell$ be the length function of $W$.
Define $W_{\min}$ and $W_{\max}$ to be the set of elements of~$W$ whose length is minimal and maximal in their conjugacy class, respectively.
Geck and Pfeiffer introduce in \cite{Geck1993} a relation $\to$ on $W$ known as \emph{cyclic shift relation}.
It is the reflexive and transitive closure of the relations $\overset{s}{\to}$  for $s\in S$ where we have $w\overset{s}{\to}w'$ if $w' = sws$ and $\ell(w')\leq \ell(w)$.

In the case where $W$ is a Weyl group, Geck and Pfeiffer show that $W_{\min}$ in conjunction with $\to$ has remarkable properties and how these  properties can be used in order to define a character table of $\Heq$ with $q\neq 0$ \cite{Geck1993}.
Since then their results have been generalized to finite \cite{Geck1996}, affine \cite{He2014} and finally to all Coxeter groups \cite{Marquis2021}.
The relation $\to$ can also be used to describe the conjugacy classes of Coxeter groups \cite{Geck2000a,Marquis2020structure} in particular for computational purposes \cite{Geck1996,Geck2000a}.
Geck, Kim and Pfeiffer introduce a twisted version $\to_\delta$ of the relation belonging to twisted conjugacy classes of $W$ in \cite{Geck2000}. 
Building on the results of \cite{Geck1993}, Geck and Rouquier define a basis of $Z(\Heq)$ for $q\neq 0$ and $W$ a finite Weyl group, which is naturally indexed by the conjugacy classes of $W$ \cite{Geck1997}.

By setting $w\approx w'$ if and only if $w\to w'$ and $w' \to w$ one obtains an equivalence relation $\approx$ on $W$. 
The $\approx$-equivalence classes of $W$ are known as \emph{cyclic shift classes}.
On $\faktid{W}{}$ the relation $\to$ gives rise to a partial order. 
Gill considers the corresponding subposets $\faktid{\mc O}{}$ where $\mc O$ is a conjugacy class of $W$ \cite{Gill2000}.

For an element $\Sigma$ of the quotient set $\wmaxa$, He defines the element $\ngen_{\leq \Sigma} := \sum_x \ngen_x$ where $x$ runs over  the order ideal in Bruhat order of $W$ generated by $\Sigma$ \cite{He2015}.
Then he shows that the elements $\ngen_{\leq \Sigma}$ for $\Sigma \in \wmaxa$ form a basis of $Z(\Hec)$.
We consider He's approach in \Cref{sec:preliminaries}.

For each composition $\alpha \vDash n$, Kim defines the \emph{element in stair form} $\sigma_\alpha \in \SG_n$ \cite{Kim1998}. 
Moreover, she calls $\alpha \vDash n$ \emph{maximal} if there is $k \geq 0$ such that the first~$k$ parts of $\alpha$ are even and the remaining parts are odd and weakly decreasing. 
In this case we write $\alpha \vDash_e n$.
We show in \Cref{thm:parametrizations_of_kim} that the elements in stair form $\sigma_\alpha$ for $\alpha \vDash_e n$ form a system of representatives 
of $\smaxa$. 
For $\alpha \vDash_e n$ let $\Sigma_\alpha \in \smaxa$ be the equivalence class of the element in stair form $\sigma_\alpha$.
It follows that $\alpha \mapsto \Sigma_\alpha$ is a parametrization of $\smaxa$ by the maximal compositions of $n$.
That is, the elements $\ngen_{\leq \Sigma_\alpha}$ for $\alpha \vDash_e n$ form a basis of $Z(\He)$.
This leads to an alternative proof of Brichards dimension formula from \cite{Brichard2008}, which she obtained by considering braid diagrams on the Möbius strip.

Since $\ngen_{\leq\Sigma_\alpha}$ depends on the order ideal generated by $\Sigma_\alpha$, a description of the elements of $\Sigma_\alpha$ is desirable.
We obtain combinatorial characterizations of the equivalence classes $\Sigma_{(n)}$ (\Cref{thm:characterization_of_Sigma_(n)}) and $\Sigma_{(k,1^{n-k})}$ with~$k$ odd (\Cref{thm:characterization_of_Sigma_for_odd_hook}) and a decomposition rule $\Sigma_{(\alpha_1,\dots, \alpha_l)} =\Sigma_{(\alpha_1)} \iprod \Sigma_{(\alpha_2,\dots,\alpha_l)}$ if $\alpha_1$ is even given by an injective operator $\iprod$ which we call the \emph{inductive product} (\Cref{thm:indcutive_product_bijection}).
This allows us to describe $\Sigma_\alpha$ for all $\alpha \vDash_e n$ whose odd parts form a hook.
Moreover, we will see how these $\Sigma_\alpha$ can be computed recursively. 
In particular, we infer a characterization of $\Sigma_{(k,1^{n-k})}$ for even $k$ as well (\Cref{thm:characterization_of_Sigma_for_arbitrary_hook}).

The structure is as follows.
In \Cref{sec:preliminaries} we present the background material and review He's basis of the center of $\Hec$ and the connection to the quotient set $\wmaxa$.
In \Cref{sec:parametrization} we obtain the parametrization by maximal compositions and a set of representatives of $\smaxa$.
In \Cref{sec:equivalence_classes} we combinatorially describe the elements $\Sigma_\alpha \in\smaxa$ indexed by maximal compositions whose odd parts form a hook.

\section{Preliminaries}
\label{sec:preliminaries}

Let~$\K$ be an arbitrary field.
We set $\N := \set{1,2,\dots}$ and always assume that $n\in \N$.
For $a,b\in \Z$ we define the \emph{discrete interval} $[a,b]:=\set{c\in \Z\mid a\leq c \leq b}$ and use the shorthand $[a] := [1,a]$.

\subsection{Coxeter groups}

We consider basic concepts from the theory of finite Coxeter groups. 
Our motivation is the application to the symmetric groups.
Refer to \cite{Humphreys1990,Bjorner2006} for details.

Let~$S$ be a set.
A \emph{Coxeter matrix} is a map $m \colon S \times S \to \N \cup \set{\infty}$ such that
\begin{enumerate*} 	
	\item $m(s,s') = 1$ if and only if $s' = s$ and 
	\item $m(s,s') = m(s',s)$
\end{enumerate*}
for all $s,s'\in S$.
The corresponding \emph{Coxeter graph} is the undirected graph with vertex set $S$ containing the edge $\set{s,s'}$ if and only if $m(s,s') \geq 3$. If $m(s,s') \geq 4$ then the edge $\set{s,s'}$ is labeled with $m(s,s')$. 
A group $W$ is called \emph{Coxeter group} with \emph{Coxeter generators} $S$ if $W$ is generated by $S$ subject to the relations
\begin{enumerate}
	\item 
	$s^2 = 1$ for all $s\in S$,
	\item
	$(ss's \cdots)_{m(s,s')} = (s'ss' \cdots)_{m(s,s')}$ for all $s,s'\in S$ with $s\neq s'$ and\\ $m(s,s') < \infty$
\end{enumerate}
where $(ss's\cdots)_p$ denotes the the alternating product of $s$ and $s'$ with $p$ factors.

Let $W$ be a Coxeter group with Coxeter generators $S$.
We always assume that $W$ is finite.
For $I\subseteq S$ the \emph{parabolic subgroup} $W_I$ is the subgroup of $W$ generated by $I$.
It is a Coxeter group with Coxeter generators $I$.

Each $w\in W$ can be written as a product $w = s_1 \cdots s_k$ with $s_i \in S$.
Then $s_1 \cdots s_k$ is called a \emph{word} for $w$.
If~$k$ is minimal among all words for $w$, $s_{1} \cdots s_{k}$ is a \emph{reduced word} for $w$ and $\ell(w) := k$ is the \emph{length} of $w$.
The \emph{left} and the \emph{right descent set} of
$w\in W$ are given by
\begin{align}
	\begin{aligned}
	D_L(w) 
	&:= \set{s\in S \mid \ell(sw) < \ell(w)}, \\
	D_R(w) 
	&:= \set{s\in S \mid \ell(w s) < \ell(w)}.
	\end{aligned}
	\label{eq:descent_sets_W}
\end{align}

The \emph{Bruhat order} $\leq$ is the partial order on $W$ given by $u \leq w$ if and only if there exists a reduced word for $w$ which contains a reduced word of $u$ as a subsequence.
The Bruhat poset is graded by the length function $\ell$.
Since $W$ is finite, there exists a greatest element $w_0 \in W$ in Bruhat order. %
This element is called the \emph{longest element} of $W$.
It has the following useful properties.

\begin{lem}[{\cite[Proposition 2.3.2 and Corollary 2.3.3]{Bjorner2006}}]
	\label{thm:w0_properties}
	 Let $w_0$ be the longest element of $W$.
	 Then we have
	\begin{enumerate}
		\item $w_0^2 = 1$,
		\item $\ell(ww_0) = \ell(w_0w) = \ell(w_0) - \ell(w)$ for all $w \in W$,
		\item $\ell(w_0ww_0) = \ell(w)$ for all $w\in W$.
	\end{enumerate}
\end{lem}

\begin{lem}[{\cite[Propositions 2.3.4 and 3.1.5]{Bjorner2006}}]
\label{thm:w0_maps}
For the Bruhat order on $W$, we have that
\begin{enumerate}
\item  $w\mapsto ww_0$ and $w\mapsto w_0w$ are antiautomorphisms,
\item $w\mapsto w_0 w w_0$ is an automorphism.
\end{enumerate}
\end{lem}

We now define the 0-Hecke algebra of $W$.
Refer to Chapter~1 of \cite{Mathas1999} for background information on $\Hec$.

\begin{defi}
	\label{thm:0-Hecke algebra_Coxeter_group}
	The \emph{0-Hecke algebra}~$\Hec$ of $W$ is the unital associative $\K$-algebra generated by the elements $T_s$ for $s\in S$ subject to the relations
	\begin{enumerate}
	\item 
	$\ngen_s^2 = -\ngen_s,$
	\item
	$(\ngen_{s}\ngen_{s'}\ngen_{s} \cdots)_{m(s,s')} = (\ngen_{s'}\ngen_{s}\ngen_{s'} \cdots)_{m(s,s')}$ for all $s,s'\in S$ with $s\neq s'$.
	\end{enumerate}
\end{defi}

For $w \in W$ define $\ngen_w := \ngen_{s_1} \cdots \ngen_{s_k}$ where $\fs_{1} \cdots \fs_{k}$ is a reduced word for $w$.
The word property \cite[Theorem 3.3.1]{Bjorner2006} ensures that this is well defined. 
By \cite[Theorem 1.13]{Mathas1999}, we have that $\set{\ngen_{w} \mid w \in W}$ is a $\K$-basis of $\Hec$ with multiplication given by
\begin{align*}
\ngen_s \ngen_w = 
 \begin{cases}
	 \ngen_{s w} & \text{if}\ \ell(s w ) > \ell(w) \\
	 - \ngen_{w}    & \text{if}\ \ell(s w ) < \ell(w) \\
 \end{cases}
\end{align*}
 for $w\in W$ and $s\in S$.

\subsection{The symmetric group}

For a finite set $X$ we define $\SG(X)$ to be the group formed by all bijections from $X$ to itself.
The \emph{symmetric group} $\SG_n$ is the group $\SG([n])$.
Its elements are called \emph{permutations}.

Let $S$ be the set of adjacent transpositions $s_i := (i, i+1) \in \SG_n$ for $i = 1,\dots, n-1$.
The elements of $S$ generate $\SG_n$ as a Coxeter group subject to the relations
\begin{alignat*}{3}
\fs_i^2 &= 1,   			\\
\fs_i\fs_{i+1}\fs_i &= \fs_{i+1}\fs_i\fs_{i+1}, 	\\
\fs_i\fs_j &= \fs_j \fs_i 		\ \text{if } |i-j|\geq 2
\end{alignat*}
\cite[Proposition 1.5.4]{Bjorner2006}.
For $n\geq 2$, $\SG_n$ is an irreducible Coxeter group of type $A_{n-1}$.
While considering the symmetric group $\SG_n$, we always assume that $S$ is the corresponding set of adjacent transpositions.
For $\sigma \in \SG_n$ we have
\begin{align}
	\begin{aligned}
		D_L(\sigma) &= \set{s_i \in S \mid \sigma^{-1}(i) > \sigma^{-1}(i+1)}, \\
		D_R(\sigma) &= \set{s_i \in S \mid \sigma(i) > \sigma(i+1)}
	\end{aligned}
	\label{eq:descent_sets_of_Sn}
\end{align}
\cite[Proposition 1.5.3]{Bjorner2006}.
The longest element $w_0$ of $\SG_n$ is given by $w_0(i) = n-i+1$ for $i\in [n]$.
We denote the \emph{$0$-Hecke algebra of the symmetric group} $\SG_n$ with $\He := H_{\SG_n}(0)$ and use the shorthand $\ngen_i := \ngen_{s_i}$ for $i\in [n-1]$.

\subsection{Combinatorics}

A \emph{composition} $\alpha = \parts{\alpha}l$ is a finite sequence of positive integers. 
The \emph{length} and the \emph{size} of $\alpha$ are given by $\ell(\alpha):= l$ and $|\alpha| := \sum_{i=1}^l \alpha_i$, respectively.
The $\alpha_i$ are called \emph{parts} of $\alpha$.
If $\alpha$ has size~$n$, $\alpha$ is called \emph{composition of~$n$} and we write $\alpha \vDash n$.
A \emph{weak} composition of $n$ is a finite sequence of nonnegative integers that sum up to $n$. We write $\alpha \vDash_0 n$ if $\alpha$ is a weak composition of $n$.
The \emph{empty composition} $\emptyset$ is the unique composition of length and size~$0$.
A \emph{partition} is a composition whose parts are weakly decreasing.
We write $\lambda \vdash n$ if $\lambda$ is a partition of size~$n$.
For example, $(1,4,3)\vDash 8$ and  $(4,3,1)\vdash 8$.
Partitions of $n$ of the form $(k, 1^{n-k})$ with $k\in [n]$ are called \emph{hooks}.

A permutation $\sigma\in \SG_n$ can be represented in cycle notation where
cycles of length one may be omitted.
The \emph{cycle type} (or simply \emph{type}) of a permutation $\sigma\in \SG_n$ is the partition of~$n$ whose parts are the sizes of all the cycles of $\sigma$.
If $\sigma$ has cycle type $(k,1^{n-k})$ for a $k\in [n]$ we also call it a $k$-cycle.
A $k$-cycle is \emph{trivial} if $k = 1$.
Writing $\sigma$ in cycle notation is the same as expanding $\sigma$ into a product $\sigma_1 \cdots \sigma_r$ of disjoint cycles where the trivial cycles may be omitted  in the expansion.
On the other hand, in order to describe the cycle notation of a permutation combinatorially, it can be useful to include them.
In \Cref{sec:equivalence_classes} we will characterize the elements of certain equivalence classes of $\SG_n$ by considering them in cycle notation.

\subsection{Centers of 0-Hecke algebras}

\label{sec:preliminaries:center}

In this \namecref{sec:preliminaries:center} we introduce He's basis of the center of $\He$.
Following his approach in \cite{He2015}, we take a more general point of view and consider the center of $\Hec$ for a finite Coxeter group $W$ twisted by an automorphism $\delta$.
This enables us to prove a useful invariance property in  \Cref{thm:conjugation_with_w0_and_equivalence_classes}.
By setting $W = \SG_n$ and $\delta = \id$ we recover the desired results on the center of $\He$.

Let~$W$ be a finite Coxeter group with Coxeter generators $S$ and $\delta$ be a~$W$-automorphism with $\delta(S) = S$.
For instance, we can chose $\delta=\id$.
Another example is given by the conjugation with $w_0$.
For $u,w\in W$ we use the shorthand $w^u = u w u\inv$.
Define $\nu\colon W \to W$, $w\mapsto w^{w_0}$.
Then $\nu$ is a group automorphism and from \Cref{thm:w0_properties} it follows that $\ell(\nu(w)) = \ell(w)$ for all $w\in W$ so that $\nu(S) = S$.
In general, each graph automorphism of the Coxeter graph of~$W$ gives rise to a~$W$-automorphism that fixes~$S$. 
By the next \namecref{thm:delta_is_diagram_automorphism}, the converse direction is also true.
The result is not new.
For instance, it was already used implicitly in \cite[Section 2.10]{Geck2000}.

\begin{lem}
	\label{thm:delta_is_diagram_automorphism}
	Let $\delta$ be a group automorphism of $W$ with $\delta(S) = S$.
	\begin{enumerate}
		\item $\delta$ is an automorphism of the Coxeter graph of~$W$.
		\item $\delta$ is an automorphism of the Bruhat order of~$W$.
	\end{enumerate}
\end{lem}

\begin{proof}
	For $w\in W$  denote the order of~$w$ with $\ord(w)$.
	Let $m$ be the Coxeter matrix and $\Gamma$ be the Coxeter graph of~$W$.
	Then $m(s,s') = \ord(ss')$ for all $s,s' \in S$.
	Since $\delta$ is a group automorphism, we have $\ord(\delta(w)) = \ord(w)$ for all $w\in W$.
	Hence for all $s,s'\in S$ 
	\begin{align*}
		m(\delta(s),\delta(s')) 
		= \ord(\delta(s)\delta(s'))
		= \ord(s s')
		= m(s,s').
	\end{align*}
	Thus, $\delta$ is an automorphism of $\Gamma$.

	 By a comment following \cite[Proposition 2.3.4]{Bjorner2006}, we have that  multiplicatively extending a graph automorphism  of $\Gamma$ yields an Bruhat order automorphism of $W$.
	 Hence, $\delta$ is such an automorphism.
\end{proof}

\begin{exa}
	\label{thm:delta_for_Sn}
	
	The Coxeter graph of $\SG_n$ is shown below.
		
		\begin{center}
			\begin{tikzpicture}		
			\filldraw (0,0) circle (2pt) node (s1) {};
			\draw (s1)  node[align=center, above] {$s_1$};
			
			\filldraw (1,0) circle (2pt) node (s2) {};
			\draw (s2)  node[align=center, above] {$s_2$};
			
			\filldraw (2,0) circle (2pt) node (s3) {};
			\draw (s3)  node[align=center,   above] {$s_3$};
			
			\node (s4) at (2.5,0) {};
			\node (s5) at (3.0,0) {};
			
			\filldraw (3.5,0) circle (2pt) node (sn-2) {};
			\draw (sn-2)  node[align=center,   above] {$s_{n-2}$};
			
			\filldraw (4.5,0) circle (2pt) node (sn-1) {};
			\draw (sn-1)  node[align=center,   above] {$s_{n-1}$};
			
			\draw (s1) -- (s2) -- (s3) -- (s4);
			\draw[dotted] (s4) -- (s5);
			\draw (s5) -- (sn-2) -- (sn-1);;		
			\end{tikzpicture}
		\end{center}
		
	This graph has at most two automorphisms: the identity and the mapping given by $s_i \mapsto s_{n-i}$.
	For $n\geq 3$ these maps are distinct.
	Let $w_0$ be the longest element of $\SG_n$.
	Then $w_0(j) = n-j+1$ for all $j\in [n]$ and therefore $s_i^{w_0} = (n-i+1, n-i) = s_{n-i}$.
	Hence the second map is $\nu$.
	Thus, $\id$  and $\na$ are the only possibilities for $\delta$ if $W = \SG_n$.
\end{exa}

Two elements $w,w' \in W$ are called \emph{$\delta$-conjugate} if there is an $x\in W$ such that $w' = x w \delta(x)^{-1}$.
The set of $\delta$-conjugacy classes of~$W$ is denoted by $\cl(W)_\delta$.
For $\mc O \in \cl(W)_\delta$ the set of elements of minimal length in $\mc O$ and the set of elements of maximal length in $\mc O$ is denoted by $\mc O_{\min}$ and $\mc O_{\max}$, respectively.

We want to decompose these sets using an equivalence relation.
Let $w, w' \in W$. For $s \in S$ we write $w \overset{s}{\to}_\delta w'$ if $w' = s w \delta(s)$ and $\ell(w') \leq \ell(w)$.
We write $w \to_\delta w'$ if there is a sequence $w = w_1,w_2, \dots, w_{k+1} = w'$ of elements of~$W$ such that for each $i\in [k]$ there exists an $s\in S$ such that $w_{i} \overset{s}{\to}_\delta w_{i+1}$.
If $w\to_\delta w'$ and $w'\to_\delta w$ we write $w \approx_\delta w'$.
Then $\approx_\delta$ is an equivalence relation.
If $w \approx_\delta w'$ then $\ell(w) = \ell(w')$.
Thus, for all $\mc O \in \cl(W)_\delta$, $\mc O_{\min}$ and $\mc O_{\max}$ decompose into equivalence classes of $\approx_\delta$.
Define  $W_{\delta,\min} := \bigcup_{\mc O \in \cl(W)_\delta} \mc O_{\min}$ and $\wfakt{\delta}{\min}$
to be the quotient set of $W_{\delta,\min}$ by $\approx_\delta$.
Analogously, define the sets $W_{\delta,\max} := \bigcup_{\mc O \in \cl(W)_\delta} \mc O_{\max}$ and $\wfakt{\delta}{\max}$.
In the case $\delta = \id$ we may omit the index $\delta$.

\begin{exa}
	\label{exa:delta_quotient_sets}
	We have  $(1,2,3) \overset{(1,2)}{\to} (1,3,2) \overset{(1,2)}{\to} (1,2,3)$ so that $(1,2,3) \approx (1,3,2)$.
	Moreover, $\ell((1,2)) = \ell((2,3)) = 1$ and $\ell((1,3)) = 3$.
	Hence, 
	\begin{align*}
	\set{1},\ \set{(1,2,3),(1,3,2)} \ \text{and} \ \set{(1,3)}
	\end{align*}
	are the elements of $\faktor{(\SG_3)_{\max}}{\approx}$.
\end{exa}

Since $\delta$ is a Bruhat order automorphism of $W$ by  \Cref{thm:delta_is_diagram_automorphism},  we obtain an algebra automorphism of $\Hec$ by setting $\ngen_{s} \mapsto \ngen_{\delta(s)}$ for all $s\in S$ and extending multiplicatively and linearly.
This algebra automorphism is also denoted by $\delta$. 
The \emph{$\delta$-center} of $\Hec$ is given by
\begin{align*}
Z(\Hec)_\delta := \set{z \in \Hec \mid a z = z \delta(a) \text{ for all } a \in \Hec}.
\end{align*}

We now come to He's basis of  $Z(\Hec)_\delta$.
For $\Sigma \in \wfakt{\delta}{\max}$ set 
\begin{align*}
W_{\leq \Sigma} &:= \set{x\in W \mid x \leq w \text{ for some } w \in \Sigma}
\end{align*}
and
\begin{align*}
\ngen_{\leq \Sigma} &:= \sum_{x\in W_{\leq \Sigma}} \ngen_x.
\end{align*}

\begin{thm}[{\cite[Theorem 5.4]{He2015}}]
	\label{thm:basis_of_center_general}
	The elements $\ngen_{\leq \Sigma}$  for $\Sigma \in \wfakt{\delta}{\max}$ form a $\K$-basis of $Z(\Hec)_\delta$.
\end{thm}

We are concerned with the following special case.

\begin{cor}
	\label{thm:basis_of_center}
	The elements $\ngen_{\leq \Sigma}$  for $\Sigma \in \smaxa$ form a basis of $Z(\He)$.
\end{cor}

\begin{exa}
	Note that in $\SG_3$
	\begin{align*}
		(1,2,3) = s_1s_2, \ (1,3,2) = s_2s_1 \ \text{and} \ (1,3) = w_0.
	\end{align*} 
	Thus, \Cref{exa:delta_quotient_sets} and \Cref{thm:basis_of_center} yield that the elements
	\begin{align*}
		1,\ 1+ \ngen_1 + \ngen_2 + \ngen_1\ngen_2 + \ngen_2\ngen_1 \ \text{and} \ \sum_{w \in \SG_3} \ngen_w
	\end{align*}
	form a basis of $Z(H_3(0))$.
\end{exa}

The basis of $Z(\He)$ from \Cref{thm:basis_of_center} depends on $\smaxa$. 
This is the motivation for considering $\smaxa$ in this paper.
The remainder of this \namecref{sec:preliminaries} is devoted to show that $\na(\Sigma) = \Sigma$ for all $\Sigma\in \smaxa$.
This result will be useful  in \Cref{sec:equivalence_classes}.
In order to obtain it, we further study the quotient sets of $W_{\delta, \min}$ and $W_{\delta, \max}$ under $\approx_\delta$.

Define $\delta' := \na \circ \delta$.
Then $\delta'$ is a $W$-automorphism with $\delta'(S) = S$ as well.
The Bruhat order antiautomorphism $w\mapsto ww_0$ from \Cref{thm:w0_maps} relates  $\faktor{W_{\delta, \min}}{\approx_\delta}$ to $\faktor{W_{\delta', \max}}{{\approx_{\delta'}}}$.

\begin{lem}[{\cite[Section 2.2]{He2015}}]
	\label{thm:delta_min_and_delta_max}
	The map $\faktor{W_{\delta, \min}}{\approx_\delta} \to \faktor{W_{\delta', \max}}{{\approx_{\delta'}}}$, $\Sigma \mapsto \Sigma w_0$ is a bijection.
\end{lem}

We now come to parametrizations of $\wfakt{\delta}{\min}$ and $\wfakt{\delta}{\max}$ which are due to He. 
A $\delta$-conjugacy class $\mc O \in \cl(W)_\delta$ is called \emph{elliptic} (or \emph{cuspidal}) if $\mc O \cap W_I = \emptyset$ for all $I \subsetneq S$ such that $\delta(I) = I$.
Define
\begin{align*}
\Gamma_\delta := \set{(I,C) \mid \text{$I\subseteq S$, $I = \delta(I)$ and $C\in \cl(W_I)_\delta$ is elliptic} }.
\end{align*}

\begin{prop}[{\cite[Corollaries 4.2 and 4.3]{He2015}}]
	\label{thm:parametrizations_of_he}
	The maps
	\begin{align*}
	\begin{aligned}
	\Gamma_\delta &\to \wfakt{\delta}{\min} \\
	(I,C) &\mapsto C_{\min}
	\end{aligned}
	\quad \text{and} \quad
	\begin{aligned}
	\Gamma_{\delta'} &\to \wfakt{\delta}{\max} \\
	(I,C) &\mapsto C_{\min}w_0
	\end{aligned}
	\end{align*}
	are bijections. 
\end{prop}

The quotient sets $\wfakt{\delta}{\min}$ and $\wfakt{\delta}{\max}$ are given rather implicit both in their definitions and in \Cref{thm:parametrizations_of_he}.
The aim of this paper is to describe $\smaxa$ more explicitly.
In \Cref{sec:parametrization} we obtain a parametrization of $\smaxa$ by certain kinds of compositions and a set of representatives of $\smaxa$. 
Then in \Cref{sec:equivalence_classes} we combinatorially describe some of the elements of $\smaxa$.
Such matters are not discussed in \cite{He2015}.

From \Cref{thm:parametrizations_of_he} we deduce the following invariance properties.

\begin{lem} 
	\label{thm:delta_fixes_equivalence_classes}
	\hfil	
	\begin{enumerate}
		\item 
		We have $\delta(\Sigma) = \Sigma$ for each $\Sigma \in \wfakt{\delta}{\min}$.
		\item 
		We have $\delta'(\Sigma) = \Sigma$ for each $\Sigma \in \wfakt{\delta}{\max}$.
	\end{enumerate}
\end{lem}	

\begin{proof}
	\begin{wideenumerate}
	\item
	Let $\Sigma \in \wfakt{\delta}{\min}$ and $w\in \Sigma$.
	By \Cref{thm:parametrizations_of_he} there exists a tuple $(I,C) \in \Gamma_\delta$ such that $C \in \cl(W_I)_\delta$ and $\Sigma = C_{\min}$.
	Hence $w \in W_I$ and therefore $w^{-1} \in W_I$.
	It follows that
	\begin{align*}
		\delta(w) = w^{-1} w \delta(w^{-1})^{-1} \in C.
	\end{align*}
	Moreover,
	$\ell(\delta(w)) = \ell(w)$ because $\delta$ is a Bruhat order automorphism by \Cref{thm:delta_is_diagram_automorphism}.
	Therefore, $\delta(w) \in C_{\min{}} = \Sigma$.
	Hence, $\delta(\Sigma) = \Sigma$.
	\item
	Let $\Sigma \in \wfakt{\delta}{\max}$.
	From \Cref{thm:delta_min_and_delta_max} it follows that $\Sigma w_0 \in \wfakt{\delta'}{\min}$.
	Hence,
	\begin{align*}
		 \delta'(\Sigma)w_0 = \delta'(\Sigma w_0) = \Sigma w_0,
	\end{align*}
	where we use that $\delta'$ is a group homomorphism with $\delta'(w_0) = w_0$ for the first and Part~(1) for the second equality.
	Now multiply from the right with $w_0$.
	\qedhere
	\end{wideenumerate}
\end{proof}

Setting $W = \SG_n$ and $\delta = \id$ in the second part of  \Cref{thm:delta_fixes_equivalence_classes} yields the desired result on $\smaxa$ and $\na$.

\begin{cor}
	\label{thm:conjugation_with_w0_and_equivalence_classes}
	We have $\na(\Sigma) = \Sigma$ for each $\Sigma \in \smaxa$.
\end{cor}

\section{A parametrization by compositions}
\label{sec:parametrization}

The goal of this \namecref{sec:equivalence_classes}
is to obtain a new parametrization  of $\smaxa$ and, by \Cref{thm:basis_of_center}, of He's basis of $Z(\He)$.
The parameters will be compositions of the following kind.

\begin{defi}[{Kim, \cite{Kim1998}}]
	\label{def:maximal_composition}
	Let $\alpha = \parts{\alpha}{l} \vDash n$.
	We call $\alpha$ \emph{maximal} and write $\alpha \vDash_e n$ if there exists a $k$ with $0\leq k \leq l$ such that $\alpha_i$ is even for $i\leq k$, $\alpha_i$ is odd for $i>k$ and $\alpha_{k+1} \geq \alpha_{k+2} \geq \dots \geq \alpha_l$.
\end{defi}

For example, among the two compositions $(4, 6, 2, 3, 1, 1)$ and $(6, 4, 3, 2, 1, 1)$ of $17$ only the first one is maximal.
The parametrization corresponds to a set of representatives of $\smaxa$ given by
\emph{elements in stair form}:

\begin{defi}[{Kim, \cite{Kim1998}}]
	\label{thm:element_in_stair_form}
	Let $\alpha = \parts{\alpha}{l} \vDash n$. We define the list $(x_1, x_2, \dots, x_n)$ by setting $x_{2i-1} := i$ and  $x_{2i} := n-i+1$. The \emph{element in stair form  $\sigma_\alpha \in \SG_n$ corresponding to $\alpha$} is given by
	\begin{align*}
	\sigma_\alpha := \sigma_{\alpha_1} \sigma_{\alpha_2} \cdots \sigma_{\alpha_l}
	\end{align*}
	where $\sigma_{\alpha_i}$ is the $\alpha_i$-cycle
	\begin{align*}
	\sigma_{\alpha_i} := \left( x_{\alpha_1 + \dots + \alpha_{i-1} + 1}, x_{\alpha_1 + \dots + \alpha_{i-1} + 2}, \dots , x_{\alpha_1 + \dots + \alpha_{i-1} + \alpha_i} \right ).
	\end{align*}
\end{defi}

For instance, $\sigma_{(4,2)} = (1, 6, 2, 5)(3,4)$. 
We obtain $\sigma_\alpha$ for  $\alpha = \parts{\alpha}l\vDash n$ as follows.
 Let $d_i := \sum_{j=1}^{i}\alpha_i$ for $i = 1,\dots, l$ and consider the list $(x_1, x_2, \dots, x_n)$ given as above.
 Then split  the list between $x_{d_i}$ and $x_{d_i+1}$ for $i=1,\dots, l-1$. The resulting sublists are the cycles of $\sigma_\alpha$.
 In particular, if $\alpha$ and $\beta$ are compositions with $\sigma_\alpha = \sigma_\beta$ then $\alpha = \beta$.

The term \emph{maximal} is justified by the following result, which goes back to Kim \cite{Kim1998} and was proven by Geck, Kim and Pfeiffer \cite[Theorem 3.3]{Geck2000}.

\begin{lem}
	\label{thm:elememt_in_stair_form_maximality}
	Let $\alpha \vDash n$. Then $\sigma_\alpha \in (\SG_n)_{\max}$ if and only if $\alpha$ is a maximal composition.
\end{lem}

\begin{defi}
	For $\alpha \vDash_e n$ define $\Sigma_\alpha \in \smaxa$ to be the $\approx$-equivalence class of the element in stair form $\sigma_\alpha$.
\end{defi}

Thanks to \Cref{thm:elememt_in_stair_form_maximality} this is well defined.
We now come to the parametrization of $\smaxa$.

\begin{thm}
	\label{thm:parametrizations_of_kim}
	The map
	\begin{align*}
	\set{\alpha \vDash_e n}  \to \smaxa, \quad
		\alpha \mapsto \Sigma_\alpha
	\end{align*}
	is a bijection. 
	In particular, $\set{\sigma_\alpha \mid \alpha \vDash_e n}$ is a complete set of representatives of $\smaxa$.
\end{thm}

\begin{exa}
	For $n=3$ we have
	\begin{align*}
	\begin{array}{r|ccc}
	\alpha \vDash_e 3   & (3) 		& (2,1) & (1^3) \\ \hline
	\sigma_\alpha		& (1,3,2)	& (1,3)	& 1 	\\
	\end{array}
	\end{align*}
	which by \Cref{exa:delta_quotient_sets} is a complete set of representatives of  $\faktor{(\SG_3)_{\max}}{\approx}$.
\end{exa}

Before we begin with the proof of \Cref{thm:parametrizations_of_kim}, we discuss some immediate consequences.
By combining \Cref{thm:basis_of_center} and \Cref{thm:parametrizations_of_kim}, we obtain the following.

\begin{cor}
	\label{thm:basis_of_center_from_elements_in_stair_form}
	The elements $\ngen_{\leq \Sigma_\alpha}$  for $\alpha \vDash_e n$ form a basis of $Z(\He)$.
\end{cor}

This leads to an alternative proof of Brichards dimension formula.

\begin{cor}[{\cite[Section 5.1]{Brichard2008}}]
	\label{thm:dimension_of_center_and_cocenter}
	 The dimension of $Z(\He)$ equals
	\begin{align*}
	 \sum_{\lambda \vdash n} \frac{n_\lambda !}{m_\lambda}
	\end{align*}
	where for $\lambda = (1^{k_1}, 2^{k_2},\dots) \vdash n$, $m_\lambda := \prod_{i\geq 1} k_{2i}!$ and $n_\lambda := \sum_{i\geq 1} k_{2i}$ is the number of even parts of $\lambda$.
\end{cor}
\begin{proof}
	Each summand is the number of maximal compositions that have the same multiset of parts as $\lambda \vdash n$.
	Hence, the sum is the number of maximal compositions of~$n$.
	By \Cref{thm:basis_of_center_from_elements_in_stair_form} this is the dimension of $Z({\He})$.
\end{proof}

\begin{rem}
	From \Cref{thm:parametrizations_of_kim} and  \Cref{thm:delta_min_and_delta_max} it follows that we have a bijection  $\set{\alpha \vDash_e n}  \to \sfakt\na\min, \alpha \mapsto [\sigma_\alpha w_0]_\na$, where $[\sigma_\alpha w_0]_\na$ denotes the $\approx_\na$-equivalence class of $\sigma_\alpha w_0$.
	That is, the elements $\sigma_\alpha w_0$ for $\alpha \vDash_e n$ form a system of representatives of $\sfakt\na\min$.
	By {\cite[Theorem 6.5]{He2015}}, a basis of the $\na$-cocenter of $\He$ is given by such a system.
\end{rem}

We now come to the proof of \Cref{thm:parametrizations_of_kim}.
Because of \Cref{thm:elememt_in_stair_form_maximality}, it remains to show the following.
\begin{enumerate}[label = \textup{(}\alph*\textup{)}]
	\item
	\label{enum:parametrization_of_kim_surjectivity}
	For each $\Sigma \in \smaxa$ there is an $\alpha \vDash_e n$ such that $\sigma_\alpha \in \Sigma$.
	\item
	\label{enum:parametrization_of_kim_injectivity}
	If $\alpha,\beta \vDash_e n$ and $\sigma_\alpha \approx \sigma_\beta$ then $\alpha = \beta$.
\end{enumerate}

\Cref{thm:basis_of_center} and Brichards dimension formula, \Cref{thm:dimension_of_center_and_cocenter}, imply that
\begin{align*}
 \card{\smaxa} = \dim Z(\He) = \card{\set{\alpha \vDash_e n}}.
\end{align*}
Therefore, by simply citing the dimension formula from \cite{Brichard2008}, it would suffice to prove either \ref{enum:parametrization_of_kim_surjectivity} or \ref{enum:parametrization_of_kim_injectivity}.  
However, we chose to show both statements here as both proofs involve intermediate results that will be useful in the next \namecref{sec:equivalence_classes}.

In order to prove Statement~\ref{enum:parametrization_of_kim_surjectivity} we need the following result.

\begin{lem}
	\label{thm:arrow_and_length}
	Let $W$ be a finite Coxeter group and $w,w'\in W$ be such that $w\rightarrow w'$ and $\ell(w) = \ell(w')$.
	Then $w\approx w'$.
\end{lem}

\begin{proof}
	Let $S$ be the set of Coxeter generators of $W$.
	It suffices to consider the case where $w \overset{s}\rightarrow w'$ for some $s\in S$ because by definition $\rightarrow$ is the transitive closure of all the relations $\overset{t}\rightarrow$ with $t \in S$.
	Then $w' = s w s$.
	Thus, $w = s w' s$ and since $\ell(w) = \ell(w')$, we have $w'\overset{s}\rightarrow w$.
	Hence $w\approx w'$. 
\end{proof}

\begin{proof}[Proof of Statement~\ref{enum:parametrization_of_kim_surjectivity}]
	Let $\Sigma \in \smaxa$ and $\sigma\in \Sigma$.
	In \cite[Section 3]{Kim1998} it is shown that there is a $\beta \vDash n$ such that $\sigma_\beta \rightarrow \sigma$. 
	Moreover, Statement~(a$''$) of Section~3.1 in \cite{Geck2000} provides the existence of an $\alpha \vDash_e n$ such that $\sigma_\alpha \rightarrow \sigma_\beta$.
	Therefore, $\sigma_\alpha \rightarrow \sigma$.
	Hence, $\sigma_\alpha$ and $\sigma$ are conjugate and $\ell(\sigma_\alpha) \geq \ell(\sigma)$.
	But the length of $\sigma$ is maximal in its conjugacy class.
	Hence, $\ell(\sigma_\alpha) = \ell(\sigma)$ and \Cref{thm:arrow_and_length} yields $\sigma_\alpha \approx  \sigma$.	
\end{proof}

We begin working towards Statement \ref{enum:parametrization_of_kim_injectivity}.
As before, we will trace the relation $\approx$ back to the elementary steps $\overset{s_i}{\rightarrow}$ with $i\in [n-1]$.
Consider $\sigma\in \SG_n$ and $\tau = s_i \sigma s_i$.
Then we have $\tau \overset{s_i}{\rightarrow} \sigma$ or $\sigma \overset{s_i}{\rightarrow} \tau$ depending on $\ell(s_i \sigma s_i) - \ell(\sigma)$.
Moreover $\sigma \approx \tau$ if and only if the difference vanishes.
Thus our first goal is to determine $\ell(s_i \sigma s_i) - \ell(\sigma)$ depending on $\sigma$ and $s_i$.

\begin{lem}	
	\label{thm:length_and_conjugation_lem}
	Let $\sigma \in \SG_n$ and $i, j\in [n-1]$.
	Then $\set{\sigma(i),\sigma(i+1)} \neq \set{j,j+1}$ if and only if $\left( s_j\in D_L(\sigma) \iff s_j\in D_L(\sigma s_i)\right)$.
\end{lem}

\begin{proof}
	We consider all permutations in one-line notation.  
	From \Cref{eq:descent_sets_of_Sn} it follows for each $\sigma \in \SG_n$  that $j\in D_L(\sigma)$ if and only if $j+1$ is left of~$j$ in $\sigma$.
	
	Now fix $\sigma \in \SG_n$.
	Observe that we obtain $\sigma s_i$ from $\sigma$ by swapping $\sigma(i)$ and $\sigma(i+1)$.
	Since these are two consecutive characters in the the one-line notation of $\sigma$, the relative positioning of~$j$ and $j+1$ is affected by this interchange if and only if $\set{\sigma(i),\sigma(i+1)} = \set{j,j+1}$.
	Now use the note on left descents from the beginning to deduce the claim.
\end{proof}

\begin{lem}
	\label{thm:length_and_conjugation}
	Let $\sigma \in \SG_n$ and $i\in [n-1]$.
	\begin{enumerate}
		\item If $\set{\sigma(i),\sigma(i+1)} \neq \set{i,i+1}$ then
		\begin{align*}
		\ell(s_i \sigma s_i) = 
		\begin{cases}
		\ell(\sigma) -2 &\myif \sigma(i) > \sigma(i+1) \text{ and } \sigma\inv(i) > \sigma\inv(i+1), \\
		\ell(\sigma) +2 &\myif \sigma(i) < \sigma(i+1) \text{ and } \sigma\inv(i) < \sigma\inv(i+1), \\
		\ell(\sigma) &\text{else.}
		\end{cases}
		\end{align*}
		\item If $\set{\sigma(i),\sigma(i+1)} = \set{i,i+1}$ then $i$ and $i+1$ either~are fixpoints of $\sigma$ or form a~$2$-cycle in $\sigma$.
		
	In particular,
		$s_i \sigma s_i = \sigma$. 
	\end{enumerate}
\end{lem}

\begin{proof} 	
	Part~$(2)$ should be clear.
	For Part~$(1)$ assume that $\set{\sigma(i),\sigma(i+1)} \neq \set{i,i+1}$. 
	We have that
	\begin{align*}
	\ell(s_i \sigma s_i) - \ell(\sigma)= \ell(s_i \sigma s_i) - \ell(\sigma s_i) + \ell(\sigma s_i) - \ell(\sigma).
	\end{align*}
	Equation~\Cref{eq:descent_sets_W} yields that each of the two differences on the right hand side is $-1$ or~$1$ depending the truth value of the statements $s_i \in D_L(\sigma s_i)$ and $s_i \in D_R(\sigma)$, respectively.
	From \Cref{thm:length_and_conjugation_lem} we have that $s_i\in D_L(\sigma s_i)$ if and only if $s_i\in D_L(\sigma)$. 
	That is, the first difference depends on whether $s_i \in D_L(\sigma)$ or not.
	Thus, \Cref{eq:descent_sets_of_Sn} implies the claim.
\end{proof}

We now show for each $\alpha\vDash_e n$ that all elements of $\Sigma_\alpha$ have the same orbits of even length on $[n]$.

\begin{lem}
	\label{thm:orbits_of_maximal_classes}
	Let $\alpha \vDash_e n$ and $\sigma \in \SG_n$ such that $\sigma_\alpha \approx \sigma$. Then we have the following. 
	\begin{enumerate}
		\item The orbits of even length of $\sigma$ and $\sigma_\alpha$ on $[n]$ coincide.
		\item Let $\mc O$ be an $\sigma$-orbit on $[n]$ of even length.
		Then the orbits of $\sigma^2$ and $\sigma_\alpha^2$ on $\mc O$ coincide.
	\end{enumerate}
\end{lem}

\begin{proof}
	Since $\sigma_\alpha \approx \sigma$, we have $\sigma_\alpha \rightarrow \sigma$ and $\ell(\sigma_\alpha) = \ell(\sigma)$.
	Using induction on the minimal number of elementary steps $w \overset{s}\rightarrow w'$ (with some $w,w'\in \SG_n$ and $s \in S$) necessary to relate $\sigma_\alpha$ to $\sigma$,
	we may assume that there are $\tau \in \SG_n$ and $s_i\in S$ such that 
	$\sigma_\alpha \rightarrow  \tau \overset{s_i}\rightarrow \sigma$ and $\tau$ satisfies~$(1)$ and~$(2)$ ($\sigma_\alpha$ certainly does).
	Then $\ell(\sigma_\alpha) \geq \ell(\tau) \geq \ell(\sigma)$ so that in fact 
	$\ell(\sigma_\alpha) = \ell(\tau) = \ell(\sigma)$ and $\sigma_\alpha \approx \tau \approx \sigma$ by \Cref{thm:arrow_and_length}.

	It remains to show that $\overset{s_i}\rightarrow$ transfers Properties~$(1)$ and~$(2)$ from $\tau$ to $\sigma$.
	Because $\sigma = s_i \tau s_i$, we obtain $\sigma$ from $\tau$ by interchanging~$i$ and $i+1$ in the cycle notation of $\tau$. If~$i$ and $i+1$ both appear in orbits of uneven length of $\tau$ then~$(1)$ and~$(2)$ are not affected by this interchange.
	Thus, we are left with two cases.
	
	\textbf{Case 1.} 
	Assume that~$i$ and $i+1$ appear in different orbits of $\tau$, say $\mc O_1$ and $\mc O_2$ such that at least one of them, say $\mc O_1$, has even length. We show that this case does not occur. 
	To do this, let $m_1$ and $m_2$ be the minimal elements of $\mc O_1$ and $\mc O_2$, respectively. 
	If $\mc O_2$ also has even length, we assume $m_1 < m_2$. 
	
	For $w\in \SG_n$ and $j\in [n]$ let $\langle w \rangle$ denote the subgroup of $\SG_n$ generated by~$w$ and $\langle w \rangle j$ be the orbit of~$j$ under the natural action of $\langle w \rangle$ on $[n]$.
	Since $\tau$ satisfies Property~$(2)$ and $\mc O_1$ has even length, there is a $p_1 \geq m_1$ such that 
	\begin{align}
	\label{eq:orbits_of_tau}
	\begin{aligned}
	\mc O_1^< 	&	:= \langle \tau^2 \rangle m_1 
	= \langle \sigma_\alpha^2 \rangle m_1 
	= \left\{m_1, m_1+1,\dots, p_1 \right\},\\
	\mc O_1^>& 	:= \langle \tau^2 \rangle \tau(m_1)
	= \langle \sigma_\alpha^2 \rangle \sigma_\alpha(m_1) 
	= \left\{n-m_1+1, n-m_1,\dots, n-p_1+1 \right\}.
	\end{aligned}
	\end{align}

	\begin{cla*}
		Let $a \in \mc O_1^<, b \in \mc O_2$ and $c\in \mc O_1^>$. Then $a<b<c$. 
	\end{cla*}
	
	To prove the claim, consider the positions of elements of $[n]$ in the cycle notation $\sigma_\alpha =\sigma_{\alpha_1} \cdots \sigma_{\alpha_l}$ given by \Cref{thm:element_in_stair_form}.
	The elements on odd positions $1, 2, 3,\dots$ form an strictly increasing sequence.
	The elements on even positions $n, n-1, \dots $ form an strictly decreasing sequence but they are always greater than the entries on odd positions.
	
	We want to show that the elements of $\mc O_2$ all appear right of the cycle consisting of the elements of $\mc O_1$. 
	If $\mc O_2$ has even length this is clear.
	If $\mc O_2$ has odd length, we can use that by Property~$(1)$, the unions of odd orbits of $\tau$ and $\sigma_\alpha$ coincide and that in $\sigma_\alpha$ the elements of odd orbits are all located right of the elements of the even orbits. 
	
	Let $a\in \mc O^<_1$. Then~$a$ is on an odd position and thus it is smaller than any entry right of it.
	On the other hand, $c\in \mc O^>_1$ implies that~$c$ is on an even position and thus is greater then any entry right of it.
	Finally, in the last paragraph we have shown that each $b \in \mc O_2$ is located right of~$a$ and~$c$. 
	This establishes the claim.
	
	Now, we have to deal with two cases. 
	
	If $i\in \mc O_1$ and $i+1 \in \mc O_2$ then the claim implies $i\in \mc O^<_1$.
	Then $\tau^{-1}(i),\tau(i) \in \mc O^>_1$.
	Since $\tau^{-1}(i+1),\tau(i+1) \in \mc O_2$, our claim yields $\tau^{-1}(i) >  \tau^{-1}(i+1)$ and $\tau(i) > \tau(i+1)$. In addition, since $\mc O_1$ has even length and $i+1\not \in \mc O_1$, $\tau(i) \neq i,i+1$. Thus, we obtain from  \Cref{thm:length_and_conjugation} that $\ell(\sigma) < \ell(\tau)$, a contradiction to $\ell(\tau) = \ell(\sigma)$.
	
	If $i+1\in \mc O_1$ and $i \in \mc O_2$ then the claim implies $i+1\in \mc O^>_1$ and similarly as before we obtain $\tau^{-1}(i) >  \tau^{-1}(i+1)$ and $\tau(i) > \tau(i+1)$ and thus the same contradiction using \Cref{thm:length_and_conjugation}.
	That is, we have shown that~$i$ and $i+1$ cannot appear in two different orbits if one of the latter has even length.

	\textbf{Case 2}. Assume that~$i$ and $i+1$ appear in the same orbit with even length $\mc O_1$ of $\tau$. Then~$(1)$ also holds for $\sigma$. 
	
	To show~$(2)$, assume $i+1 \in \langle \tau^2 \rangle i$ first. Then both elements appear in the same cycle of $\tau^2$. As we obtain $\sigma^2$ from $\tau^2$ by swapping~$i$ and $i+1$ in cycle notation,~$(2)$ also holds for $\sigma$. 
	
	Lastly, we show that $i+1  \in \langle \tau^2 \rangle i$ is always true.
	For the sake of contradiction, assume $i+1 \not \in \langle \tau^2 \rangle i$.
	
	Suppose in addition that $\card{\mc O_1} = 2$.
	Then $\set{\tau(i),\tau(i+1)} = \set{i,i+1}$ and from \Cref{thm:length_and_conjugation} we obtain $\sigma = s_i\tau s_i = \tau$.
	This contradicts the minimality of the sequence of arrow relations from $\sigma_\alpha$ to $\sigma$.
	
	Now suppose $\card{\mc O_1} > 2$.
	Then $\set{\tau(i),\tau(i+1)} \neq \set{i,i+1}$.
	Since $i+1 \not \in \langle \tau^2 \rangle i$, it follows from \Cref{eq:orbits_of_tau} that $i = \max \mc O^<_1$ and $i+1= \min \mc O^>_1$.
	Consequently,  $\tau^{-1}(i), \tau(i) \in \mc O_1^>$ and $\tau^{-1}(i+1), \tau(i+1) \in \mc O_1^<$.
	But this means that
	\begin{align*}
	\tau^{-1}(i)> \tau^{-1}(i+1) \text{ and  }\tau(i)> \tau(i+1). 
	\end{align*}
	Because $\set{\tau(i),\tau(i+1)} \neq \set{i,i+1}$, we can now apply \Cref{thm:length_and_conjugation} and obtain that $\ell(\sigma) < \ell(\tau)$.
	Again, we end up with a contradiction.
\end{proof}

Let $\sigma\in \SG_n$. Then the set of orbits of $\sigma$ on $[n]$ is a set partition of $[n]$. We denote this partition by $\partition(\sigma)$.
The set of even orbits of $\sigma$ is given by
\begin{align*}
\partition_e(\sigma) := \set{\mc O \in P(\sigma) \mid \text{$|\mc O|$ is even}}
\end{align*}

If $\partition(\sigma) = \partition(\sigma')$ for $\sigma,\sigma'\in \SG_n$ then $\sigma$ and $\sigma'$ have the same type, \ie they are conjugate.

\begin{lem}
	\label{thm:element_in_stair_form_same_even_orbits_implies_equality}
	Let $\alpha,\beta \vDash_e n$ such that $\sigma_\alpha$ and $\sigma_\beta$ are conjugate. If $\partition_e(\sigma_\alpha) = \partition_e(\sigma_\beta)$ then $\alpha = \beta$. 
\end{lem}

\begin{proof}
	Let $\alpha=\parts{\alpha}l, \beta=\parts{\beta}{l'} \vDash_e n$ and $(x_1, x_2, \dots, x_n)$ be the sequence with $x_{2i-1} = i$ and $x_{2i} = n-i+1$.
	Since $\alpha$ is maximal, there is a $k\in [0,l]$ such that $\alpha_i$ is even for $i\leq k$ and odd for $i>k$.
	Assume that  $\sigma_\alpha$ and $\sigma_\beta$ are conjugate and $\partition_e(\sigma_\alpha) = \partition_e(\sigma_\beta)$.

	Because $\sigma_\alpha$ and $\sigma_\beta$ are conjugate, $\alpha$ and $\beta$ have the same multiset of parts.
	In particular, $l = l'$.
	Since $\alpha$ and $\beta$ are maximal, the odd parts of $\alpha$ and $\beta$ form an weakly decreasing sequence at the end of $\alpha$ and $\beta$, respectively.
	As both compositions have the same length and multiset of parts, it follows that $\alpha_i = \beta_i$ for $i = k+1, \dots, l$. 
	
	We show that $\alpha_i = \beta_i$ for $i=1,\dots, k$ with induction.
	Assume that $i\in [k]$ and $\alpha_j = \beta_j$ for all $1\leq j<i$. Define  $d := \sum_{j=1}^{i-1} \alpha_i$. Then by assumption $d = \sum_{j=1}^{i-1} \beta_i$.
	Moreover, let $\mc O_{\alpha_{i}}$ and $\mc O_{\beta_{i}}$ be the orbits of $x_{d+1}$ under $\sigma_\alpha$ and $\sigma_\beta$, respectively.
	From the definition of elements in stair form it follows that 
	\begin{align*}
	\mc O_{\alpha_{i}} = \set{ x_{d+1}, x_{d+2} , \dots, x_{d+\alpha_i}}, \\
	\mc O_{\beta_{i}} = \set{ x_{d+1}, x_{d+2} , \dots, x_{d+\beta_i}}.
	\end{align*}
	In particular $|\mc O_{\alpha_{i}}| = \alpha_i$ and $|\mc O_{\beta_{i}}| = \beta_i.$
	Since $i\leq k$, $\alpha_i$ and $\beta_i$ are even.
	Consequently,  $\mc O_{\alpha_{i}}$ and $\mc O_{\beta_{i}}$ both have even length.
	Moreover, they have the element $x_{d+1}$ in common. 
	Hence, $\partition_e(\sigma_\alpha) = \partition_e(\sigma_\beta)$ implies $\mc O_{\alpha_{i}} = \mc O_{\beta_{i}}$.
	Thus, $\alpha_i = |\mc O_{\alpha_{i}}|= |\mc O_{\beta_{i}}| = \beta_i$.
\end{proof}

We are now in the position to prove  Statement~\ref{enum:parametrization_of_kim_injectivity}.
This finishes the proof of \Cref{thm:parametrizations_of_kim}.

\begin{proof}[Proof of Statement~\ref{enum:parametrization_of_kim_injectivity}]
	Let $\alpha,\beta \vDash_e n$ such that  $\sigma_\alpha \approx \sigma_\beta$. Then $\sigma_\alpha$ and $\sigma_\beta$ are conjugate. Moreover, \Cref{thm:orbits_of_maximal_classes} implies $\partition_e(\sigma_\alpha) = \partition_e(\sigma_\beta)$. Hence $\alpha = \beta$ by \Cref{thm:element_in_stair_form_same_even_orbits_implies_equality}.
\end{proof}

We use some of the intermediary results that lead to \Cref{thm:parametrizations_of_kim} in order to prepare a result for later use in \Cref{sec:equivalence_classes:inductive_product}.

\begin{prop}
	\label{thm:characterization_of_Sigma_using_length}
	Let $\alpha \vDash_e n$ and $\sigma \in \SG_n$. Then $\sigma \in \Sigma_\alpha$ if and only if
	\begin{enumerate}
		\item $\sigma$ and $\sigma_\alpha$ are conjugate in $\SG_n$,
		\item $\ell(\sigma) = \ell(\sigma_\alpha)$,
		\item $\partition_e(\sigma) = \partition_e(\sigma_\alpha)$.
	\end{enumerate}
\end{prop}
\begin{proof}
	First, assume $\sigma \in \Sigma_\alpha$.
	Because $\sigma_\alpha \in \Sigma_\alpha$ and $\Sigma_\alpha\in \smaxa$, $\sigma$ satisfies $(1)$ and $(2)$. By \Cref{thm:orbits_of_maximal_classes}, $(3)$ holds as well. 
	
	Second, assume that $\sigma$ satisfies $(1) - (3)$. By $(1)$, $\sigma$ is in the same conjugacy class as $\sigma_\alpha$. From $(2)$ it follows, that $\sigma$ is maximal in its conjugacy class. 
	Then \Cref{thm:parametrizations_of_kim} provides the existence of a $\beta\vDash_e n$ such that $\sigma \in \Sigma_\beta$. 
	Using the already proven implication from left to right, we obtain that $\sigma$ and $\sigma_\beta$ are conjugate and  $\partition_e(\sigma) = \partition_e(\sigma_\beta)$.
	But as $\sigma$ satisfies $(1)$ and $(3)$, it follows that $\sigma_{\beta}$ and $\sigma_\alpha$ are conjugate and $\partition_e(\sigma_\beta) = \partition_e(\sigma_\alpha)$. 
	Thus, \Cref{thm:element_in_stair_form_same_even_orbits_implies_equality} yields $\beta = \alpha$ as desired.
\end{proof}

We end this \namecref{sec:parametrization} with a remark on conjugacy classes.

\begin{rem}
The conjugacy classes of $\SG_n$ are parametrized by the partitions of $n$ via the cycle type.
For a composition $\alpha$ we denote the partition obtained by sorting the parts of $\alpha$ in decreasing order by $\widetilde \alpha$.
Let $\lambda \vdash n$ and $\mc O$ be the conjugacy class whose elements have cycle type $\lambda$.
From \Cref{thm:element_in_stair_form} it follows that for $\alpha \vDash_e n$
the element in stair form $\sigma_\alpha$ is contained in $\mc O$ if and only if $\widetilde \alpha = \lambda$.
Hence, \Cref{thm:parametrizations_of_kim} implies that $\set{\sigma_\alpha \mid \alpha \vDash_e n, \widetilde \alpha = \lambda}$ is a complete set of representatives of 
$\faktid{\mc O}{\max}$.
In particular, we have that
\begin{center}
$\left|\faktid{\mc O}{\max}\right| = 1$ if and only if the even parts of $\lambda$ are all equal.
\end{center}
\end{rem}

\section{\texorpdfstring{Equivalence classes of $(\SG_n)_{\max}$ under $\approx$}{Equivalence classes}}

\label{sec:equivalence_classes}

Recall that for $\alpha \vDash_e n$,  $\Sigma_\alpha$ is the $\approx$-equivalence class of the element in stair form $\sigma_\alpha$.
From \Cref{thm:parametrizations_of_kim} we have that $\smaxa = \set{ \Sigma_\alpha \mid \alpha \vDash_e n}$.
In \Cref{thm:basis_of_center_from_elements_in_stair_form} we concluded that the elements $\ngen_{\leq \Sigma_\alpha}$ for $\alpha \vDash_e n$ form a basis of $Z(\He)$.
We emphasize that  $\ngen_{\leq \Sigma_\alpha}$ directly depends on $\Sigma_\alpha$ since $\ngen_{\leq \Sigma_\alpha} = \sum_x \ngen_x$ where $x$ runs over the order ideal in Bruhat order generated by $\Sigma_\alpha$.
This \namecref{sec:equivalence_classes} is devoted to the description of equivalence classes $\Sigma_\alpha$ and bijections between them.

In \Cref{sec:equivalence_classes:one_part} we consider the case where $\alpha$ has only one part. The first result is the characterization of the elements of $\Sigma_{(n)}$ by properties of their cycle notation.
From this we obtain bijections relating $\Sigma_{(n-1)}$ with $\Sigma_{(n)}$ for $n\geq 4$ and a closed formula for the cardinality of $\Sigma_{(n)}$.

In \Cref{sec:equivalence_classes:odd_hook} we generalize the characterization of $\Sigma_{(n)}$ to odd hooks, where a hook $\alpha := (k,1^{n-k})$ is called \emph{odd} if~$k$ is odd and \emph{even} otherwise.
Moreover, we define a bijection $\Sigma_{(k)} \times  [m+1, n-m] \to \Sigma_{(k,1^{n-k})}$ where $m := \frac{k-1}2$.
From this we obtain the cardinality of $\Sigma_{(k,1^{n-k})}$.

In \Cref{sec:equivalence_classes:inductive_product} we consider the inductive product $\iprod$ that allows the decomposition $\Sigma_{(\alpha_1,\dots, \alpha_l)} =\Sigma_{(\alpha_1)} \iprod \Sigma_{(\alpha_2,\dots,\alpha_l)}$ if $\alpha_1$ is even.
Using the results of the previous \namecref{sec:equivalence_classes:inductive_product}s, we infer a description of $\Sigma_\alpha$ for all $\alpha \vDash_e n$ whose odd parts form a hook. 
In particular, we obtain a characterization of $\Sigma_\alpha$ in the case where $\alpha$ is an even hook.

\subsection{\texorpdfstring{Equivalence classes of~$n$-cycles}{Equivalence classes of n-cycles}}
\label{sec:equivalence_classes:one_part}

In this \namecref{sec:equivalence_classes:one_part} we seek a combinatorial description of the elements of $\Sigma_{(n)}$.
Examples are given in \Cref{tbl:Sigma_(n)}.
\begin{table}
	\caption{The elements of $\Sigma_{(n)}$ for small~$n$ with the element in stair form $\sigma_{(n)}$ in the top row.\\}
	\label{tbl:Sigma_(n)}
	$
		\begin{array}{c|cccccc}
		\alpha &	(1) & (2) & (3) & (4) & (5) & (6) \\ \hline
		\multirow{6}{*}{$\Sigma_\alpha$} 
	 	& (1) & (1,2) & (1,3,2) & (1,4,2,3) 
	 			   & (1,5,2,4,3)& (1,6,2,5,3,4) \\
		& 		&     & (1,2,3) & (1,3,2,4) 
				   & (1,5,2,3,4) & (1,6,2,4,3,5) \\
	   	& &  &  &  & (1,5,3,2,4) & (1,6,3,4,2,5) \\
		& &  &  &  & (1,4,2,3,5) & (1,5,2,4,3,6) \\
		& &  &  &  & (1,4,3,2,5) & (1,5,3,4,2,6) \\ 
		& &  &  &  & (1,3,4,2,5) & (1,4,3,5,2,6)\\
		\end{array}
	$	
\end{table}
The description is given by two properties: \emph{being oscillating} and \emph{having connected intervals}.
We begin with the property of being oscillating.

\begin{defi}
	\label{def:oscillating_n-cycle}
	 We call the $n$-cycle $\sigma \in \SG_n$ \emph{oscillating} if there exists a positive integer $m \in \set{\frac{n-1}2, \frac{n}2, \frac{n+1}2}$ such that $\sigma([m]) = [n-m+1, n].$
\end{defi}

In \Cref{thm:oscillating_n-cycle_cycle_notation} we will obtain a more descriptive characterization of oscillating $n$-cycles.
It turns out that the $n$-cycle $\sigma$ of $\SG_n$ (represented in cycle notation) is oscillating if $n$ is even and the entries of $\sigma$ alternate between the sets $[1, \frac{n}2]$ and $[\frac{n}2 + 1, n]$ or $n$ is odd and after deleting the entry $\frac{n+1}2$ from $\sigma$ the remaining entries alternate between the sets $[\frac{n-1}2]$ and $[\frac{n+3}2,n]$.

\begin{samepage}
\begin{exa}
	\begin{enumerate}
	\item
	Recall that for $n\in \N$ the element in stair form $\sigma_{(n)}$ is an $n$-cycle of $\SG_n$.
	For
	\begin{align*}
		\text{
		$\sigma_{(5)} = (1,5,2,4,3)$, \quad
		$\sigma_{(5)}\inv = (1,3,4,2,5)$ \quad
		and \quad
		$\sigma_{(6)} = (1,6,2,5,3,4)$
		}
	\end{align*}
	we have
	\begin{align*}
		\text{
		$\sigma_{(5)}([2]) = [4,5]$, \quad
		$\sigma_{(5)}\inv([3]) = [3,5]$ \quad
		and \quad
		$\sigma_{(6)}([3]) = [4,6]$.
		}
	\end{align*}
	Hence, they are oscillating and the integer $m$ used in \Cref{def:oscillating_n-cycle} is given by
	\begin{align*}
	\text{
		$m = 2 = \frac{5-1}2$, \quad
		$m = 3 = \frac{5+1}2$ \quad
		and \quad
		$m = 3 = \frac{6}2$,
	}
	\end{align*}
	respectively. 
	Note that the entries in the cycles alternate as described after \Cref{def:oscillating_n-cycle}.
	
	\item
	All the elements shown in \Cref{tbl:Sigma_(n)} are oscillating.
	\end{enumerate}
\end{exa}
\end{samepage}

We explicitly write down the three cases for $m$ in \Cref{def:oscillating_n-cycle}.

\begin{rem}
Let $\sigma$ be an oscillating $n$-cycle $\sigma \in \SG_n$ with parameter $m$ from \Cref{def:oscillating_n-cycle}.
Then we have
\begin{enumerate}
\item 
$n$ is even and $\sigma([\frac{n}2]) = [\frac{n}2 +  1, n]$ if $m = \frac{n}2$,
\label{enum:def_oscillating_n-cycle_even}
 	   	\item 
 	   	$n$ is odd and $\sigma([\frac{n-1}2]) = [\frac{n+3}2, n]$ if $m = \frac{n-1}2$,
 		\label{enum:def_oscillating_n-cycle_odd_sigma(n+1/2)<n+1/2}
  	\item $n$ is odd and $\sigma([\frac{n+1}2]) = [\frac{n+1}2, n]$ if $m = \frac{n+1}2$.
	\label{enum:def_oscillating_n-cycle_odd_sigma(n+1/2)>n+1/2}
\end{enumerate}
\end{rem}

Our next aim is to give a characterization of the term \emph{oscillating} in \Cref{thm:oscillating_n-cycle_local}.
By considering complements in $[n]$ we obtain the following.

\begin{lem}
	\label{thm:oscillating_n_cycle_complements}
	Let $\sigma \in \SG_n$ be an $n$-cycle and $m \in [n]$.
	Then $\sigma([m]) = [n-m+1, n]$  if and only if  $\sigma([m+1, n]) = [n-m]$.
\end{lem}

\Cref{thm:oscillating_n_cycle_complements} implies that an $n$-cycle $\sigma\in \SG_n$ is oscillating with parameter $m$ if and only if $\sigma([m+1, n]) = [n-m]$.

\begin{lem}
	\label{thm:oscilating_n-cycle_and_inversion}
	Let $\sigma\in \SG_n$ be an $n$-cycle.
	Then $\sigma$ is oscillating if and only if $\sigma\inv$ is oscillating.
\end{lem}

\begin{proof}
	Let $M := \N \cap  \set{\frac{n-1}2, \frac{n}2, \frac{n+1}2}$.
	If $n = 1$ then $\sigma = \id = \sigma\inv$ (which is oscillating).
	Thus assume $n\geq 2$.
	It suffices to show the implication from left to right.
	Suppose that $\sigma$ is oscillating.
	Then there is an $m \in M$ such that $\sigma([m]) = [n-m+1, n].$
	Consequently, $\sigma([m+1, n]) = [n-m]$ by \Cref{thm:oscillating_n_cycle_complements} and hence
	\begin{align*}
		\sigma\inv([n-m]) = [m+1, n].
	\end{align*}
	 Moreover,  $m+1 = n-(n-m)+1$ and we have $n-m \in M$ since $m \in M$ and $n\geq 2$.
	 Therefore, $\sigma\inv$ is oscillating.
\end{proof}

In the following we rephrase \Cref{def:oscillating_n-cycle} from a more local point of view.

\begin{lem}
	\label{thm:oscillating_n-cycle_local}
	Let $\sigma\in \SG_n$ be an $n$-cycle.
	We consider the four implications for all $i \in [n]$
		\begin{enumerate}
			\item[]
			\begin{enumerate}
			\item
			$i < \frac {n+1}2 \implies \sigma(i) \geq \frac {n+1}2$,
			\label{enum:oscillating_n-cycle_characterization_implication_i<}
			\item
			$i < \frac {n+1}2 \implies \sigma^{-1}(i)\geq \frac {n+1}2$,
			\label{enum:oscillating_n-cycle_characterization_implication_i<_inv}
			\item
			$i > \frac {n+1}2 \implies \sigma(i) \leq \frac {n+1}2$,
			\label{enum:oscillating_n-cycle_characterization_implication_i>}
			\item
			$i > \frac {n+1}2 \implies \sigma^{-1}(i) \leq\frac {n+1}2 $,
			\label{enum:oscillating_n-cycle_characterization_implication_i>_inv}
			\end{enumerate} 
		\end{enumerate}
	and if $n$ is odd the statement
	\begin{enumerate}
		\item[]
	\begin{enumerate}[label = \textup{(}\Alph*\textup{)}]
		\item either $\sigma^{-1}(\frac{n+1}2) >  \frac{n+1}2$ or $\sigma(\frac{n+1}2) >   \frac{n+1}2$.
		\label{enum:oscillating_n-cycle_characterization_xor_contidion}
	\end{enumerate}
	\end{enumerate}
	Then the following are equivalent.
	\begin{enumerate}
		\item $\sigma$ is oscillating.
		\item  One of \ref{enum:oscillating_n-cycle_characterization_implication_i<} -- \ref{enum:oscillating_n-cycle_characterization_implication_i>_inv} is true and if $n$ is odd and $n\geq 3$ then also \ref{enum:oscillating_n-cycle_characterization_xor_contidion} is true.
		\item  Each one of \ref{enum:oscillating_n-cycle_characterization_implication_i<} -- \ref{enum:oscillating_n-cycle_characterization_implication_i>_inv} is true and if $n$ is odd and $n\geq 3$ then also \ref{enum:oscillating_n-cycle_characterization_xor_contidion} is true.
	\end{enumerate}
\end{lem}

\begin{proof}
	
	First suppose that $n$ is odd.
	If $n=1$ then $\sigma= \id$ is oscillating and the implications \ref{enum:oscillating_n-cycle_characterization_implication_i<} -- \ref{enum:oscillating_n-cycle_characterization_implication_i>_inv} are trivially  satisfied.
	
	Assume $n\geq 3$.
	We show for each of the implications (x) that \ref{enum:oscillating_n-cycle_characterization_xor_contidion} and (x) is true if and only if $\sigma$ is oscillating.
	As $n$ is odd and $n\geq 3$,  Statement~\ref{enum:oscillating_n-cycle_characterization_xor_contidion} can be expanded as
	\begin{center}
	\begin{tabular}{cl}
	either 
	&$\sigma\inv(\frac{n+1}2) > \frac{n+1}2$ and $\sigma(\frac{n+1}2) < \frac{n+1}2$	
	\\
	or 
	&$\sigma\inv(\frac{n+1}2) < \frac{n+1}2$ and $\sigma(\frac{n+1}2) > \frac{n+1}2$ .
	\end{tabular}
	\end{center}
	Moreover, \ref{enum:oscillating_n-cycle_characterization_implication_i<} can be rephrased as $\sigma([\frac{n-1}2]) \subseteq [\frac{n+1}2,n]$.
	Hence, we have \ref{enum:oscillating_n-cycle_characterization_xor_contidion} and \ref{enum:oscillating_n-cycle_characterization_implication_i<} if and only if
	\begin{center}
	\begin{tabular}{ccl}
		either 
		& $\sigma([\frac{n-1}2]) = [\frac{n+3}2, n]$
		&(if $\sigma\inv(\frac{n+1}2) > \frac{n+1}2$ and $\sigma(\frac{n+1}2) < \frac{n+1}2$)
		\\
		or 
		&$\sigma([\frac{n+1}2]) = [\frac{n+1}2, n]$ 
		&(if $\sigma\inv(\frac{n+1}2) < \frac{n+1}2$ and $\sigma(\frac{n+1}2) > \frac{n+1}2$).
		\end{tabular}
	\end{center}
	In other words, $\sigma([m]) = [n-m+1, n]$ for either $m = \frac{n-1}2$ or $m =\frac{n+1}2$, \ie $\sigma$ is oscillating.
	
	Similarly, we have \ref{enum:oscillating_n-cycle_characterization_xor_contidion} and \ref{enum:oscillating_n-cycle_characterization_implication_i>} if and only if
		\begin{center}
			\begin{tabular}{llll}
			either 
			&$\sigma([\frac{n+1}2,n]) = [\frac{n+1}2]$
			& or 
			& $\sigma([\frac{n+3}2,n]) = [\frac{n-1}2]$.
			\end{tabular}
		\end{center}
	That is, $\sigma([m+1,n]) = [n-m]$ for either $m = \frac{n-1}2$ or $m = \frac{n+1}2$.
	This is equivalent to $\sigma$ being oscillating by \Cref{thm:oscillating_n_cycle_complements}.
	
	So far we have shown that
	\begin{align}
		\text{
		\ref{enum:oscillating_n-cycle_characterization_xor_contidion} and \ref{enum:oscillating_n-cycle_characterization_implication_i<} $\iff$ 
		$\sigma$ is oscillating
		$\iff$
		\ref{enum:oscillating_n-cycle_characterization_xor_contidion} and \ref{enum:oscillating_n-cycle_characterization_implication_i>}.
		}
		\label{eq:oscilating_n-cycle_and_implications}
	\end{align}
	By \Cref{thm:oscilating_n-cycle_and_inversion} we therefore also have
		\begin{align}
		\text{
			\ref{enum:oscillating_n-cycle_characterization_xor_contidion} and \ref{enum:oscillating_n-cycle_characterization_implication_i<_inv} $\iff$ 
			$\sigma$ is oscillating
			$\iff$
			\ref{enum:oscillating_n-cycle_characterization_xor_contidion} and \ref{enum:oscillating_n-cycle_characterization_implication_i>_inv}.
		}
		\label{eq:oscilating_n-cycle_and_implications_inv}
		\end{align}
	This finishes the proof for odd $n$.
	
	Suppose now that $n$ is even.
	Note that $\frac{n+1}2 \not \in [n]$ as it is not an integer.
	It is not hard to see that the equivalences from \Cref{eq:oscilating_n-cycle_and_implications} and therefore those from \Cref{eq:oscilating_n-cycle_and_implications_inv} hold if  we drop Statement~\ref{enum:oscillating_n-cycle_characterization_xor_contidion}.
	\end{proof}

	We continue with two consequences of \Cref{thm:oscillating_n-cycle_local}.
	We first infer the description of oscillating $n$-cycles mentioned at the beginning of the \namecref{sec:equivalence_classes:one_part}.

	\begin{cor}
		\label{thm:oscillating_n-cycle_cycle_notation}
		Let $\sigma \in \SG_n$ be an $n$-cycle.
		We consider $\sigma$ in cycle notation.
		Then $\sigma$ is oscillating if and only if one of the following is true.
		\begin{enumerate}
			\item $n$ is even and the entries of $\sigma$ alternate between the sets $\left[\frac{n}2\right]$ and $\left[\frac{n}2+1,n\right]$.
			\item $n$ is odd and after deleting the entry $\frac{n+1}2$ from $\sigma$, the remaining entries alternate between the sets $\left[\frac{n-1}2\right]$ and $\left[\frac{n+3}2, n\right]$.
		\end{enumerate}
	\end{cor}
	
	\begin{proof}
	With \ref{enum:oscillating_n-cycle_characterization_xor_contidion}, \ref{enum:oscillating_n-cycle_characterization_implication_i<} and \ref{enum:oscillating_n-cycle_characterization_implication_i>} we refer to the statements of \Cref{thm:oscillating_n-cycle_local}.
	
	Suppose that $n$ is even. 
	By \Cref{thm:oscillating_n-cycle_local}, $\sigma$ is oscillating if and only if the implications \ref{enum:oscillating_n-cycle_characterization_implication_i<} and \ref{enum:oscillating_n-cycle_characterization_implication_i>} are satisfied which is the case if and only if the entries of $\sigma$ alternate between $[\frac{n}2]$ and $[\frac{n}2+1,n]$.
	
	Suppose that $n$ is odd. 
	If $n\geq 3$ then property \ref{enum:oscillating_n-cycle_characterization_xor_contidion} states that one of the neighbors $\sigma\inv(\frac{n+1}2)$ and $\sigma(\frac{n+1}2)$ of $\frac{n+1}2$ in $\sigma$ is an element of $[\frac{n-1}2]$ and the other one is an element of $[\frac{n+3}2, n]$.
	Therefore, $\sigma$ satisfies \ref{enum:oscillating_n-cycle_characterization_xor_contidion}, \ref{enum:oscillating_n-cycle_characterization_implication_i<} and \ref{enum:oscillating_n-cycle_characterization_implication_i>}
	if and only if
	after deleting $\frac{n+1}2$ from the cycle notation of $\sigma$, the remaining entries alternate between the sets $[\frac{n-1}2]$ and $[\frac{n+3}2, n]$.
	Thus, \Cref{thm:oscillating_n-cycle_local} yields that the latter property is satisfied if and only if $\sigma$ is oscillating.
	\end{proof}

By considering $\sigma$ in cycle notation beginning with $1$, we can rephrase \Cref{thm:oscillating_n-cycle_cycle_notation} in a more formal way.

\begin{cor}
	\label{thm:reformulation_of_oscillating}
	 Let $\sigma \in \SG_n$ be an~$n$-cycle. If~$n$ is odd, let $0\leq l \leq n-1$ be such that $\sigma^l(1) = \frac{n+1}2$. If~$n$ is even, set $l := \infty$. Then $\sigma$ is oscillating if and only if
	 for all $0\leq k \leq n-1$ we have
\begin{align*}
\begin{aligned}
\sigma^k(1) &< \frac{n+1}2 \quad \text{if $k<l$ and~$k$ is even or $k>l$ and~$k$ is odd}, \\
\sigma^k(1) &> \frac{n+1}2 \quad \text{if $k<l$ and~$k$ is odd or $k>l$ and~$k$ is even}.
\end{aligned}
\end{align*}
\end{cor}

We now consider the second property in the characterization of $\Sigma_{(n)}$: the property of \emph{having connected intervals}.
Roughly speaking, an $n$-cycle of $\SG_n$ has connected intervals if in its cycle notation for each $1\leq k\leq \frac{n}2$ the elements of the interval $[k, n-k+1]$ are grouped together.

	\begin{defi}
		\label{def:connected_intervals_n-cycle}
		\begin{enumerate}
			\item
		Let $\sigma \in \SG_n$ and $M\subseteq [n]$. We call~$M$ \emph{connected} in $\sigma$ if there is an $m\in M$ such that 
		\begin{align*}
		M = \set{m, \sigma(m), \sigma^2(m),\dots, \sigma^{|M|-1}(m)}.
		\end{align*}
			\item Let $\sigma \in \SG_n$ be an~$n$-cycle. We say that $\sigma$ has \emph{connected intervals} if the interval ${[k, n-k+1]}$ is connected in $\sigma$ for all integers $k$ with $1\leq k \leq \frac n2$.
	\end{enumerate}
	\end{defi}
	\begin{exa}
		All elements shown in \Cref{tbl:Sigma_(n)} have connected intervals. 
		In particular, the element in stair form $\sigma_{(6)} = (1,6,2,5,3,4)$ has connected intervals. In contrast, in $(1,5,2,6,3,4)$ the set $[2,5]$ is not connected.
	\end{exa}

	 The main result of this \namecref{sec:equivalence_classes:one_part} is that an~$n$-cycle $\sigma \in \SG_n$ is an element of $\Sigma_{(n)}$ if and only if $\sigma$ is oscillating and has connected intervals.
	 We now begin working towards this result.

	\begin{lem}
		\label{thm:sigma_(n)_osc_and_c.I.}
		 The element in stair form $\sigma_{(n)} \in \SG_n$ is oscillating and has connected intervals.
	\end{lem}
	
	\begin{proof}
		By \Cref{thm:element_in_stair_form},
		\begin{align*}
		 \sigma_{(n)} = 
		 \begin{cases}
	  (1, n, 2, n-1, \dots, \frac n2, n - \frac n2 +1) & \text{if~$n$ is even} \\
   (1, n, 2, n-1, \dots, \frac{n-1}2,  n - \frac{n-1}2+1, \frac {n+1}2) & \text{if~$n$ is odd}.
		 \end{cases}
		\end{align*}
		Thus, 	
		$\sigma_{(n)}([\tfrac{n}2])  = [\tfrac{n}2 +1, n]$
		if $n$ is even and 
		$\sigma_{(n)}([\tfrac{n-1}2]) = [\tfrac{n+3}2, n]$ if $n$ is odd.
		That is, $\sigma_{(n)}$ is oscillating.

		For all $k\in \N$ with $1\leq k \leq \frac{n}2$ the rightmost ${|[k,n-k+1]|}$ elements in the cycle of $\sigma_{(n)}$ from above form ${[k,n-k+1]}$. Thus, $\sigma_{(n)}$ has connected intervals.
	\end{proof}

	Let $\sigma \in \SG_n$. 
	Sometimes it will be convenient to consider $\sigma^{w_0}$ instead of $\sigma$.
	We will now show that conjugation with the longest element $w_0$ of $\SG_n$ preserves the properties of being oscillating and having connected intervals.
	
	\begin{lem}
		\label{thm:conjudation_with_w0_oscillating_and_c.I.}
		Let $\sigma\in \SG_n$ be an~$n$-cycle.
		\begin{enumerate}
		\item If $\sigma$ is oscillating then $\sigma^{w_0}$ is oscillating.
		\item If $\sigma$ has connected intervals then $\sigma^{w_0}$ has connected intervals.
		\end{enumerate}
	\end{lem}
	
	\begin{proof}
		If $n = 1$ the result is trivial.
		Thus suppose $n\geq 2$.
		\begin{wideenumerate}
			\item Set $M := \N \cap \set{\frac{n-1}2, \frac{n}2, \frac{n+1}2}$ and assume that $\sigma$ is oscillating.
			Then there is an $m\in M$ such that $\sigma([m]) = [n-m+1,n]$ and from \Cref{thm:oscillating_n_cycle_complements} it follows that $\sigma([m+1,n]) = [n-m]$.
			Using $w_0(i) = n-i+1$ for $i\in [n]$, we obtain
			\begin{align*}
			\sigma^{w_0}([n-m]) 
			&= w_0\sigma w_0([n-m]) \\
			&= w_0\sigma ([m+1,n]) \\
			&= w_0([n-m]) \\
			&= [n-(n-m)+1, n].
			\end{align*}

		As $n-m\in M$, it follows that $\sigma^{w_0}$ is oscillating.
		\item
		Let $I := [k,n-k+1]$ be given by an integer $k$ with $1\leq k \leq \frac{n}2$. Then  $w_0(I)= I$.
		Hence, if~$I$ is connected in $\sigma$ then it is also connected in $\sigma^{w_0}$.
		\qedhere
		\end{wideenumerate}
	\end{proof}

In the following result we study the interplay between the conjugation with $w_0$ and the relation $\approx$.
The generalization to all finite Coxeter groups is straight forward.

	\begin{lem} 
			\label{thm:conjudation_with_w0_and_arrow}
			Let $w,w'\in \SG_n$ 
			and $\na$ be the automorphism of $\SG_n$ given by $x \mapsto x^{w_0}$.
		\begin{enumerate}
			\item If $w\overset{s_i}{\to}w'$ then $\na(w)\overset{s_{n-i}}{\to}\na(w')$.
			\item If $w\approx w'$ then $\na(w) \approx \na(w')$.
		\end{enumerate}	
	\end{lem}

	\begin{proof} 
	Assume $w\overset{s_i}{\to}w'$. Then $w' = s_i w s_i$ and $\ell(w') \leq \ell(w)$.
	Since $\na(s_i) = s_{n-i}$, we have $\na(w') = s_{n-i} \na(w) s_{n-i}$.
	Moreover, $\ell(\na(w')) \leq \ell(\na(w))$ because $\ell(x) = \ell(\na(x))$ for all $x\in \SG_n$. Thus, $\na(w)\overset{s_{n-i}}{\to}\na(w')$. Now, use the definition of $\approx$ to obtain $(2)$ from $(1)$.
	\end{proof}

		Consider $n=5$, the oscillating~$n$-cycle $\sigma = (\bs{1},4,2,3,\bs{5})$ and its connected interval $I = \set{2,3,4}$. 
		In the cycle notation of $\sigma$, this interval is enclosed by the two elements $a = 1$ and $b = 5$.
		Note that $\frac{n+1}2 = 3$, $a<3$ and $b>3$.
		This illustrates a property of oscillating~$n$-cycles which is the subject of the next \namecref{thm:oscillation_around_intervals}.
	
	\begin{lem} 
		\label{thm:oscillation_around_intervals}
		Assume that $\sigma\in \SG_n$ is an oscillating~$n$-cycle with a connected interval $I := [i,n-i+1]$ such that $i\in \N$ and $2\leq i \leq \frac{n+1}2$. Let $r := |I|$ and $m\in I$ be such that $I = \set{\sigma^k(m) \mid k = 0,\dots, r-1}$. Moreover, set $a := \sigma^{-1}(m)$ and $b:=\sigma^r(m)$.
		Then $a,b \neq \frac{n+1}2$ and
		\begin{align*}
		a <\frac{n+1}2 \iff b>\frac{n+1}2.
		\end{align*}
	\end{lem}
	\begin{proof}
		Let $p\in [n-1]$ be such that $\sigma^p(1) = a$. Then $\sigma^{p+r+1}(1) = b$.
		Since $i>1$, $1\not  \in I$ and thus $p+r+1\leq n-1$.
		We have $r = n-2i+2$. Hence,~$r$ has the same parity as~$n$.
		
		We want to apply \Cref{thm:reformulation_of_oscillating}.
		If~$n$ is odd, let $l\in [0,n-1]$ be such that $\sigma^l(1) = \frac{n+1}2$.
		Then $\frac{n+1}2 \in I$ so that $p< l < p + r+1$.
		In particular, $a,b \neq \frac{n+1}2$.
		Clearly, if~$n$ is even then $a,b \neq \frac{n+1}2$.
		
		Therefore,
		\begin{align*}
		a=\sigma^p(1) < \frac{n+1}2 
		&\iff \text{$p$ is even} \\
		&\iff 
		\begin{cases}
			\text{$p+r+1$ is odd} & \text{if~$n$ even} \\
			\text{$p+r+1$ is even} & \text{if~$n$ odd} 
		\end{cases}\\
		&\iff b=\sigma^{p+r+1}(1) > \frac{n+1}2.
		\end{align*}
		where we use \Cref{thm:reformulation_of_oscillating} (and $p<l<p+r+1$ if~$n$ is odd) for the first and third equivalence.
	\end{proof}

Since the $\to$ relation is the transitive closure of the $\overset{s_i}{\to}$ relations, we are interested in the circumstances under which the conjugation with $s_i$ preserves the property of being oscillating with connected intervals.

\begin{lem}
	\label{thm:characterization_osc_and_ci}
	Let $\sigma\in \SG_n$ be an oscillating~$n$-cycle with connected intervals, $i\in [n-1]$ with  $i\leq \frac {n+1}2$ and  $\sigma' := s_i \sigma s_i$. Then $\sigma'$ is oscillating and has connected intervals if and only if
	\begin{enumerate}
		\item if $i  = \frac n2$ then $n=2$,
		\item if $i =\frac{n-1}2$ or $i = \frac{n+1}2$ then $\sigma(i) = i+1$ or $\sigma^{-1}(i) = i+1$,
		\item if $i < \frac{n-1}2$ then
		\begin{align*}
		\text{$\sigma(i) \in I$ and $\sigma(i+1) \not \in I$ or $\sigma\inv(i)\in I$ and $\sigma\inv(i+1) \not \in I$}
		\end{align*} where $I := [i+1, n-i]$. 
	\end{enumerate}
\end{lem}

\begin{proof}
	We will use \Cref{thm:oscillating_n-cycle_local} without further reference.
	Note that $\sigma' = s_i \sigma s_i$ means that we obtain $\sigma'$ from $\sigma$ by interchanging $i$ and $i+1$ in cycle notation.
	We show the equivalence case by case, depending on~$i$.
	\begin{caseenum}
	\item
	 Suppose $i = \frac{n}2$. In this case~$n$ is even. If $n=2$ then $(1,2)$ is the only~$2$-cycle in $\SG_n$. 
	Thus, $\sigma = \sigma' = (1,2)$.
	This element is oscillating and has connected intervals.
	
	Assume now that $n>2$. Since $\sigma$ is oscillating,
	\begin{align*}
	\text{$\sigma(i) > \frac n2$ and $\sigma\inv(i) > \frac n2$.}
	\end{align*}
	Moreover as $n>2$, at most one of  $\sigma(i)$ and $\sigma\inv(i)$ equals $i+1$. Since we obtain $\sigma'$ from $\sigma$ by swapping~$i$ and $i+1$ in cycle notation we infer
	\begin{align*}
	\text{$\sigma'(i+1) > \frac n2$ or ${\sigma'}^{-1}(i+1) > \frac n2$.}
	\end{align*}
	As $i+1 > \frac{n}2$, this means that $\sigma'$ is not oscillating 
	
	\item
	Suppose $i = \frac{n-1}2$ or $i =  \frac{n+1}2$. In this case~$n$ is odd and $n\geq 3$. Moreover, $i,i+1 \in [k,n-k+1]$ for $k= 1, \dots, \frac{n-1}2$.
	Hence, each of the intervals remains connected if we interchange~$i$ and $i+1$.
	Therefore, $\sigma'$ has connected intervals.
	It remains to determine in which cases $\sigma'$ oscillates. 
	We do this for $i= \frac{n-1}2$. 
	The proof for $i = \frac{n+1}2$ is similar.
	
	For $i= \frac{n-1}2$ we have $i+1 = \frac{n+1}2$.
	Since $\sigma$ is oscillating, 
	\begin{align*}
		\text{$\sigma(i) \geq \frac{n+1}2$ and $\sigma^{-1}(i) \geq \frac{n+1}2$}.
	\end{align*}
	Because $n\geq 3$, there is at most one equality among these two inequalities.
	Assume that there is no equality at all. Then
	\begin{align*}
	\text{$\sigma'\left(\frac{n+1}2\right) > \frac{n+1}2$ and ${\sigma'}^{-1}\left(\frac{n+1}2\right) > \frac{n+1}2$}
	\end{align*}
	since  $\sigma' = s_i\sigma s_i$.
	Hence, $\sigma'$ is not oscillating. 
	
	Conversely, assume that $\sigma(i) = i+1$ or $\sigma^{-1}(i) = i+1$.
	In other words, there exists an $\varepsilon \in \set{-1,1}$ such that $\sigma^\varepsilon(i) = i+1$.
	Since $i+1 = \frac{n+1}2$ and $\sigma$ is oscillating, we then have $a:= \sigma^{-\varepsilon}(i) > \frac{n+1}2$.
	Moreover, $\sigma^{-\varepsilon}(i+1) = i < \frac{n+1}2$.
	Thus $\sigma$ being oscillating implies that $b := \sigma^{\varepsilon}(i+1) > \frac{n+1}2$.
	By definition of~$a$ and~$b$,
	\begin{align*}
		\sigma^{\varepsilon} = ( a , i, i+1, b, \dots ).
	\end{align*}
	As a consequence,
	\begin{align*}
	{\sigma'}^{\varepsilon} = (  a , i+1, i, b, \dots )
	\end{align*}
	and $\sigma^{\varepsilon}$ and ${\sigma'}^{\varepsilon}$ coincide on the part represented by the dots because $\sigma' = s_i \sigma s_i$.
	From $a> \frac{n+1}2$,
	$i+1 = \frac{n+1}2$,
	$i < \frac{n+1}2$
	and
	$b > \frac{n+1}2$
	it now follows that $\sigma'$ is oscillating.

\item
Suppose $i< \frac{n-1}2$. Note that then $n\geq 4$.
Define $I := [i+1, n-i]$ as in the theorem and set $r := |I|$. Since $i+1< \frac{n+1}2$, we have $r>1$. 
We show the implication from left to right first.
Assume that $\sigma'$ is oscillating and has connected intervals. 
Note that
\begin{align*}
\tau^\varepsilon(j) \neq i,i+1  \text{ for all } \tau \in \set{\sigma,\sigma'}, \varepsilon\in \set{-1,1} \text{ and } j\in \set{i,i+1}
\end{align*}
since $\sigma$ and $\sigma'$ are oscillating and $i,i+1<\frac{n+1}2$.
Because~$I$ is connected in $\sigma'$, $i+1\in I$ and $r>1$, we have that
\begin{align*}
\exists \varepsilon \in \set{-1,1} \text{ such that }  {\sigma'}^\varepsilon(i+1) \in I.
\end{align*}
Therefore, 
\begin{align*}
\exists \varepsilon \in \set{-1,1} \text{ such that }  \sigma^\varepsilon(i) \in I
\end{align*}
as $\sigma' = s_i\sigma s_i$ and ${\sigma'}^\varepsilon(i+1) \neq i,i+1$.
In fact, the statement
\begin{align}
\label{eq:xor_for_i}
\exists \varepsilon \in \set{-1,1} \text{ such that }  \sigma^\varepsilon(i) \in I \text{ and } \sigma^{-\varepsilon}(i)\not\in I
\end{align}  
is true since otherwise we would have
\begin{align*}
\sigma = (n+i-1 , \dots, \sigma\inv(i), i, \sigma(i) , \dots )
\end{align*}
with $\sigma\inv(i),  \sigma(i) \in I$ and $i, n+i-1\not \in I$ in which case~$I$  would not be connected in $\sigma$.

By interchanging the roles played by $\sigma$ and $\sigma'$ in the argumentation leading to \Cref{eq:xor_for_i}, we get that
\begin{align*}
\exists \varepsilon \in \set{-1,1} \text{ such that }  {\sigma'}^\varepsilon(i) \in I \text{ and } {\sigma'}^{-\varepsilon}(i)\not\in I.
\end{align*}  
From this we obtain that
\begin{align}
\label{eq:xor_for_i+1}
\exists \varepsilon \in \set{-1,1} \text{ such that }  \sigma^\varepsilon(i+1) \in I \text{ and } \sigma^{-\varepsilon}(i+1)\not\in I
\end{align}
by swapping~$i$ and $i+1$ in cycle notation and using that $\sigma'(i),{\sigma'}^{-1}(i) \neq i,i+1$.

Now, let $\varepsilon\in\set{-1,1}$ be such that $\sigma^\varepsilon(i) \in I$ and $\sigma^{-\varepsilon}(i)\not\in I$.
Then
\begin{align}
\label{eq:description_of_I}
I = \set{\sigma^{\varepsilon k}(i) \mid k = 1,\dots, r}
\end{align}
since~$I$ is connected in $\sigma$ and $i\not \in I$.
From \Cref{eq:xor_for_i+1} it follows that $i+1$ appears at the border of~$I$ in the cycle notation of $\sigma$.
Hence,  \Cref{eq:description_of_I} implies that
\begin{align*}
\text{$\sigma^{\varepsilon}(i) = i+1$ or $\sigma^{\varepsilon r}(i) = i+1$.}
\end{align*}
As $\sigma^\varepsilon(i) \neq i+1$, it follows that $i+1 = \sigma^{\varepsilon r}(i)$.
Thus, \Cref{eq:description_of_I} yields that $\sigma^{-\varepsilon}(i+1) \in I$ and $\sigma^{\varepsilon}(i+1)\not \in I$.
Therefore, we have $\sigma^{\varepsilon}(i) \in I$ and $\sigma^{\varepsilon}(i+1)\not \in I$ for an $\varepsilon \in \set{-1,1}$ as desired.

Lastly, we prove the direction from right to left of the equivalence.
We are still in the case $i < \frac{n-1}2$.
Thus, assume that there is an $\varepsilon \in \set{-1,1}$ such that $\sigma^{\varepsilon}(i)\in I$ and $\sigma^{\varepsilon}(i+1) \not \in I$.
Since $\sigma$ is oscillating and we interchange two elements $i,i+1<\frac{n+1}2$ in $\sigma$ in order to obtain $\sigma'$ from $\sigma$, $\sigma'$ is also oscillating. 

It remains to show that $\sigma'$ has connected intervals.
Since  $i\not \in I$, $\sigma^{\varepsilon}(i)\in I$ and~$I$ is connected in $\sigma$, we have \Cref{eq:description_of_I}.
Moreover, from $i+1\in I$, $\sigma^{\varepsilon}(i+1)\not\in I$ and~$I$ being connected in $\sigma$, it follows that $\sigma^{\varepsilon r}(i) = i+1$.
Thus,
\begin{align*}
I = \set{{\sigma'}^{\varepsilon k}(i+1) \mid k = 0,\dots, r-1}
\end{align*}
because $\sigma' = s_i \sigma s_i$.
That is,~$I$ is connected in $\sigma'$.
Let $J := [k, n-k+1]$ for $k\in \N$ with $1\leq k \leq \frac{n}2$ and $k\neq i+1$ be an interval different from~$I$.
Then either $i,i+1 \in J$ or $i,i+1 \not \in J$. As~$J$ is connected in $\sigma$ and $\sigma' = s_i\sigma s_i$, it follows that~$J$ is connected in $\sigma'$. 
Therefore, $\sigma'$ has connected intervals.
\qedhere
	\end{caseenum}
\end{proof}
	
\begin{exa}
	Consider $\sigma = \sigma_{(6)} = (1, 6,2,5,3,4)$ and $\sigma_i := s_i\sigma s_i$ for $i = 1,2$.
	Then $\sigma$ is oscillating with connected intervals.
	
	Since $\sigma^{-1}(1)\in [2,5]$ and $\sigma^{-1}(2) \not \in [2,5]$, 
	\Cref{thm:characterization_osc_and_ci} yields that $\sigma_1$ is oscillating with connected intervals.
	In contrast, $\sigma_2$ is not oscillating with connected intervals because of $\sigma(2),\sigma^{-1}(2)\not \in [3,4]$ and \Cref{thm:characterization_osc_and_ci}.
	This can also be checked directly. We have
	\begin{align*}
		\sigma_1 = (1,5,3,4,2,6) \quad \text{and} \quad \sigma_2 = (1,6,3,5,2,4).
	\end{align*}
	For instance, $[3,4]$ is not connected in $\sigma_2$.
\end{exa}

In the next result we show that the relation $\approx$ is compatible with the concept of oscillating~$n$-cycles with connected intervals.
	
\begin{lem}
	\label{thm:equivalence_implies_oscillating_and_ci}
	Let $\sigma\in \SG_n$ be an oscillating~$n$-cycle with connected intervals, $i \in [n-1]$  and $\sigma' := s_i \sigma s_i$.
	If $\sigma \approx \sigma'$ then $\sigma'$ is oscillating and has connected intervals.
\end{lem}	
\begin{proof} 
	We do a case analysis depending on~$i$.
	
	\textbf{Case 1.} 
	Suppose $i=\frac n2$.
	Then $n$ is even.
	By \Cref{thm:characterization_osc_and_ci}, $\sigma'$ is oscillating with connected intervals if and only if $n=2$.
	Thus, we have to show that $\sigma \not \approx \sigma'$ if $n\geq 4$.
	In this case we have $\sigma(i),\sigma\inv(i) > \frac n2$ and $\sigma(i+1),\sigma\inv(i+1) \leq \frac n2$ because $\sigma$ is oscillating.
	But then \Cref{thm:length_and_conjugation} yields $\ell(\sigma') < \ell(\sigma)$ so that $\sigma'\not \approx \sigma$.

	\textbf{Case 2.} Suppose $i =\frac {n-1}2$ or $i =  \frac {n+1}2$. We only do the case $i =\frac {n-1}2$.
	The other one is similar.
	Let $I := [i, n-i+1] = \set{i,i+1,i+2}$.
	We show the contraposition and assume that $\sigma'$ is not oscillating or that it does not have connected intervals.
	Then from \Cref{thm:characterization_osc_and_ci} it follows that   $\sigma(i) \neq i+1$ and $\sigma^{-1}(i) \neq i+1$.
	Furthermore, there is an $m\in I$ such that
	\begin{align*}
 I = \set{\sigma^{-1}(m), m, \sigma(m)}
	\end{align*}
	since~$I$ is connected in $\sigma$.
	Thus, $m= i+2$.
	Assume $\sigma\inv(i+2) = i$ and $\sigma(i+2) = i$ (the proof of  the other case with $\sigma(i+2) = i$ is analogous). Then $\sigma^{-1}(i) > i+2$ as $\sigma$ is oscillating and $\sigma\inv(i) \neq i+1,i+2$. Moreover, \Cref{thm:oscillation_around_intervals} applied to~$I$ in $\sigma$ and $\sigma^{-1}(i)> \frac{n+1}2$ yields $\sigma(i+1) < \frac{n+1}2 = i+1$. 
	Therefore,
	\begin{align*}
	\sigma(i) = i+2 > \sigma(i+1) \text{\quad and \quad} \sigma\inv(i) > i+2 = \sigma\inv(i+1)
	\end{align*}
	so that $\ell(\sigma') < \ell(\sigma)$ by \Cref{thm:length_and_conjugation} and hence $\sigma'\not\approx \sigma$.
	
	\textbf{Case 3}. Suppose $i<\frac{n-1}2$. Then for all $j\in\set{i,i+1}$ we have $\sigma(j),\sigma\inv(j) \geq \frac{n+1}2$ since $j< \frac{n+1}2$ and $\sigma$ is oscillating. 
	We assume $\sigma\approx \sigma'$ and show that $\sigma'$ is oscillating and has connected intervals. Define $I_k := [k,n-k+1]$ for all $k\leq \frac{n+1}2$ and $I := I_{i+1} = [i+1,n-i]$. 
	Thanks to \Cref{thm:characterization_osc_and_ci} it suffices to show
	\begin{align*}
	\text{$\sigma(i) \in I$ and $\sigma(i+1) \not \in I$ or $\sigma\inv(i)\in I$ and $\sigma\inv(i+1) \not \in I$}.
	\end{align*}
	Since $\sigma \approx \sigma'$, $\ell(\sigma) = \ell(\sigma')$.
	Hence,  \Cref{thm:length_and_conjugation} implies that either $\sigma(i) < \sigma(i+1)$ or $\sigma\inv(i) < \sigma\inv(i+1)$.
	We assume  $\sigma(i) < \sigma(i+1)$ and $\sigma\inv(i) > \sigma\inv(i+1)$. The other case is similar. 
	
	First we show $\sigma(i) \in I$.
	Assume $\sigma(i) \not \in I$ instead.
	Then $\sigma(i) \geq \frac{n+1}2$ implies $\sigma(i) > \max I$. Now we use that $\sigma(i) < \sigma(i+1)$ to obtain $\sigma(i+1)\not \in I$. 
	From this it follows that
	\begin{align*}
	I = \set{\sigma^{-k}(i+1) \mid k = 0,\dots, r-1}
	\end{align*}
	where $r := |I|$ since~$I$ is connected in $\sigma$ and $i+1\in I$. 
	Now we consider the interval $I_i = [i,n-i+1]$ in $\sigma$.
	Because $\sigma$ is oscillating, $\sigma(i+1) > \frac{n+1}2$. 
	An application of \Cref{thm:oscillation_around_intervals} to~$I$ in $\sigma$ yields $\sigma^{-r}(i+1) < \frac{n+1}2$. 
	In particular, $\sigma^{-r}(i+1) \neq n-i+1$. 
	But we also have $i \neq \sigma^{-r}(i+1)$ because $\sigma(i)\not \in I$. That is $\sigma^{-r}(i+1) \not \in I_i$.
	As a consequence,
	\begin{align*}
	I_i = \set{\sigma^{-k}(i+1) \mid k = 0,\dots, r-1} \cup \set{\sigma(i+1), \sigma^2(i+1)} 
	\end{align*}
	since $I \subseteq I_i$ and $I_i$ is connected in $\sigma$.
	Hence
	\begin{align*}
		\set{\sigma(i+1), \sigma^2(i+1)} = \set{i, n-i+1}. 
	\end{align*}
	As $\sigma(i+1) > \frac{n+1}2$, it follows that $\sigma(i+1) = n-i+1$ and $\sigma^2(i+1) = i$. Consequently,  
	\begin{align*}
		\sigma(i) > \max I_i = n-i+1 = \sigma(i+1). 
	\end{align*}
	This is a contradiction to $\sigma(i) < \sigma(i+1)$ and shows that $\sigma(i)\in I$. 
	
	It remains to show that $\sigma(i+1)\not \in I$. 
	Because $i\not \in I$, $\sigma(i) \in I$ and~$I$ is connected,
	\begin{align*}
	 I = \set{\sigma^{k}(i) \mid k = 1,\dots, r}.
	\end{align*}
	We can apply \Cref{thm:oscillation_around_intervals} to~$I$ in $\sigma$ and $i<\frac{n+1}2$ to obtain $\sigma^{r+1}(i) > \frac{n+1}2$.
 Thus $\sigma^{r}(i) \leq \frac{n+1}2$. In particular, $\sigma^r(i) \neq n-i$.
 
	If $i = \frac n2 -1$ then $I =\set{i+1, n-i}$ and it follows that  $\sigma(i) = n-i$ and $\sigma^2(i) = i+1$. That  is, $\sigma(i+1) \not \in I$ as desired.
	
	Now suppose $i< \frac n2 - 1$.
	Then $i+2 \leq \frac{n+1}2$ and we consider $I_{i+2} = [i+2, n-i-1]$. Assume for the sake of contradiction that $\sigma(i+1)\in I$. This means that $\sigma^r(i) \neq i+1$. In addition, we have already seen that $\sigma^r(i) \neq n-i$. Therefore, $\sigma^r(i) \in I_{i+2}$. Since $I_{i+2}$ is connected in $\sigma$ and $I_{i+2} \subseteq I$, we have
	\begin{align*}
	I_{i+2} = \set{\sigma^{k}(i) \mid k = 3,\dots, r}.
	\end{align*}
	and hence $\set{\sigma(i), \sigma^2(i)} = \set{i+1, n-i}$. As $i < \frac{n+1}2$, it follows that $\sigma(i) = n-i$ and $\sigma^2(i) = i+1$.
	But then
	\begin{align*}
	\sigma(i) = n-i > n-i-1 = \max I_{i+2} \geq \sigma(i+1)
	\end{align*}
	which again contradicts the assumption $\sigma(i) < \sigma(i+1)$ and thus shows that $\sigma(i+1) \not \in I$.
	
	\textbf{Case 4.} Suppose $i>\frac {n+1}2$. Assume $\sigma \approx \sigma'$ and let $\na\colon \SG_n \to \SG_n, x\mapsto x^{w_0}$, $\tau := \na(\sigma)$ and $\tau':=\na(\sigma')$.
	Since $\sigma$ is oscillating and has  connected intervals, \Cref{thm:conjudation_with_w0_oscillating_and_c.I.} implies that $\tau$ is oscillating and has connected intervals. 
	In addition, from \Cref{thm:conjudation_with_w0_and_arrow} we have $\tau \approx \tau'$. 
	Because
	$\tau' = s_{n-i}\tau s_{n-i}$ with $n-i < \frac{n+1}2$, we now obtain from the already proven cases that $\tau'$ is oscillating and has connected intervals. Hence, $\sigma' = \na(\tau')$ and \Cref{thm:conjudation_with_w0_oscillating_and_c.I.} yield that $\sigma'$ is oscillating with connected intervals.
\end{proof}

In order to show that each oscillating~$n$-cycle with connected intervals is $\approx$-equivalent to $\sigma_{(n)}$, we use
 an algorithm that takes an oscillating~$n$-cycle $\sigma\in \SG_n$ with connected intervals as input and successively conjugates $\sigma$ with simple reflections until we obtain $\sigma_{(n)}$.
This algorithm has the property that all permutations appearing as interim results are oscillating with connected intervals and $\approx$-equivalent to $\sigma$.
Eventually, it follows that $\sigma \approx \sigma_{(n)}$.

The mechanism of the algorithm is due to Kim \cite{Kim1998}.
She used it in order to show that for each $\alpha\vDash_e n$ the element in stair form $\sigma_\alpha$ has maximal length in its conjugacy class.
The next lemma corresponds to one step of the algorithm.

	\begin{lem}
		\label{thm:kim_algorithm_and_n-cycles}
		Let $\alpha = (n)$ and $\sigma\in \SG_n$ be an oscillating~$n$-cycle with connected intervals which is different from the element in stair form $\sigma_{\alpha}$. Then there exists a minimal integer $p$ such that $1\leq p\leq n-1$ and $\sigma^p(1) \neq \sigma_\alpha^p(1)$. Set $a := \sigma^p(1)$, $b := \sigma_\alpha^p(1)$ and
		\begin{align*}
		\sigma' := \begin{cases}
			s_{a-1} \sigma s_{a-1} &\myif a > b \\
			s_a \sigma s_a & \myif a < b. \\
		\end{cases}
		\end{align*}
		Then $\sigma' \approx \sigma$ and $\sigma'$ is oscillating and has connected intervals.
	\end{lem}

	\begin{proof}
	Set $I_k := [k, n-k+1]$ for all $k\in \N$ with $k\leq \frac{n+1}2$. 	
	Because $\sigma \neq \sigma_\alpha$ and both permutations are~$n$-cycles, we have $p\leq n-2$. 
	Recall that by \Cref{thm:element_in_stair_form},
	\begin{align*}
		 \sigma_{\alpha} = 
	\begin{cases}
		(1, n, 2, n-1, \dots, \frac n2, \frac n2 +1) & \text{if~$n$ is even} \\
		(1, n, 2, n-1, \dots, \frac{n-1}2,  \frac {n+3}2, \frac {n+1}2) & \text{if~$n$ is odd}.
	\end{cases}
	\end{align*}
	
	If~$n$ is odd then  $\frac{n+1}2 = \sigma^{n-1}(1)$ and hence $p \leq n-2$ implies $b \neq \frac{n+1}2$.
	If~$n$ is even then $b\neq \frac{n+1}{2}$ anyway.
	
	We assume $b < \frac{n+1}2$. The proof in the case $b>\frac{n+1}2$ is similar and therefore omitted. By the choice of~$p$, we have $b\neq 1$ so that $1<b<\frac{n+1}2$. The definition of $\sigma_\alpha$ implies
	\begin{align}
	\label{eq:cycle_of_1_and_I_b}
	\begin{aligned}
	\set{\sigma^k_\alpha(1)\mid k = 0,\dots, p-1} &= [n] \setminus I_b, \\
	\set{\sigma^k_\alpha(1)\mid k = p,\dots, n-1} &= I_b.
	\end{aligned}			
	\end{align}
	Again by the choice of~$p$, the same equalities hold for $\sigma$.
	Hence, $b < a$ as $a\in I_b$ and $b= \min I_b$. Therefore, we consider $\sigma' = s_{a-1}\sigma s_{a-1}$ and show that $\sigma \approx \sigma'$. Then \Cref{thm:equivalence_implies_oscillating_and_ci} implies that $\sigma'$ also is oscillating and has connected intervals.
	
	It follows from the definition of $\sigma_\alpha$ and $b<\frac{n+1}2$ that
	\begin{align}
	\label{eq:thm:kim_algorithm_and_n-cycles:sigma_inv(a)}
	\sigma^{-1} (a) = \sigma_\alpha^{-1}(b) = n-b+ 2 > \frac{n+1}{2}.
	\end{align}
	As $\sigma$ is oscillating, we obtain that $a\leq \frac{n+1}2$ from \Cref{thm:oscillating_n-cycle_local}.
	Since \Cref{eq:cycle_of_1_and_I_b} holds for $\sigma$ and $p>0$,
	\begin{align*}
	\sigma\inv(a) \not \in I_b \supseteq I_{a-1} \supseteq I_{a}.
	\end{align*}
	Let $r := |I_a|$. Because $I_a$ is connected in $\sigma$, $a\in I_a$ and $\sigma^{-1}(a) \not \in I_a$, we have 
	\begin{align*}
	\set{\sigma^k(a) \mid k = 0,\dots, r-1} &= I_{a}.
	\end{align*}
	Now we can use that $I_{a-1} = I_a \cup  \set{a-1, n-a+2}$ is connected in $\sigma$ and that $\sigma\inv(a) \not\in I_{a-1}$ to obtain
	\begin{align*}
	\set{\sigma^k(a) \mid k = 0,\dots, r+1} &= I_{a-1}
	\end{align*}
	The descriptions of $I_a$ and $I_{a-1}$ imply that
	\begin{align*}
	\set{\sigma^{r}(a), \sigma^{r+1}(a)} = \set{a-1, n-a+2}.
	\end{align*}
	\Cref{thm:oscillation_around_intervals} applied to $I_a$ in $\sigma$ and $\sigma^{-1}(a) > \frac{n+1}2$ now imply that $\sigma^{r}(a) < \frac{n+1}2$.
	Thus, $\sigma^r(a) = a-1$ and $\sigma^{r+1}(a) = n-a+2$.
	That is,
	\begin{align}
	\label{eq:thm:kim_algorithm_and_n-cycles:sigma(a-1)}
	\sigma(a-1) = n-a+2
	\end{align}
	Moreover, $\sigma\inv(a-1) \in I_a$ implies
	\begin{align}
	\label{eq:thm:kim_algorithm_and_n-cycles:sigma_inv(a-1)}
	\sigma\inv(a-1) \leq n-a+1.
	\end{align}
	
	We now show
	\begin{align}
	\label{eq:thm:kim_algorithm_and_n-cycles:sigma(a)}
	\sigma(a) \leq n-a+1.
	\end{align}
	and deal with two cases. 
	If $a= \frac{n+1}2$ then $n-a+1 = a$. Furthermore, we then have $r=1$ and therefore $\sigma(a) = a-1 < n-a+1$. 
	If $a< \frac{n+1}2$ then $r>1$ so that $\sigma(a) \in I_a$ and thus $\sigma(a) \leq n-a+1$ as desired.

	From~\cref{eq:thm:kim_algorithm_and_n-cycles:sigma_inv(a-1),eq:thm:kim_algorithm_and_n-cycles:sigma_inv(a)} it follows that
	\begin{align*}
	\sigma^{-1}(a-1) &\leq n-a+1 < n-b+2 = \sigma^{-1}(a).
	\end{align*}
	Moreover, \cref{eq:thm:kim_algorithm_and_n-cycles:sigma(a-1),eq:thm:kim_algorithm_and_n-cycles:sigma(a)} imply
	\begin{align*}
	\sigma(a-1) &= n-a+2 > n-a+1 \geq \sigma(a).
	\end{align*}
	Since $\sigma' = s_{a-1}\sigma s_{a-1}$, \Cref{thm:length_and_conjugation} now yields $\ell(\sigma') = \ell(\sigma)$.
	Hence, $\sigma' \approx \sigma$ by \Cref{thm:arrow_and_length}.
\end{proof}

\begin{exa}
	Let $n = 5$ and $\alpha = (n)$.
	The $n$-cycle $\sigma = (1,3,4,2,5)\in\SG_n$ is oscillating and has connected intervals.
	We can successively use \Cref{thm:kim_algorithm_and_n-cycles} in order to obtain the sequence
	\begin{align*}
		\sigma = \sigma^{(0)} &= (1,3,4,2,5), \\
		\sigma^{(1)}&= (1,4,3,2,5) = s_3 \sigma^{(0)} s_3,\\
		\sigma^{(2)}&= (1,5,3,2,4) = s_4 \sigma^{(1)} s_4,\\
		\sigma^{(3)}&= (1,5,2,3, 4) = s_2 \sigma^{(2)} s_2, \\
				\sigma_\alpha = \sigma^{(4)}&= (1,5,2,4,3) = s_3 \sigma^{(3)} s_3.
	\end{align*}
	Moreover, \Cref{thm:kim_algorithm_and_n-cycles} ensures that each $\sigma^{(j)}$ is oscillating with connected intervals and all $\sigma^{(j)}$ are $\approx$-equivalent.
	Therefore, $\sigma \in \Sigma_\alpha$ by \Cref{thm:parametrizations_of_kim}.
\end{exa}

We now come to the characterization of $\Sigma_{(n)}$.

\begin{thm}	
	\label{thm:characterization_of_Sigma_(n)}
	Let $\sigma\in \SG_n$ be an~$n$-cycle.
	Then $\sigma \in \Sigma_{(n)}$ if and only if $\sigma$ is oscillating and has connected intervals.
\end{thm}
\begin{proof}
	
	Let $\sigma\in \SG_n$ be an~$n$-cycle.
	Recall that $\sigma \in \Sigma_{(n)}$ if and only if $\sigma \approx \sigma_{(n)}$ by \Cref{thm:parametrizations_of_kim}.
	Assume that $\sigma \in \Sigma_{(n)}$. Then $\sigma \approx \sigma_{(n)}$ which by definition of $\approx$ implies that there are sequences $\sigma_\alpha = \sigma^{(0)}, 	\sigma^{(1)}, \dots, \sigma^{(m)} = \sigma \in \SG_n$
	and $i_1, \dots, i_m \in [n-1]$ such that $\sigma^{(j-1)} \approx \sigma^{(j)}$ and $\sigma^{(j)} = s_{i_j}\sigma^{(j-1)}s_{i_j}$ for $j\in [m]$.
	From \Cref{thm:sigma_(n)_osc_and_c.I.} we have that $\sigma_{(n)}$ is oscillating and has connected intervals.
	Moreover, \Cref{thm:equivalence_implies_oscillating_and_ci} yields that $\sigma^{(j)}$ is oscillating with connected intervals if 
	  $\sigma^{(j-1)}$ is oscillating with connected intervals.
	  Hence, $\sigma$ is oscillating and has connected intervals by induction.
	
	Conversely, assume that $\sigma$ is oscillating and has connected intervals. Then we can use \Cref{thm:kim_algorithm_and_n-cycles} iteratively to obtain a sequence of $\approx$-equivalent~$n$-cycles starting with $\sigma$ and eventually ending with $\sigma_{\alpha}$. Thus $\sigma \approx \sigma_\alpha$.
\end{proof}

\begin{figure}	
	\begin{tikzcd}[column sep=huge, row sep = huge]
		  (1,3,4,2,5) \arrow[shift left]{dr}{\del_3} 
		& (1,4,3,2,5) \arrow[shift left]{d}{\del_3} 
		& (1,4,2,3,5) \arrow[shift left]{dl}{\del_3} \\
		&(1,3,2,4) 	\arrow[shift left]{d}{\del_3}
					\arrow[shift left]{ul}{\ins_{3,1}}
					\arrow[shift left]{u}{\ins_{3,2}}
					\arrow[shift left]{ur}{\ins_{3,3}}  &\\
		&(1,2,3) 	\arrow[shift left]{u}{\ins_{3,1}} &
	\end{tikzcd}
	\caption{%
		Examples for the operators $\del_k$ and $\ins_{k,p}$ appearing in \Cref{thm:recursion_for_Sigma_(n)} and its proof. 
		The lower part of the picture serves as an example for the operators used in the case when~$n$ is even. The upper part is an example for those used in the case when~$n$ is odd.
		Note that for the integer~$m$ from the \namecref{thm:recursion_for_Sigma_(n)} we have $m = \frac n2 +1 = 3$ if $n=4$ and $m=\frac{n+1}2=3$ if $n=5$.}
	\label{fig:del_and_ins}
\end{figure}
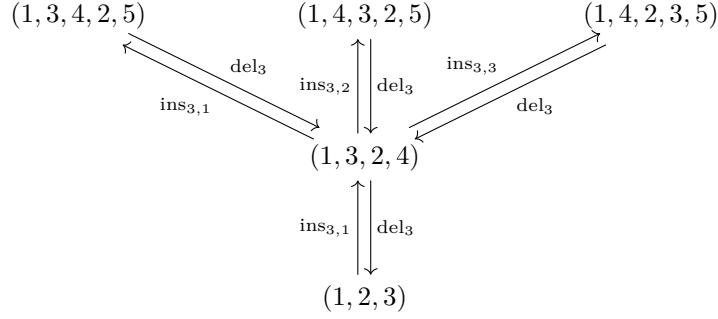

The goal of the remainder of this \namecref{sec:equivalence_classes:one_part} is to find bijections that relate $\Sigma_{(n-1)}$ to $\Sigma_{(n)}$. 
From this we will obtain a recursive description of $\Sigma_{(n)}$ and a formula for the cardinality of $\Sigma_{(n)}$.
To achieve our goal, we define two operators $\ins$ and $\del$.

Assume that the $n$-cycle $\sigma \in \SG_n$ is given in cycle notation starting with~$1$.
Then for $k\in [2,n+1]$ $\ins_{k,p}(\sigma)\in \SG_{n+1}$ is the $(n+1)$-cycle obtained from $\sigma$ by adding~$1$ to each element greater or equal to~$k$ in $\sigma$ and then inserting~$k$ behind the~$p$th element in the resulting cycle.
Likewise, for $k\in [2,n]$, $\del_{k}(\sigma) \in \SG_{n-1}$ is the $(n-1)$-cycle obtained by first deleting~$k$ from $\sigma$ and then decreasing each element greater than~$k$ by~$1$.
See \Cref{fig:del_and_ins} for examples.

We now define $\ins$ and $\del$ more formally.
Let $\sigma\in \SG_n$ be an~$n$-cycle and $k\in \N$.
Set
\begin{align*}
	\varepsilon_r := 
	\begin{cases}
		0 &\myif \sigma^r(1) < k \\
		1 &\myif \sigma^r(1) \geq k
	\end{cases}
\end{align*}
for $r = 0 ,\dots, n-1$.
In the following we will assume $k>1$. 
The operators could also be defined for $k=1$ but this is not necessary for our purposes and would only make the exposition less transparent.

For $k\in [2,n+1]$ and $p\in [n]$, define $\ins_{k,p}(\sigma)$ to be the $(n+1)$-cycle of $\SG_{n+1}$ given by 
\begin{align*}
\ins_{k,p}(\sigma)^r(1) :=
\begin{cases}
\sigma^r(1) + \varepsilon_r &\myif r < p \\
k 							&\myif r = p \\
\sigma^{r-1}(1) + \varepsilon_{r-1} &\myif r > p
\end{cases}
\end{align*}
for $r = 0, \dots, n$.
For $k\in [2,n]$, define $\del_k(\sigma)$ to be the $(n-1)$-cycle of $\SG_{n-1}$ given by
\begin{align*}
	\del_k(\sigma)^r(1) := 
	\begin{cases}
		\sigma^r(1) - \varepsilon_r &\myif r<p \\
		\sigma^{r+1}(1) - \varepsilon_{r+1} &\myif r\geq p
	\end{cases}
\end{align*}
for $r = 0,\dots, n-2$ where~$p$ is the element of $[0,n-1]$ with $\sigma^p(1) = k$.

The next results relates $\Sigma_{(n)}$ with $\Sigma_{(n-1)}$ via a bijection for $n \geq 4$. 

\begin{thm}
	\label{thm:recursion_for_Sigma_(n)}
	Suppose $n\geq 4$.
	If~$n$ is even then set $m := \frac{n}2 + 1$ and 
	\begin{align*}
		\psi \colon \Sigma_{(n-1)} \to \Sigma_{(n)}, \quad \sigma \mapsto \ins_{m,p}(\sigma)
	\end{align*}
	where~$p$ is the element of $[n-1]$ with $\sigma^{p-1}(1) = \min\set{\sigma\inv(\frac n2), \frac n2 }$.
	If~$n$ is odd then set $m := \frac{n+1}2$ and
	\begin{align*}
	\psi \colon \Sigma_{(n-1)} \times \set{0,1,2} \to \Sigma_{(n)}, \quad (\sigma,q) \mapsto \ins_{m, p+q}(\sigma)
	\end{align*}
	where~$p$ is the element of $[n-3]$ with $\sigma^{p-1}(1) \not\in \set{m-1,m}$ and $\sigma^p(1) \in \set{m-1,m}$.
	Then $\psi$ is a bijection.
\end{thm}

\begin{cor}
	\label{thm:recurrences_of_Sigma_n}
	Suppose $n\geq 4$. Then
	\begin{align*}
	\card{\Sigma_{(n)}} = 
	\begin{cases}
	\card{\Sigma_{(n-1)}} & \myif \text{$n$ is even} \\
	3\card{\Sigma_{(n-1)}} & \myif \text{$n$ is odd}.
	\end{cases}
	\end{align*}
\end{cor}

\begin{proof}[Proof of \Cref{thm:recursion_for_Sigma_(n)}]
	\Cref{thm:characterization_of_Sigma_(n)} states that for all $n\in \N$, $\Sigma_{(n)}$ is the set of oscillating~$n$-cycles of $\SG_n$ with connected intervals.
	In this proof we repeatedly use this result without further notice.

	Let $n\geq 4$. 
	We consider all permutations in the cycle notation where~$1$ is the leftmost entry in its cycle. 
	In particular, deleting an entry from a permutation or inserting an entry into a permutation means that we do this in the chosen cycle notation.
	We distinguish two cases depending on the parity of~$n$.
	
	\begin{caseenum}	
		\item 
		Assume that~$n$ is even. 
		Then $m= \frac{n}2+1$.
		For $\tau \in \Sigma_{(n-1)}$ let~$p$ be given as in the definition of $\psi$.
		Then $\min\set{\tau\inv(\frac n2), \frac n2 }$ is the~$p$th element in the cycle notation of $\tau$.
		Hence, we obtain $\psi(\tau)$ by increasing each element in $\tau$ greater or equal to~$m$ by one and then inserting~$m$ behind the element at position~$p$.
		
		Set 
		$
		\varphi \colon \Sigma_{(n)} \to \Sigma_{(n-1)}, \sigma \mapsto \del_m(\sigma).$
		That is, for $\sigma \in \Sigma_{(n)}$ we obtain $\varphi(\sigma)$ by first deleting~$m$ from $\sigma$ and then decreasing each entry greater than~$m$ by~$1$.
		
		We show that $\varphi$ and $\psi$ are well defined and inverse to each other.
		\begin{wideenumerate}
			\item
			We prove that $\varphi$ is well defined. Let $\sigma \in \Sigma_{(n)}$ and $\tau := \varphi(\sigma)$. 
			We have to show that $\tau\in \Sigma_{(n-1)}$.
			That is, we have to prove that $\tau$ is oscillating and has connected intervals. 
			
			To show the latter, let $1\leq i \leq \frac{n-1}2 < \frac{n}2$. As $[i,n-i+1]$ is connected in $\sigma$ there is a $0\leq q\leq n-1$ such that
			\begin{align*}
			\set{\sigma^{q+1}(1), \dots, \sigma^{q+r}(1)} = [i,n-i+1]
			\end{align*}
			where $r:= \card{[i,n-i+1]}$.
			Moreover, $m \in [i,n-i+1]$.
			Thus, $\tau = \del_m(\sigma)$ implies 
			\begin{align*}
			\set{\tau^{q+1}(1), \dots, \tau^{q+r-1}(1)} = [i,n-i].
			\end{align*}
			Hence, $[i,(n-1)-i+1]$ is connected in $\tau$.
			It follows that $\tau$ has connected intervals.
			
			We now show that $\tau$ is oscillating.
			Note that $n-1$ is odd and $\frac{(n-1)+1}2 = \frac{n}2$.
			By \Cref{thm:oscillating_n-cycle_local}, it suffices to show that $\tau(i) \geq \frac{n}2$ for all $i\in [\frac{n}2-1]$ and that
			either  $\tau^{-1}\left(\frac n2\right) > \frac n2$ or $\tau\left(\frac n2\right) > \frac n2$.

			Let $i\in [\frac{n}2-1]$.
			Since  $i < \frac{n}2$ and $\sigma$ is oscillating, we infer $\sigma(i) > \frac n2$ from \Cref{thm:oscillating_n-cycle_local}.
			If $\sigma(i) \neq m$ then $\tau(i) = \sigma(i) -1 \geq \frac n2$. 
			If $\sigma(i) = m$ then $\sigma^2(i) = \frac{n}2$ since 
			$m = \frac{n}2+1$, $\set{\frac n2, \frac n2 + 1}$ is connected in $\sigma$ and $i\not \in \set{\frac n2, \frac n2 + 1}$.
			Thus, $\tau(i) = \frac n2$.

			We now show that either $\tau^{-1}\left(\frac n2\right) > \frac n2$ or $\tau\left(\frac n2\right) > \frac n2$.
			Since $\set{\frac n2, \frac n2 + 1}$ is connected in $\sigma$ there is a $0\leq q \leq n-1$ such that 
			\begin{align*}
			\set{\sigma^q(1), \sigma^{q+1}(1)} = \set{\frac n2, \frac n2 + 1}.
			\end{align*}
			Hence, $\tau = \del_{\frac{n}2+1}(\sigma)$ implies $\tau^q(1) = \frac n2$.
			Because $n\geq 4$, we can apply \Cref{thm:oscillation_around_intervals} to $\set{\frac n2, \frac n2 + 1}$ in $\sigma$ and obtain that there are $a<\frac n2$ and $b>\frac n2+1$ such that 
			\begin{align*}
			\set{\sigma^{q-1}(1), \sigma^q(1), \sigma^{q+1}(1),  \sigma^{q+2}(1)} = \set{a,b,\frac n2, \frac n2 + 1}.
			\end{align*}
			Therefore, $\tau^q(1) = \frac n2$ and $\tau = \del_{\frac{n}2+1}(\sigma)$ yield
			$
			\set{\tau^{-1}\left(\frac n2\right), \tau\left(\frac n2\right) } = \set{a,b-1}.
			$
			That is, either $\tau^{-1}\left(\frac n2\right) > \frac n2$ or $\tau\left(\frac n2\right) > \frac n2$.
			Thus, $\tau$ is oscillating.
			
			\item
			We check that $\psi$ is well defined.
			Let $\tau \in \Sigma_{(n-1)}$ and $\sigma := \psi(\tau)$. 
			We have to show $\sigma \in \Sigma_{(n)}$.

			The definition of $\psi$ implies that $\frac n2+1$ is a neighbor of $\frac n2$ in $\sigma$.
			In addition, $[i, n-i]$ is connected in $\tau$ for $i\in [\frac{n}2-1]$.
			Therefore, $[i, n-i+1]$ is connected in $\sigma$ for $i\in [\frac{n}2]$.
			That is, $\sigma$ has connected intervals.

			We now show that $\sigma$ is oscillating. 
			By \Cref{thm:oscillating_n-cycle_local}, it suffices to show that $\sigma(i) > \frac{n}2$ for all $i \in [\frac{n}2]$.
			For $i <\frac{n}2$ this can be done as before.
			Thus, we only consider $i=\frac{n}2$.
			As $\tau$ is oscillating, \Cref{thm:oscillating_n-cycle_local} implies that one of the neighbors of $\frac n2$ is smaller than $\frac n2$ and the other one is greater than $\frac n2$.
			Let~$a$ be the smaller and~$b$ be the bigger neighbor of $\frac n2$. 
			In the definition of $\psi$,~$p$ is chosen such that $\frac{n}2+1$ is inserted in $\tau$ between~$a$ and $\frac n2$. 
			Thus, $\frac n2$ has neighbors $\frac n2 +1$ and $b+1$ in $\sigma$.
			Consequently, $ \sigma\left(\frac n2\right) > \frac n2$.
			
			\item
			We now show that $\psi \circ \varphi = \id$.
			Let $\sigma \in \Sigma_{(n)}$. Since $\set{\frac n2, \frac n2+1}$ is connected in $\sigma$, these two elements are neighbors in $\sigma$.
			As $\sigma$ is oscillating, there is an $a<\frac n2$ such that $\frac n2 +1$ has neighbors~$a$ and $\frac n2$.
			We obtain 
			$\varphi(\sigma)$ from $\sigma$ by deleting $\frac n2 +1$ so that~$a$ and $\frac n2$ are neighbors in $\varphi(\sigma)$. On the other hand, we obtain $\psi(\varphi(\sigma))$ from  $\varphi(\sigma)$ by inserting $\frac n2 +1$ between~$a$ and $\frac n2$. Thus $\psi(\varphi(\sigma))= \sigma$.
			
			\item
			Finally, we show that $  \varphi \circ \psi= \id$. Let $\tau \in \Sigma_{(n-1)}$. 
			Then we obtain $\psi(\tau)$ from $\tau$ by inserting $\frac n2+1$ at some position and get $\varphi(\psi(\tau))$ from $\psi(\tau)$ by deleting it again. Hence, $\varphi(\psi(\tau)) = \tau$.
			
		\end{wideenumerate}
		
		\item 
		Assume that~$n$ is odd.
		Then $m = \frac{n+1}2$.
		For $\tau \in \Sigma_{(n-1)}$ the set $\set{m-1,m}$ is connected.
		Thus, there is a unique integer $p$ with $1\leq p \leq n-3$ such that $\tau^{p-1}(1) \not\in \set{m-1,m}$ and $\tau^p(1) \in \set{m-1,m}$.
		That is, the integer $p$ from the definition of $\psi$ in the \namecref{thm:recursion_for_Sigma_(n)} is well defined.
		Note that~$p$ is the position of the left neighbor of the set $\set{m-1,m}$ in $\tau$.
		
		Conversely, for $\sigma \in \Sigma_{(n)}$, $I := \set{m-1, m, m+1}$ is connected in $\sigma$.
		Hence, there is a unique $0\leq p \leq n-1$ such that $I = \set{\sigma^{p+k}(1) \mid k= 0,1,2}$ and a unique $q\in \set{0,1,2}$ such that $\sigma^{p+q}(1) = m$.
		We define the map	$\varphi \colon \Sigma_{(n)} \to \Sigma_{(n-1)} \times \set{0,1,2}$ by setting $\varphi(\sigma) := (\del_m(\sigma), q)$.
		Again, we show that $\varphi$ and $\psi$ are well defined and inverse to each other.

		\begin{wideenumerate}
			\item
			First we show that the two maps are inverse to each other.
			Let $\sigma \in \Sigma_{(n)}$ and $\varphi(\sigma) = (\tau,q)$. Then we have
			\begin{align*}
			q = 
			\begin{cases}
			0 &\myif \text{$m$ is the left neighbor of $\set{m-1, m+1}$ in $\sigma$}, \\
			1 &\myif \text{$m$ is located between $m-1$ and $m+1$ in $\sigma$}, \\
			2 &\myif \text{$m$ is the right neighbor of $\set{m-1, m+1}$ in $\sigma$}. \\
			\end{cases}
			\end{align*}
			Conversely, let $\tau \in \Sigma_{(n-1)}$, $q\in \set{0,1,2}$ and $\sigma= \psi(\tau,q)$ then
			\begin{align}
			\label{eq:description_of_ins_m_p+q}
			\text{$m$ is }
			\begin{cases}
			\text{the left neighbor of $\set{m-1, m+1}$ in $\sigma$} &\myif q = 0,\\
			\text{located between $m-1$ and $m+1$ in $\sigma$} &\myif q=1,\\
			\text{the right neighbor of $\set{m-1, m+1}$ in $\sigma$}  &\myif q = 2.\\
			\end{cases}
			\end{align}
			From this it follows that $\varphi$ and $\psi$ are inverse to each other.
			
			\item
			In order to prove that $\varphi$ is well defined one has to show that $\del_m(\sigma) \in \Sigma_{(n-1)}$. This can be done similarly as in Case 1.
			
			\item
			To see that $\psi$ is well defined, let $\tau \in \Sigma_{(n-1)}$, $q\in \set{0,1,2}$ and $\sigma := \psi(\tau,q)$. 
			We first show that $\sigma$ has connected intervals.
			Recall that $m = \frac{n+1}2$. Let $i\leq \frac{n-1}2 = m-1$. Then $[i, n-i]$ is connected in $\tau$ since $\tau$ has connected intervals. 
			By the definition of $\psi$, we obtain the entries $[i,n-i+1]$ in $\sigma$ by adding~$1$ to each entry $\geq m$ of $[i, n-i]$ in $\tau$ and then inserting~$m$ such that by \Cref{eq:description_of_ins_m_p+q} at least one of the neighbors of~$m$ is $m-1$ or $m+1$.
			Since  $m-1,m,m+1 \in [i,n-i+1]$ it follows that $[i, n-i+1]$ is connected in $\sigma$.
			Therefore, $\sigma$ has connected intervals.
							
			In order to show that $\sigma$ is oscillating, let $\tau'$ be the $(n-1)$-cycle of $\SG_n$ obtained by adding~$1$ to each entry of $\tau$ which is greater or equal than~$m$. Since $\tau$ is oscillating, the entries in $\tau'$ alternate between the sets $[m-1]$ and $[m+1,n]$. Furthermore, we obtain $\sigma$ from $\tau'$ by inserting~$m$ somewhere in $\tau'$.
			Thus, \Cref{thm:oscillating_n-cycle_cycle_notation} implies that $\sigma$ is oscillating.
			\qedhere
		\end{wideenumerate}
	\end{caseenum}
\end{proof}

From \Cref{tbl:Sigma_(n)} we know $\Sigma_{(n)}$ for $n = 1,2,3$. That is, \Cref{thm:recursion_for_Sigma_(n)} allows us to compute $\Sigma_{(n)}$ recursively for each $n\in \N$. This is illustrated in the following.

\begin{exa}
	We want to compute $\Sigma_{(n)}$ for $n=4,5$. To do this we use the bijections $\psi$ and the related notation introduced in \Cref{thm:recursion_for_Sigma_(n)}.
	\begin{wideenumerate}
	\item
	Consider $n=4$.
	We have 
	\begin{align*}
		\Sigma_{(4)} = \set{\psi(\sigma) \mid \sigma \in \Sigma_{(3)}}
	\end{align*}
	by \Cref{thm:recursion_for_Sigma_(n)}.
	From \Cref{tbl:Sigma_(n)} we obtain $\Sigma_{(3)} = \set{(1,3,2),(1,2,3)}$.
	
	For $\sigma = (1,3,2)$ we have $p=3$ since 
	\begin{align*}
		\sigma^{3-1}(1) = 2 = \min\set{2,3} = \min\set{\sigma\inv\left(\frac{4}2\right), \frac{4}2}.
	\end{align*}
	Thus,
	\begin{align*}
		\psi(\sigma) = \ins_{3,3}((1,3,2)) = (1,3+1,2,3) = (1,4,2,3).
	\end{align*}
	For $\sigma = (1,2,3)$ we have $p=1$ and 
	\begin{align*}	
		\psi(\sigma) = \ins_{3,1}((1,2,3)) = (1,3,2,3+1) = (1,3,2,4).
	\end{align*}
	Therefore, $\Sigma_{(4)} = \set{(1,4,2,3), (1,3,2,4)}$.
	
	\item
	Consider $n=5$. \Cref{thm:recursion_for_Sigma_(n)} yields
	\begin{align}
	\label{eq:Sigma_(5)_recursion}
	\Sigma_{(5)} = 
	\set{\psi(\sigma, q) \mid \sigma \in \Sigma_{(4)}, q\in \set{0,1,2}}.
	\end{align}
	Let $m = \frac{5+1}2 = 3$ and $I = \set{m-1,m} = \set{2,3}$.
	
	For $\sigma = (1,4,2,3)$ we have $p=2$ since 
	$\sigma^{2-1}(1) = 4 \not \in I$ and $\sigma^{2}(1) = 2 \in I$.
	Thus, for instance we have
	\begin{align*}
		\psi(\sigma,1) = \ins_{3,3}((1,4,2,3)) = (1,4+1,2,3,3+1) = (1,5,2,3,4). 
	\end{align*}
	
	For $\sigma = (1,3,2,4)$ we have $p=1$.
	Computing $\psi(\sigma,q)$ for all $\sigma\in \Sigma_{(4)}$ and $q\in \set{0,1,2}$, we obtain the following table. By \Cref{eq:Sigma_(5)_recursion}, it lists all elements of $\Sigma_{(5)}$.
	\begin{align*}
	\begin{array}{c|ccc}
	\psi(\sigma,q)	& 0           & 1           & 2           \\ \hline
		(1,4,2,3)   & (1,5,3,2,4) & (1,5,2,3,4) & (1,5,2,4,3) \\
		(1,3,2,4)   & (1,3,4,2,5) & (1,4,3,2,5) & (1,4,2,3,5)
	\end{array}
	\end{align*}
		\end{wideenumerate}
\end{exa}

\begin{cor}
	\label{thm:sizes_of_Sigma_n}
	Let $n\in \N$. Then
	\begin{align*}
	\card{\Sigma_{(n)}} = 
	\begin{cases}
	1 & \myif \text{$n\leq2$} \\
	2 \cdot 3^{\left \lfloor{\frac{n-3}2}\right \rfloor } & \myif \text{$n\geq 3$}. \\
	\end{cases}
	\end{align*}
\end{cor}

\begin{proof}
	
		Let $x_n := |\Sigma_{(n)}|$ for $n\geq 1$, $y_1 := y_2 := 1$ and $y_n := 2\cdot 3^{\left \lfloor{\frac{n-3}2}\right \rfloor }$ for $n\geq 3$. We show that both sequences have the same initial values and recurrence relations. First note that 
		\begin{align*}
		(x_1,x_2,x_3) = (1,1,2) = (y_1,y_2,y_3).
		\end{align*}
		where we obtain the $x_i$ from \Cref{tbl:Sigma_(n)}.
		Now let $n\geq 4$. By \Cref{thm:recurrences_of_Sigma_n} we have to show that $y_n = y_{n-1}$ if~$n$ is even and $y_n = 3y_{n-1}$ if~$n$ is odd. If~$n$ is even, we have
		\begin{align*}
		\left\lfloor{\frac{n-3}2}\right \rfloor
		= \left\lfloor{\frac{n-4}2 + \frac 12}\right \rfloor
		= \frac{n-4}2 
		= \left\lfloor{\frac{n-1-3}2 }\right \rfloor
		\end{align*}
		and thus $y_n = y_{n-1}$.
		If~$n$ is odd, we have
		\begin{align*}
		\left\lfloor{\frac{n-3}2}\right \rfloor
		= \frac{n-3}2
		= \frac{n-5}2 + 1
		= \left\lfloor{\frac{n-5}2 + \frac 12}\right \rfloor + 1
		= \left\lfloor{\frac{n-4}2 }\right \rfloor + 1
		\end{align*}
		and hence $y_n = 3y_{n-1}$.
\end{proof}

\subsection{Equivalence classes of odd hook type}
\label{sec:equivalence_classes:odd_hook}

Let $\alpha = (k, 1^{n-k}) \vDash n$ be a hook.
Then $\alpha$ is a maximal composition.
Recall that
a hook $\alpha$ is called \emph{odd} if~$k$ is odd and called \emph{even} otherwise.
The main result of this \namecref{sec:equivalence_classes:odd_hook} is a combinatorial characterization of $\Sigma_\alpha$ provided that $\alpha$ is an odd hook in \Cref{thm:characterization_of_Sigma_for_odd_hook}.

We will use the inductive product which is the topic of \Cref{sec:equivalence_classes:inductive_product} in order to deal with the even hooks.
Therefore, the characterization of $\Sigma_\alpha$ for even hooks $\alpha$ is postponed until \Cref{thm:characterization_of_Sigma_for_arbitrary_hook}.

We want to generalize the concept of being oscillating and having connected intervals from $n$-cycles to arbitrary permutations.
In order to do this, we standardize cycles in the following way.
Let $\sigma := (c_1,\dots, c_k)\in \SG_n$ be a $k$-cycle.
Replace the smallest element among $c_1, \dots, c_k$ by $1$, the second smallest by $2$ and so on. 
The result is a $k$-cycle with entries $1,2,\dots, k$ which can be regarded as an element $\SG_k$.
This permutation is called the \emph{cycle standardization} $\cst(\sigma)$ of $\sigma$.

\begin{exa}	Consider $\sigma = (3,11,4,10,5)\in \SG_{11}$. Then $\cst(\sigma) = (1,5,2,4,3)\in \SG_5$ which is oscillating with connected intervals.
\end{exa}

We formally define the cycle standardization as follows.

\begin{defi}
	\label{thm:cycle_standartization}
	\begin{enumerate}
		\item Given $\sigma \in \SG_n$ and $i\in [n]$, there is a cycle $(c_1,\dots, c_k)$ of $\sigma$ containing~$i$. Then we define
		\begin{align*}
		\rnk_\sigma(i) := \card{\set{j \in [k] \mid c_j\leq i}}.
		\end{align*}
		\item  Let $\sigma = (c_1, \dots, c_k)\in \SG_n$ be a~$k$-cycle.
		The \emph{cycle standardization} of $\sigma$ is the~$k$-cycle of $\SG_k$ given by
		\begin{align*}
		\cst(\sigma) := (\rnk_\sigma(c_1), \rnk_\sigma(c_2), \dots, \rnk_\sigma(c_k)).
		\end{align*}
	\end{enumerate}	
\end{defi}

Note that the permutation $\cst(\sigma)$  is independent from the choice of the cycle notation $\sigma = (c_1, c_2, \dots, c_k)$ in  \Cref{thm:cycle_standartization}.

\begin{rem}
	\label{thm:properties_of_cycle_standartization}
	Let $\sigma = (c_1, c_2, \dots, c_k)\in \SG_n$ be a~$k$-cycle.
	\begin{enumerate} 
		\item  The anti-rank of $i\in [n]$ among the elements in its cycle in $\sigma$ is $\rnk_\sigma(i)$.
		\item
		\label{enum:cycle_standartization_same_relative_order}
		For all $i,j\in [k]$ we have $c_i < c_j$ if and only if  $\rnk_\sigma(c_i) < \rnk_\sigma(c_j)$.
		\item Let~$i$ be an element appearing in the cycle $(c_1, c_2, \dots, c_k)$. Then we have 
		\begin{align*}
		\cst(\sigma)(\rnk_\sigma(i)) = \rnk_\sigma( \sigma(i)).
		\end{align*}
	\end{enumerate}
\end{rem}

We now generalize the notions of being oscillating and having connected intervals to arbitrary permutations via the cycle decomposition and the cycle standardization.
Recall that trivial cycles are those of length~$1$.

\begin{defi}
	\label{def:oscillating_and_connected_intervals_general}
	Let $\sigma\in \SG_n$ and write $\sigma$ as a product $\sigma = \sigma_1\cdots \sigma_l$ of disjoint cycles including the trivial ones.
	\begin{enumerate} 
		\item We say that $\sigma$ is \emph{oscillating} if $\cst(\sigma_i)$ is oscillating for each cycle $\sigma_i$.
		\item We say that $\sigma$ has \emph{connected intervals} if $\cst(\sigma_i)$ has connected intervals for each cycle $\sigma_i$
	\end{enumerate}
\end{defi}

Let $(c)\in \SG_n$ be a trivial cycle.
Then $\cst((c)) = (1) \in \SG_1$ which is oscillating and has connected intervals.
Therefore, in order to show that a permutation $\sigma$ is oscillating (has connected intervals) it suffices to consider the nontrivial cycles.

\begin{exa}
	Let $\alpha = (4,5,3,1)\vDash_e 13$ and
	\begin{align*}
	\sigma_\alpha = (1, 13, 2, 12)(3, 11, 4, 10, 5)( 9, 6, 8)(7).
	\end{align*}
	The cycle standardizations of the nontrivial cycles of $\sigma_\alpha$ are
	\begin{align*}
		(1,4,2,3), (1,5,2,4,3) \text{ and } (1,2,3).
	\end{align*}
	Each of these three permutations is oscillating and has connected intervals (see \Cref{tbl:Sigma_(n)}).
	Thus, $\sigma_\alpha$ is oscillating and has connected intervals.
\end{exa}

Assume that $\sigma \in \SG_n$ is an~$n$-cycle.
Then $\sigma$ has only one cycle $\sigma$ in cycle notation and $\cst(\sigma) = \sigma$.
Thus, for~$n$-cycles our new notion of being oscillating (having connected intervals)  from \Cref{def:oscillating_and_connected_intervals_general} is equivalent to the old concept from \Cref{def:oscillating_n-cycle}  (\Cref{def:connected_intervals_n-cycle}).

We now prove some general results on oscillating permutations with connected intervals.
As in the last \namecref{sec:equivalence_classes:one_part}, we are interested in the effect of swapping entries $i$ and $i+1$ in cycle notation (that is,  conjugating with $s_i$).
This will in particular be useful to prove our results on odd hooks.
We consider the case where $i$ and $i+1$ appear in the same cycle first.

\begin{lem}
	\label{thm:equivalence_and_interchange_within_cycle}
	Let $\sigma \in \SG_n$ and write $\sigma$ as a product $\sigma = \sigma_1\cdots \sigma_l$ of disjoint cycles.
	Assume that there is an $i\in [n-1]$ and a $k\in [l]$ such that~$i$ and $i+1$ both appear in the cycle $\sigma_k$.
	Set $i' := \rnk_\sigma(i)$ and $\tau := \cst(\sigma_k)$.
	Then we have
	\begin{enumerate}
		\item $\cst(s_i\sigma_k s_i) = s_{i'} \tau s_{i'}$,
		\item 	$s_i \sigma s_i \approx \sigma$ if and only if $ s_{i'} \tau s_{i'} \approx \tau$.
	\end{enumerate}
\end{lem}

\begin{proof}
	By the definition of $\rnk_\sigma$, we have that $\rnk_\sigma(j) = \rnk_{\sigma_k}(j)$ for all entries~$j$ in the cycle $\sigma_k$.

	(1)
	We obtain $s_i\sigma_k s_i$ from $\sigma_k$ by interchanging~$i$ and $i+1$ in cycle notation.
	Since~$i$ and $i+1$ appear in $\sigma_k$, we have $\rnk_{\sigma_k}(i+1) = i'+1$.
	Thus, we obtain $\cst(s_i\sigma_k s_i)$ from $\tau = \cst(\sigma_k)$ by interchanging $i'$ and $i'+ 1$ in cycle notation.
	That is, $\cst(s_i\sigma s_i) = s_{i'} \tau s_{i'}.$
	
	(2)
	We have
	$s_i \sigma s_i \approx \sigma $ if and only if $\ell(s_i \sigma s_i) = \ell(\sigma)$.
	By \Cref{thm:length_and_conjugation}, this is the case if and only if either $\sigma(i) < \sigma(i+1)$ or $\sigma\inv(i) < \sigma\inv(i+1)$. 
	From the definition of the cycle standardization  we obtain that $\tau(\rnk_\sigma(j)) = \rnk_\sigma (\sigma(j))$ for each entry~$j$ in $\sigma_k$ (\cf \Cref{thm:properties_of_cycle_standartization}). 
	Moreover, by the definition of $\rnk_\sigma$ and the fact that~$i$ and $i+1$ appear in the same cycle of $\sigma$,
	\begin{align*}
	\sigma(i) < \sigma(i+1) \iff \rnk_\sigma(\sigma(i)) < \rnk_\sigma(\sigma(i+1)).
	\end{align*}
	Hence,
	\begin{align*}
	\sigma(i) < \sigma(i+1) \iff \tau(i') < \tau(i' + 1).
	\end{align*}
	Similarly, one shows that this equivalence is also true for $\sigma\inv$ and $\tau\inv$.
	Therefore, we have $s_i \sigma s_i \approx \sigma$ if and only if either $\tau(i') < \tau(i'+1)$ or $\tau\inv(i') < \tau\inv(i'+1)$.
	As for $\sigma$, the latter is equivalent to $ s_{i'} \tau s_{i'} \approx \tau$.
\end{proof}

We now infer from \Cref{thm:equivalence_and_interchange_within_cycle} that swaps of $i$ and $i+1$ within a cycle that preserve $\approx$ also preserve the properties of being oscillating with connected intervals.

\begin{cor}
	\label{thm:equivalence_osc+c.I_and_interchange_within_cycle}
	Let $\sigma \in \SG_n$ be oscillating with connected intervals, $i\in [n-1]$ such that~$i$ and $i+1$ appear in the same cycle of $\sigma$ and $\sigma' := s_i \sigma s_i$. If $\sigma \approx \sigma'$ then $\sigma'$ is oscillating with connected intervals.
\end{cor}

\begin{proof}
	We write $\sigma$ as a product $\sigma = \sigma_1\cdots \sigma_l$ of disjoint cycles and choose $k$ such that~$i$ and $i+1$ appear in the cycle $\sigma_k$. 
	Moreover, we set $\tau := \cst(\sigma_k)$, $\tau' := \cst(s_i\sigma_k s_i)$ and~$m$ to be the length of the cycle $\sigma_k$.
	
	As~$i$ and $i+1$ only appear in $\sigma_k$, $\sigma' = \sigma_1 \cdots \sigma_{k-1} ( s_i \sigma_k s_i) \sigma_{k+1} \cdots \sigma_l$ is the decomposition of $\sigma'$ in disjoint cycles. Since $\sigma$ is oscillating with connected intervals, $\cst(\sigma_j)$ is oscillating with connected intervals for all $j\in [l]$.
	Therefore, it remains to show that $\tau'$ has these properties.
	Since $\sigma \approx \sigma'$, \Cref{thm:equivalence_and_interchange_within_cycle} yields that $\tau \approx \tau'$.
	In addition, $\tau$ is an oscillating~$m$-cycle with connected intervals and thus $\tau \in \Sigma_{(m)}$ by \Cref{thm:characterization_of_Sigma_(n)}.
	Hence, also $\tau'\in \Sigma_{(m)}$, \ie $\tau'$ is oscillating with connected intervals.
\end{proof}

The next result is concerned with the interchange of $i$ and $i+1$ between two cycles.

\begin{lem}
	\label{thm:equivalence_osc+c.I_and_interchange_between_cycles}
	Let $\sigma \in \SG_n$ be oscillating with connected intervals, $i\in [n-1]$ such that~$i$ and $i+1$ appear in different cycles of $\sigma$ and $\sigma' := s_i \sigma s_i$. Then $\sigma'$ is oscillating and has connected intervals.
\end{lem}
\begin{proof}
	We obtain $\sigma'$ from $\sigma$ by interchanging~$i$ and $i+1$ between two cycles in cycle notation. 
	It is easy to see that this does not affect the cycle standardization of the cycles in question. In addition, all other cycles of $\sigma'$ appear as cycles of $\sigma$. 
	Since $\sigma$ is oscillating with connected intervals, it follows that the standardization of each cycle of $\sigma'$ is oscillating with connected intervals. That is, $\sigma'$ is oscillating with connected intervals.
\end{proof}

We now come to the hooks.

\begin{exa}
	\label{thm:odd_hook_example_(3_1_1)}
	Let $\alpha = (3,1,1) \hdash 5$. The elements of $\Sigma_\alpha$ are
	\begin{align*}
	(1,5,2),(1,2,5),(1,5,3),(1,3,5),(1,5,4),(1,4,5).
	\end{align*}
	Note that~$1$ and~$5$ always appear in the cycle of length~$3$.
\end{exa}

Recall that we use \emph{type} as a short form for \emph{cycle type}.

\begin{defi}
	\label{def:hook_properties}
	Let $\alpha = (k,1^{n-k})\vDash_e n$ be a hook, $\sigma\in \SG_n$ of type $\alpha$, $m := \frac{k-1}2$ if~$k$ is odd and $m := \frac{k}2$ if~$k$ is even. We say that $\sigma$ satisfies the \emph{hook properties} if
	\begin{enumerate}
		\item $\sigma$ is oscillating,
		\item $\sigma$ has connected intervals,
		\item if $k>1$ then~$i$ and $n-i+1$ appear in the cycle of length~$k$ of $\sigma$ for all $i\in [m]$.
	\end{enumerate}
\end{defi}

The permutations from \Cref{thm:odd_hook_example_(3_1_1)} satisfy the hook properties.
The main result of this \namecref{sec:equivalence_classes:odd_hook} is to show that for an odd hook $\alpha$, the elements of $\Sigma_\alpha$ are characterized by the hook properties.
In \Cref{thm:characterization_of_Sigma_for_arbitrary_hook} of \Cref{sec:equivalence_classes:inductive_product} we will see that the same is true for even hooks.

\begin{exa}  
	\label{exa:odd_hook_properties}  
	\begin{wideenumerate}
	\item 
	Let $\sigma \in \SG_n$ be of type $(1^n)$.
	Then $\sigma = \id$ and $\sigma$ satisfies the hook properties.	
	Moreover, $\Sigma_{(1^n)} = \set{\sigma}$.	
		
	\item
	Let $\sigma \in \SG_n$ be of type $(n)$.
	That is, $\sigma$ is an~$n$-cycle.
	Then the third hook property is satisfied by $\sigma$ since all elements of $[n]$ appear in the only cycle of $\sigma$.
	Thus, $\sigma$ has the hook properties if and only if $\sigma$ is oscillating with connected intervals.
	By \Cref{thm:characterization_of_Sigma_(n)}, this is equivalent to $\sigma \in \Sigma_{(n)}$.
		
	\item
	\label{exa:odd_hook_properties_3,1,1}
	Let $\alpha = (3,1,1)\vDash n$. We want to determine all permutations in $\SG_n$ of type $\alpha$ that satisfy the hook properties. 
	Let $\sigma \in \SG_n$ be of type $\alpha$, $\sigma_1$ be the cycle of length~$3$ of $\sigma$ and $\mc O_1$ be the set of elements in $\sigma_1$.
	
	Since $\sigma_1$ is the only nontrivial cycle of $\sigma$, $\sigma$ is oscillating and has connected intervals if and only if $\tau:= \cst(\sigma_1)$ has these properties.
	The type of $\tau$ is $(3)$.
	By \Cref{thm:characterization_of_Sigma_(n)}, the oscillating permutations of type $(3)$ with connected intervals form $\Sigma_{(3)}$.
	From \Cref{tbl:Sigma_(n)} we read $\Sigma_{(3)} =\set {(1,3,2), (1,2,3)}$.
	
	Let
	\begin{align*}
	 M = \set{\set{1,5}\cup\set{j} \mid j \in [2,4]} = \set{ \set{1,2,5}, \set{1,3,5}, \set{1,4,5} }
	\end{align*}
	The third hook property is satisfied by $\sigma$ if and only if $ \mc O_1 \in M$.
	
	Therefore, $\sigma$ fulfills the hook properties if and only if there is a $\tau \in \Sigma_{(3)}$ and an $\mc O_1\in M$ such that we obtain $\sigma_1$ by writing $\mc O_1$ in a cycle such that the relative order of entries matches that one in $\tau$.
	For instance, from $\tau = (1,3,2)$ and $\mc O_1 = \set{1,4,5}$ we obtain $\sigma = (1,5,4)$.
	Going through all possibilities for $\tau$ and $\mc O_1$ we obtain the desired set of permutations. 
	These are the ones shown in \Cref{thm:odd_hook_example_(3_1_1)}.
\end{wideenumerate}
\end{exa}

For the proof of the characterization of $\Sigma_\alpha$ when $\alpha$ is an odd hook, we follow the same strategy as in in the case of compositions with one part from \Cref{sec:equivalence_classes:one_part}: For any odd hook $\alpha$ we show that $\sigma_\alpha$ satisfies the hook properties, $\approx$ is compatible with the hook properties and there is an algorithm that computes a sequence of $\approx$-equivalent permutations starting with $\sigma$ and ending up with $\sigma_\alpha$ for each permutation $\sigma$ of type $\alpha$ satisfying the hook properties.

\begin{lem} 
	\label{thm:sigma_odd_hook_has_odd_hook_properties}
	Let $\alpha\hdash n$ be an odd hook. Then the element in stair form $\sigma_\alpha\in \SG_n$ satisfies the hook properties.
\end{lem}
\begin{proof}
	Let $\alpha = \parts\alpha l = (k, 1^{n-k})\hdash n$ be an odd hook.
	If $k=1$ then $\sigma_\alpha$ is the identity which satisfies the hook properties.
	Assume $k>1$ and set $m := \frac{k-1}2$.
	By definition, the cycle of length~$k$ of $\sigma_\alpha$ is given by
	\begin{align*}
	\sigma_{\alpha_1} =	
	(1, n, 2, n-1, \dots, m, n-m+1, m+1).
	\end{align*}
	Hence, $\sigma_\alpha$ satisfies the third hook property.
	In order to show that $\sigma_\alpha$ is oscillating and has connected intervals, it suffices to consider $\sigma_{\alpha_1}$ because the other cycles of $\sigma_\alpha$ are trivial.
	From the description of $\sigma_{\alpha_1}$ we obtain its cycle standardization
	\begin{align*}
	\cst(\sigma_{\alpha_1}) = 
	(1, k, 2, k-1, \dots, m, k-m+1, m+1).
	\end{align*}
	That is, $\cst(\sigma_{\alpha_1})$ is the element in stair form $\sigma_{(k)}$ which is oscillating and has connected intervals by \Cref{thm:sigma_(n)_osc_and_c.I.}.
\end{proof}

Let $\alpha\vDash_e n$ be an odd hook and $\sigma \in \SG_n$ be of type $\alpha$ satisfying the hook properties.
In order to show $\sigma_\alpha \approx \sigma$ we will successively interchange elements~$i$ and $i+1$ in the cycle notation of $\sigma$. 
The next lemma considers the case where at least one of~$i$ and $i+1$ is a fixpoint of $\sigma$.

\begin{lem}
	\label{thm:hook_interchange_with_fixpoint}
	Let $\alpha = (k, 1^{n-k})\hdash n$ be an odd hook, $m:= \frac{k-1}2$ and  $\sigma \in \SG_n$ of type $\alpha$ satisfying the hook properties.
	Furthermore, assume that there are $i,i+1\in [m+1, n-m]$ such that~$i$ or $i+1$ is a fixpoint of $\sigma$.
	Then $s_i \sigma s_i \approx \sigma$ and $s_i \sigma s_i$ satisfies the hook properties.	
\end{lem}
\begin{proof}
	If both~$i$ and $i+1$ are fixpoints of $\sigma$ then $s_i \sigma s_i = \sigma$  and there is nothing to show.
	Therefore, we assume that either~$i$ or $i+1$ is not a fixpoint and call this element~$j$. By choice of~$i$ and $i+1$, $m < j < n-m+1$.
	Since $\sigma$ satisfies the hook properties, the cycle of length~$k$ of $\sigma$ consists of the elements $1,\dots, m, j, n-m+1, \dots, n$.

	First we show that $s_i\sigma s_i$ satisfies the hook properties.
	As $\sigma$ is oscillating with connected intervals and~$i$ and $i+1$ appear in different cycles of $\sigma$, \Cref{thm:equivalence_osc+c.I_and_interchange_between_cycles} yields that $s_i\sigma s_i$ is oscillating with connected intervals too.
	As we obtain $s_i\sigma s_i$ by interchanging~$i$ and $i+1$ in cycle notation of $\sigma$ and 
	\begin{align*}
	i,i+1\not \in \set{1,\dots, m, n-m+1, \dots, n}, 
	\end{align*}
	$s_i\sigma{s_i}$ satisfies the third hook property.

	In order to show $s_i \sigma s_i\approx \sigma$, we assume that $i+1$ is a fixpoint of $\sigma$ and~$i$ is not. The other case is proven analogously.
	Let $\tau  := \cst(\sigma)$ and $i' := \rnk_\sigma(i)$. Then  $i'= m+1 = \frac{k+1}2$ by the description of the cycle of length~$k$ from above. Since $\sigma$ is oscillating, $\tau$ is oscillating.
	Thus, \Cref{thm:oscillating_n-cycle_local} implies that there is an $\varepsilon \in \set{-1,1}$ such that 
	\begin{align*}
	\text{
		$\tau^\varepsilon(i') > m+1$ and $\tau^{-\varepsilon}(i') < m+1$.}
	\end{align*}	
	Now we use that $\tau^\delta(i') = \rnk_\sigma (\sigma^\delta(i))$ for 
	$\delta = -1,1$ and obtain that 
	\begin{align*}
	\text{$\sigma^\varepsilon(i) \geq n-m+1$ and $\sigma^{-\varepsilon}(i) \leq m$.}
	\end{align*}
	As $\sigma(i+1) = i+1 \in [m+2,n-m]$, it follows that
	\begin{align*}
		\text{$\sigma^\varepsilon(i) > \sigma^\varepsilon(i+1)$ and $\sigma^{-\varepsilon}(i) < \sigma^{-\varepsilon}(i+1).$} 
	\end{align*}
	Hence, \Cref{thm:length_and_conjugation} implies $\ell(s_i\sigma s_i) = \ell(\sigma)$.
	Therefore, $s_i\sigma s_i\approx \sigma$.
\end{proof}

The following \namecref{thm:equivalence_implies_odd_hook} shows that $\approx$ preserves the hook properties.
It is an analogue to \Cref{thm:equivalence_implies_oscillating_and_ci}.

\begin{lem}
	\label{thm:equivalence_implies_odd_hook}
	Given an odd hook $\alpha = (k, 1^{n-k})\hdash n$,  $\sigma \in \SG_n$ of type $\alpha$ satisfying the hook properties and $\sigma' := s_i \sigma s_i$ with $\sigma \approx \sigma'$, we have that also $\sigma'$ satisfies the hook properties.
\end{lem}

\begin{proof}
	We show that $\sigma'$ has the hook properties.
	If $k=1$ then $\sigma = \sigma' = \id$ so that $\sigma'$ satisfies the hook properties.
	Hence, assume $k>1$.
	Set  $m:= \frac{k-1}2$, $\tau := \cst(\sigma)$ and $\tau' := \cst(\sigma')$. We deal with three cases.
	
	First, assume that neither~$i$ nor $i+1$ is a fixpoint of $\sigma$.
	Then~$i$ and $i+1$ both appear in the cycle of length~$k$ of $\sigma$.
	Since $\sigma$ satisfies the hook properties, it is oscillating and has connected intervals. Therefore, \Cref{thm:equivalence_osc+c.I_and_interchange_within_cycle} yields that also $\sigma'$ has these properties. 
	The elements $1,\dots, m,n-m+1,\dots m$ all appear in the cycle of length~$k$ of $\sigma$ because $\sigma$ satisfies the hook properties. 
	Since we interchange two entries in this cycle to obtain $\sigma'$ from $\sigma$, all the elements also appear in the cycle of length~$k$ of $\sigma'$.
	
	Second, assume that $i+1$ is a fixpoint of $\sigma$ but~$i$ is not.
	Since $\sigma \approx \sigma'$, we have $\ell(\sigma) = \ell(\sigma')$  and by \Cref{thm:length_and_conjugation}
	\begin{align}
	\label{eq:hook_equivalence_preserves_hook_properties}
	\begin{aligned}
	\text{either } &\text{$\sigma(i) > i+1$ and $\sigma\inv(i) < i+1$} \\
	\text{or } &\text{$\sigma(i) < i+1$ and $\sigma\inv(i) > i+1$}
	\end{aligned}
	\end{align}
	where we used $\sigma(i+1) = i+1$.
	The elements of the cycle of length~$k$ of $\sigma$ are $1, \dots, m, j, n-m+1 , \dots, n$ where $j \in [m+1,n-m]$. 
	We now show that $i,i+1\in [m+1,n-m]$.
	
	As $i+1$ is a fixpoint, we have $i+1\leq n-m$ and it remains to show that $i \geq m+1$.
	Assume  that $i\leq m$ instead and set $i'  := \rnk_\sigma(i)$.
	Then $i' < \frac{k+1}2$.
	Since $\tau\in \SG_k$ is an oscillating $k$-cycle, \Cref{thm:oscillating_n-cycle_local} yields that 
	$\tau\inv(i'),\tau(i') \geq \frac{k+1}2$.
	Because $\rnk_\sigma(j) =\frac{k+1}2$, it follows that
	$\sigma\inv(i),\sigma(i) \geq j$.
	Moreover, $i+1$ being a fixpoint and $i \leq m$ imply that $i+1 < j$.
	Hence, $\sigma\inv(i),\sigma(i) > i+1$ which contradicts \Cref{eq:hook_equivalence_preserves_hook_properties}.
	
	Since $i,i+1\in [m+1,n-m]$ and $i+1$ is a fixpoint of $\sigma$,
	we can apply \Cref{thm:hook_interchange_with_fixpoint} which implies that $\sigma'$ satisfies the hook properties.
	
	In the same vein, one proves the remaining case where~$i$ is a fixpoint but $i+1$ is not.
\end{proof}

We now extend \Cref{thm:kim_algorithm_and_n-cycles} to the case of odd hooks. 
That is, we consider one step of the algorithm mentioned earlier.

\begin{lem} 
	\label{thm:kim_algorithm_and_odd_hooks}
	Let $\alpha = (k, 1^{n-k})\hdash n$ be an odd hook and  $\sigma \in \SG_n$ such that $\sigma$ is of type $\alpha$, $\sigma$ satisfies the hook properties and $\sigma\neq\sigma_\alpha$. 
	Then there exists a minimal integer $p$ such that $1\leq p\leq k-1$ and $\sigma^p(1) \neq \sigma_\alpha^p(1)$. Set $a := \sigma^p(1)$, $b := \sigma_\alpha^p(1)$ and
	\begin{align*}
	\sigma' := \begin{cases}
	s_{a-1} \sigma s_{a-1} &\myif a > b \\
	s_a \sigma s_a & \myif a < b. \\
	\end{cases}
	\end{align*}
	Then $\sigma' \approx \sigma$ and $\sigma'$ satisfies the hook properties.
\end{lem}

\begin{proof}
	Set $m := \frac{k-1}2$.
	If $\alpha = (1^n)$ then the only permutation of type $\alpha$ is the identity and there is nothing to show.
	If $\alpha = (n)$ then this is \Cref{thm:kim_algorithm_and_n-cycles}.
	Therefore, assume $1<k<n$. Since $\sigma$ satisfies the hook properties,~$1$ appears in the cycle of length~$k$ of $\sigma$.
	By definition, $\sigma_\alpha$ has the form
	\begin{align*}
	\sigma_\alpha =
	\begin{cases}
	(1, n ,2 , n-1 ,\dots, m+1) (n-m) (m+2) \cdots (\frac{n+3}2) (\frac{n+1}2) & \myif \text{$n$ is odd} \\
	(1, n ,2 , n-1 ,\dots, m+1) (n-m) (m+2) \cdots (\frac{n}2) (\frac{n}2+1) & \myif \text{$n$ is even}.
	\end{cases}	 
	\end{align*}
	In particular, $[m+2, n-m]$ is the set of fixpoints of $\sigma_\alpha$
	and~$1$ also appears in the cycle of length~$k$ of $\sigma_\alpha$. 
	Thus, from $\sigma \neq \sigma_\alpha$ it follows that there exists $p$ as claimed.
	In particular, we can define $a$, $b$ and $\sigma'$ as in the theorem.
	
	If~$n$ is odd, $k<n$ implies that $\frac{n+1}2$ is a fixpoint of $\sigma_\alpha$ and hence $b \neq \frac{n+1}2$. If~$n$ is even, we have $b \neq \frac{n+1}2$ anyway.
	Let $\tau := \cst(\sigma)$ and note that $\cst(\sigma_\alpha)$ is just the element in stair form $\sigma_{(k)}$.
	Moreover set $a' := \rnk_\sigma(a)$.

	Assume $b<\frac{n+1}2$. The proof for $b>\frac{n+1}2$ is similar and hence omitted. If $b<\frac{n+1}2$ then $b\leq m+1$ by the description of $\sigma_\alpha$ from above.
	The choice of~$p$ and $1 < b \leq m+1$ imply
	\begin{align*}
	\sigma^{-1}(a) = \sigma_\alpha^{-1}(b) = n-b+2 > m+1
	\end{align*}
	and
	\begin{align*}
	\set{1, 2, \dots, b-1} \subseteq \set{\sigma_\alpha^r(1) \mid r=0,\dots, p-1} = \set{\sigma^r(1) \mid r=0,\dots, p-1}.
	\end{align*}
	The last equality and $a\neq b$ imply $b < a$. Thus, we consider $\sigma' := s_{a-1}\sigma s_{a-1}$.
	From the hook properties, we obtain that the elements in the cycle of length~$k$ of $\sigma$ are $1,\dots, m,j,n-m+1,\dots n$ where $j\in [m+1, n-m]$.
	Thus, $\sigma^{-1}(a) >m+1$ implies $\tau\inv(a') > m+1$.
	But since $\sigma$ is oscillating, $\tau$ is oscillating and therefore \Cref{thm:oscillating_n-cycle_local} implies $a' \leq m+1$.
	From the description of the elements in the~$k$-cycle of $\sigma$, it now follows that $a\leq n-m$.
	
	To sum up, we have $b < a \leq n-m$ and  $\sigma' = s_{a-1}\sigma s_{a-1}$. 
	Now we have two cases depending on $a-1$.
	If $a-1$ is a fixpoint of $\sigma$ then because of $a\leq n-m$, we can apply \Cref{thm:hook_interchange_with_fixpoint} and obtain that $\sigma'\approx \sigma$ and $\sigma'$ satisfies the hook properties. 
	
	If $a-1$ is not a fixpoint of $\sigma$ then $\rnk_\sigma(a-1) = a'-1$.
	Moreover, interchanging $a-1$ and~$a$ in $\sigma$ does not affect the third part of the hook property.
	Therefore, we obtain from \Cref{thm:equivalence_and_interchange_within_cycle} that $\sigma'\approx \sigma$ and $\sigma'$ satisfies the hook properties if  $\tau' := s_{a'-1}\tau s_{a'-1}\approx \tau$ and $\tau'$ is oscillating with connected intervals. 
	By \Cref{thm:kim_algorithm_and_n-cycles}, $\tau'$ has these properties if $\tau^r(1) = \sigma_{(k)}^r(1)$ for $0\leq r \leq p-1$, $\tau^p(1) > \sigma^p_{(k)}(1)$ and $\tau^p(1) = a'$.
	This is what remains be shown.
	
	As $\sigma^r(1) = \sigma^r_\alpha(1)$ for $0\leq r \leq p-1$, we have the following equality of tuples
	\begin{align*}
	(\tau^0(1), \tau^1(1), \dots, \tau^{p-1}(1)) 
	&= (\rnk_\sigma(1), \rnk_\sigma(n) , \rnk_\sigma(2), \rnk_\sigma(n-1),\dots , \rnk_\sigma(n-b+2)) \\
	&= (1,k,2,k-1, \dots, k-b+2) \\
	&= (\sigma_{(k)}^0(1), \sigma_{(k)}^1(1), \dots, \sigma_{(k)}^{p-1}(1)). 
	\end{align*}
	Since the cycle of length~$k$ of $\sigma$ contains exactly one element of $[m+1,n-m]$, $a-1$ and~$a$ appear in this cycle and $a\leq n-m$, we have that $a\leq m+1$.
	Moreover, $1, \dots, m$ appear in the cycle of length~$k$ of $\sigma$ and $\sigma_\alpha$. Since $b< a \leq m+1$, this implies
	\begin{align*}
	\sigma^p_{(k)}(1) = \rnk_{\sigma_{\alpha}}(b) = b \text{ and } \tau^p(1)  = \rnk_{\sigma}(a) = a.
	\end{align*}
	In particular, $a' = \tau^p(1)$. Moreover,
	we have $b< a$ so that $\sigma^p_{(k)}(1) < \tau^p(1)$ as desired.
\end{proof}

We now come to the main result of this  \namecref{sec:equivalence_classes:odd_hook}.

\begin{thm}
	\label{thm:characterization_of_Sigma_for_odd_hook}
 	Let $\alpha\hdash n$ be an odd hook and $\sigma\in \SG_n$ of type $\alpha$. Then 
	$\sigma \in \Sigma_\alpha$ if and only if $\sigma$ satisfies the hook properties.
\end{thm}
\begin{proof}
	Let $\alpha = (k, 1^{n-k})\hdash n$ be an odd hook and $\sigma_\alpha$ be the element in stair form.
	The proof is analogous to the one of \Cref{thm:characterization_of_Sigma_(n)}.
	First, $\sigma_\alpha$ satisfies the hook properties by \Cref{thm:sigma_odd_hook_has_odd_hook_properties}.
	Let $\sigma \in \SG_n$.
	
	For the direction from left to right
	assume that $\sigma \in \Sigma_\alpha$.
	Then $\sigma \approx \sigma_\alpha$.
	From the definition of $\approx$ and \Cref{thm:equivalence_implies_odd_hook}
	it follows that $\approx$ transfers the hook properties from $\sigma_\alpha$ to $\sigma$.
	
	For the converse direction, assume that $\sigma$ satisfies the hook properties.
	By using \Cref{thm:kim_algorithm_and_odd_hooks} iteratively, we obtain a sequence of $\approx$-equivalent permutations starting with $\sigma$ and ending in $\sigma_\alpha$.
	Hence $\sigma \in \Sigma_\alpha$.
\end{proof}

We continue with a rule for the construction of $\Sigma_{(k,1^{n-k})}$ from $\Sigma_{(k)}$ in the case where~$k$ is odd and $k\geq 3$.
The rule can be sketched as follows.
Given a $\tau \in \Sigma_{(k)}$ we can choose a subset of $[n]$ of size $k$ in accordance with the third hook property.
Arranging the elements of this subset in a cycle of length $k$ such that its cycle standardization is  $\tau$ (and letting the other elements of $[n]$ be fixpoints) then results in an element of $\Sigma_{(k,1^{n-k})}$.
See Part~\ref{exa:odd_hook_properties_3,1,1} of \Cref{exa:odd_hook_properties}  for an illustration.

\begin{cor}
	\label{thm:odd_hook_bijection}
	Let $\alpha = (k,1^{n-k})\vDash_e n$ be an odd hook with $k\geq 3$.
	Set $m := \frac{k-1}2$. 
	For $\tau \in \Sigma_{(k)}$ and $j \in [m+1, n-m]$ define $\varphi(\tau,j)$ to be the element $\sigma \in \SG_n$ of type $\alpha$ such that $\cst(\sigma) = \tau$ and the entries in the cycle of length~$k$ of $\sigma$ are $1,\dots, m, j ,n-m+1,\dots, n$. 
	Then
	\begin{align*}
	\varphi \colon \Sigma_{(k)} \times [m+1, n-m] \to \Sigma_\alpha, \quad (\tau,j) \mapsto \varphi(\tau,j)
	\end{align*}
	is a bijection.
\end{cor}

\begin{proof}
	
	Given a $\tau \in \Sigma_{(k)}$ and a $j\in [m+1, n-m]$ there is only one way (up to cyclic shift) to write the elements $1,2,\dots, m, j , n-m+1, \dots, n$ in a cycle of length~$k$ such that the standardization of the corresponding~$k$-cycle in $\SG_n$ is $\tau$.
	This~$k$-cycle is $\varphi(\tau,j)$.
	By construction, $\varphi(\tau,j)$ satisfies the hook properties.
	Hence, \Cref{thm:characterization_of_Sigma_for_odd_hook} yields $\varphi(\tau,j) \in \Sigma_\alpha$.
	That is, $\varphi$ is well defined.
	
	Let $\sigma\in \Sigma_\alpha$.
	Then by \Cref{thm:characterization_of_Sigma_for_odd_hook}, $\sigma$ satisfies the hook properties.
	The third hook property yields that there is a unique $j\in [m+1, n-m]$ such that the elements in the cycle of length~$k$ of $\sigma$ are $1,2,\dots, m, j , n-m+1, \dots, n$.
	From the first two hook properties it follows that $\tau:= \cst(\sigma)$ is oscillating and has connected intervals.
	Thus, $\tau \in \Sigma_{(k)}$ by \Cref{thm:characterization_of_Sigma_(n)}.
	By definition of $\varphi$, the cycles of length~$k$ of $\varphi(\tau,j)$ and $\sigma$ contain the same elements.
	Moreover, they have the same cycle standardization $\tau$.
	Consequently, $\varphi(\tau,j) = \sigma$.
	That is, $\varphi$ is surjective.
	Since $\tau$ and~$j$ uniquely depend on $\sigma$, $\varphi$ is also injective.
\end{proof}

In the last result of the \namecref{sec:equivalence_classes:odd_hook} we determine the cardinality of $\Sigma_\alpha$ for each odd hook $\alpha$.

\begin{cor}
	\label{thm:odd_hook_cardinality_of_Sigma_alpha}
	If $\alpha= (k,1^{n-k})\hdash n$ is an odd hook then
	\begin{align*}
	|\Sigma_{\alpha}| = 
	\begin{cases}
	1 & \myif \text{$k=1$} \\	
	2 (n-k+1) 3^{\frac{k-3}2} & \myif \text{$k\geq 3$}. \\
	\end{cases}
	\end{align*}
\end{cor}

\begin{proof}
	Let $\sigma \in \Sigma_\alpha$.
	If $k=1$ then $\Sigma_\alpha = \set{1}$.
	Now suppose that $k\geq 3$ and set $m := \frac{k-1}2$.
	The cardinality of $[m+1,n-m]$ is $n-k+1$.
	Hence, \Cref{thm:odd_hook_bijection} yields that $|\Sigma_\alpha| = (n-k+1) |\Sigma_{(k)}|.$
	In addition, we have $|\Sigma_{(k)}| = 	2 \cdot 3^{\frac{k-3}2}$
	from  \Cref{thm:sizes_of_Sigma_n}.
\end{proof}

\subsection{The inductive product}
\label{sec:equivalence_classes:inductive_product}

In this \namecref{sec:equivalence_classes:inductive_product}
we define the inductive product $\iprod$ and use it to obtain in \Cref{thm:indcutive_product_bijection} a recursion the rule for $\Sigma_{(\alpha_1,\dots, \alpha_l)}$ in the case where $\alpha_1$ is even.
This leads to a description of $\Sigma_\alpha$ for all maximal compositions $\alpha$ whose odd parts form a hook (see \Cref{thm:inductive_product_remark_reduction_to_odd_partitions}).
As a consequence, we show in \Cref{thm:characterization_of_Sigma_for_arbitrary_hook} that $\Sigma_\alpha$ is characterized by the hook properties if $\alpha$ is an even hook.

Recall that we write $\gamma \vDash_0 n$ if $\gamma$ is a weak composition of $n$, that is, a finite sequence of nonnegative  integers that sum up to $n$.

\begin{defi}
\label{def:inductive_product}
Let $(n_1,n_2) \vDash_0 n$.
The \emph{inductive product} on  $\SG_{n_1} \times \SG_{n_2}$ is the binary operator
\begin{align*}
	\iprod \colon \SG_{n_1} \times \SG_{n_2} &\to \SG_n \\
	(\sigma_1, \sigma_2) &\mapsto\sigma_1 \iprod \sigma_2
\end{align*}
where $\sigma_1 \iprod \sigma_2$ is the element of $\SG_n$ whose cycles are the cycles of $\sigma_1$ and $\sigma_2$ altered as follows:
\begin{enumerate}
	\item in the cycles of ${\sigma_1}$, add ${n_2}$ to each entry ${>k}$,
	\item in the cycles of ${\sigma_2}$, add ${k}$ to each entry
\end{enumerate}
where $k := \ceil{\frac{n_1}2}$.
\end{defi}

For two sets $X_1\subseteq \SG_{n_1}$ and $X_2 \subseteq \SG_{n_2}$ we define
\begin{align*}
X_1 \iprod X_2 := \set{\sigma_1 \iprod \sigma_2 \mid \sigma_1\in X_1, \sigma_2 \in X_2 }.
\end{align*}

We will see in \Cref{thm:inductive_product_cycles} that the inductive product is well-defined.

\begin{exa} 
	\label{exa:inductive_product}
	\begin{enumerate}
		\item
		Let $\emptyset\in \SG_0$ be the empty function and $\sigma \in \SG_n$. Then 
		\begin{align*}
		\emptyset\iprod \sigma = \sigma \iprod \emptyset = \sigma.
		\end{align*}
		
		\item
			Consider $n_1 = 6$, $n_2 = 4$, $n= 10$ and the elements in stair form $\sigma_{(6)} \in \SG_{n_1}$ and $\sigma_{(3,1)} \in \SG_{n_2}$.
			 Then $k = 3$ and
			\begin{align*}
				\sigma_{(6)} \iprod \sigma_{(3,1)} 
				&=  (1, 6, 2, 5, 3, 4) \iprod  (1, 4, 2)(3) \\
				&= (1, 6 + 4 , 2,  5 + 4, 3, 4 + 4)( 1  +  3, 4  +  3, 2  +  3)(3 +  3)  \\
				&= (1, 10 , 2 ,  9, 3,  8)( 4,  7,  5)( 6 ).
			\end{align*}
		
			\item
		Consider $n_1 = 5$, $n_2 = 4$ and the elements in stair form $\sigma_{(5)} = ( 1, 5, 2, 4, 3)\in \SG_{n_1}$ and  $\sigma_{(3,1)} = (1,4,2)(3)\in \SG_{n_2}$.
		Then $\sigma_{(3,1)}^{w_0} = ( 1, 3,4)(2)$ where $w_0 = (1,4)(2,3)$ is the longest element of $\SG_4$. 
		We have
		 $k = 3$ and
		\begin{align*}
		\sigma_{(5)} \iprod \sigma^{w_0}_{(3,1)} 
		&= ( 1, 5+4, 2, 4+4, 3)(1+3, 3+3, 4+3)(2+3) \\
		&= 
		(1, 9 , 2 , 8, 3)(7, 4 , 6)(5).
		\end{align*}
\end{enumerate}
Note that in Parts (2) and (3) we obtain the elements in stair form $\sigma_{(6,3,1)}$ and $\sigma_{(5,3,1)}$, respectively.
\end{exa}

In order to work with the inductive product, it is convenient to describe it more formally.
To this end we introduce the following notation which we will use throughout the \namecref{sec:equivalence_classes:inductive_product}.

\begin{notn}
\label{not:inductive_product}
Let $n\geq 0$,
$(n_1,n_2) \vDash_0 n$,
 $k := \ceil{\frac{n_1}2}$,
\begin{align*}
N_1 := [k] \cup [k+n_2+1, n] \quad \text{and} \quad N_2 := [k+1, k+n_2].
\end{align*}

We have that $|N_1| = n_1, |N_2| = n_2$, $N_1$ and $N_2$ are disjoint and $N_1\cup N_2 = [n]$.
Note that $[0] = [1,0] = \emptyset$.
Define the bijections $\varphi_{1}\colon [n_1] 
\to N_1$ and $\varphi_2 \colon [n_2]  \to N_2$ by
\begin{align*}
\varphi_1(i) :=
	\begin{cases}
		i &\myif i\leq k \\
		i+n_2 & \myif i> k
	\end{cases} 
	\quad \text{and} \quad
	\varphi_2(i)  := i + k.	
\end{align*}
The bijections $\varphi_1$ and $\varphi_2$ formalize the alteration of the cycles of $\sigma_1$ and $\sigma_2$ in \Cref{def:inductive_product}, respectively.
Their inverses are given by
\begin{align*}
\varphi\inv_1(i) := \begin{cases}
i & \myif i\leq k \\
i-n_2 & \myif i>k
\end{cases}
\quad \text{and} \quad
\varphi\inv_2(i) := i - k.
\end{align*}

For $i=1,2$ and $\sigma_i\in \SG_{n_i}$, write $\sigma_i^{\varphi_i} := \varphi_i \circ \sigma_i \circ \varphi_i\inv$.
Then $\sigma_i^{\varphi_i}\in \SG(N_i)$ and $\sigma_i^{\varphi_i}$ can naturally be identified with the element of $\SG_n$ that acts on $N_i$ as $\sigma_i^{\varphi_i}$ and fixes all elements of $[n] \setminus N_i$.

\end{notn}

We will see in \Cref{thm:inductive_product_cycles} that we obtain $\sigma_i^{\varphi_i}$ by applying $\varphi_i$ on each entry in of $\sigma_i$ in cycle notation.

\begin{exa}
	\label{exa:inductive_product_functions}
	Let $n_1 = 6$ and $n_2 = 4$ and consider the elements in stair form
	\begin{center} $\sigma_1 := \sigma_{(6)} = ( 1, 6, 2, 5, 3, 4) \in \SG_6$ \quad and \quad $\sigma_2 := \sigma_{(3,1)} = ( 1, 4, 2)(3) \in \SG_4$. 
	\end{center}
	 Then $k = 3$ and
	\begin{align*}
		\sigma_1^{\varphi_1} &= (1, 6 + 4 , 2,  5 + 4, 3, 4 + 4) = (1, 10 , 2 ,  9, 3,  8),\\ \sigma_2^{\varphi_2} &= ( 1  +  3, 4  +  3, 2  +  3)(3 +  3) = ( 4,  7,  5)( 6 ).
	\end{align*}
	Thus, from \Cref{exa:inductive_product} it follows that $\sigma_1 \iprod \sigma_2 = \sigma_1^{\varphi_1}\sigma_2^{\varphi_2}$.
	The next \namecref{thm:inductive_product_cycles} states that this is true in general.
\end{exa}

We now come to the more formal description of the inductive product.
\begin{lem}
	\label{thm:inductive_product_cycles}
	Let $\sigma_r \in \SG_{n_r}$ with cycle decomposition $\sigma_r = \sigma_{r,1}\sigma_{r,2}\cdots \sigma_{r,p_r}$ for $r= 1,2$.
	\begin{enumerate}
	
		\item
			We have 
			\begin{align*}
				\sigma_1 \iprod \sigma_2 = \sigma_1^{\varphi_1} \sigma_2^{\varphi_2}.
			\end{align*}
					
		\item 
		Let $r\in \set{1,2}$ and $\sigma_{r,j} = (c_1,\dots, c_t)$ be a cycle of $\sigma_r$.
		Then
		\begin{align*}
			\sigma_{r,j}^{\varphi_r} = (\varphi_r(c_1),\dots, \varphi_r(c_t)).
		\end{align*}
					
		\item 
		The decomposition of $\sigma_1 \iprod \sigma_2$ in disjoint cycles is given by 
			\begin{align*}
				\sigma_1 \iprod \sigma_2 = \sigma^{\varphi_1}_{1,1} \cdots \sigma^{\varphi_1}_{1,p_1} 
				\cdot
				\sigma^{\varphi_2}_{2,1} \cdots \sigma^{\varphi_2}_{2,p_2}. 
			\end{align*}	
	\end{enumerate}
\end{lem}

\begin{proof}
	Set $\sigma := \sigma_1 \iprod \sigma_2$ and $\sigma' := \sigma_1^{\varphi_1}\sigma_2^{\varphi_2}$.
	It will turn out that $\sigma = \sigma'$.

	We first show Part~(2).
	Let $r\in \set{1,2}$, $\xi$ be a cycle of $\sigma_r$ and $i\in [n_r]$.
	Then
	\begin{align*}
		\xi^{\varphi_r}(\varphi_r(i)) = (\varphi_r \circ \xi \circ \varphi_r\inv \circ \varphi_r )(i) = \varphi_r(\xi(i)).
	\end{align*}
	Hence, if $\xi = (c_1,\dots, c_t) \in \SG_{n_r}$ then $\xi^{\varphi_r} =  (\varphi_r(c_1),\dots, \varphi_r(c_t)) \in \SG(N_r)$.

	We continue with showing Part~(3) for $\sigma'$.
	For $r =1,2$ we have
	\begin{align*}
		\sigma_r^{\varphi_r} 
		&= \varphi_r \circ \sigma_r \circ \varphi_r\inv \\
		&= \varphi_r \circ \sigma_{r,1}\cdots \sigma_{r,p_r} \circ \varphi_r\inv \\
		&= (\varphi_r \circ \sigma_{r,1} \circ \varphi_r\inv) \cdots (\varphi_r \circ \sigma_{r,p_r} \circ \varphi_r\inv) \\
		&= \sigma_{r,1}^{\varphi_r} \cdots \sigma_{r,p_r}^{\varphi_r}.
	\end{align*}
	Thus,
	\begin{align}
		\label{eq:cycle_decomposition_sigma}
		\sigma'
		= \sigma_{1,1}^{\varphi_1} \cdots \sigma_{1,p_1}^{\varphi_1} \sigma_{2,1}^{\varphi_2} \cdots \sigma_{1,p_2}^{\varphi_2}.
	\end{align}
	The cycles in this decomposition are given by Part~(1). 
	As $\varphi_1$ and $\varphi_2$ are bijections with disjoint images, the cycles are disjoint.
	
	Lastly, we show $\sigma = \sigma'$.
	From \Cref{eq:cycle_decomposition_sigma}, Part~(2) and the definition of $\varphi_1$ and $\varphi_2$ it follows that we obtain the cycles of $\sigma'$ by altering the cycles of $\sigma_1$ and $\sigma_2$ as described in \Cref{def:inductive_product}.
	Hence, $\sigma = \sigma'$.
\end{proof}

\begin{cor}
	\label{thm:inductive_product_orbits}
	Let	$\sigma_1\in \SG_{n_1},
	\sigma_2\in \SG_{n_2}$ and
	$\sigma := \sigma_1\iprod \sigma_2$.
	Then
	\begin{align*}
		P(\sigma) = \varphi_1(P(\sigma_1)) \cup \varphi_2(P(\sigma_2)).
	\end{align*}
\end{cor}

We continue with basic properties of the inductive product.

\begin{lem}
	\label{thm:inductive_product_two_domains}
	Let $\sigma_1\in \SG_{n_1}, \sigma_2\in \SG_{n_2}$ and $\sigma := \sigma_1\iprod \sigma_2$. Then for all $i\in [n]$
	\begin{align*}
	\sigma(i) = 
	\begin{cases}
	\sigma_1^{\varphi_1}(i) &\myif i\in N_1 \\
	\sigma_2^{\varphi_2}(i) &\myif i\in N_2.
	\end{cases}
	\end{align*}
\end{lem}

\begin{proof}
	By \Cref{thm:inductive_product_cycles}, $\sigma = \sigma_1^{\varphi_1}\sigma_2^{\varphi_2}$.
	If $n_1= 0$ or $n_2= 0$ the claim is trivially true.
	Thus, suppose $n_1,n_2\geq 1$ and let $i\in [n]$.
	Consider $\sigma_1^{\varphi_1}$ and $\sigma_2^{\varphi_2}$ as elements of $\SG_n$.
	Since $\set{N_1,N_2}$ is a partition of $[n]$ there is exactly one $r\in \set{1,2}$ such that $i\in N_r$.
	We have that $\sigma_r^{\varphi_r}(N_r) = N_r$ and that $\sigma_{2-r+1}^{\varphi_{2-r+1}}$ fixes each element of $N_r$. Hence,
	\begin{align*}
	\sigma(i) = \sigma_1^{\varphi_1}\sigma_2^{\varphi_2}(i) = \sigma_r^{\varphi_r}(i). &\qedhere
	\end{align*}
\end{proof}

We now determine the image of the inductive product and show that it is injective.

\begin{lem}
	\label{thm:inductive_product_image_and_injectivity}
	Let $(n_1,n_2) \vDash_0 n$.
	\begin{enumerate}
		\item
		The image of  $\SG_{n_1} \times \SG_{n_2}$ under $\iprod$ is given by 
		\begin{align*}
		\SG_{n_1} \iprod \SG_{n_2} = \set{\sigma \in \SG_n \mid \sigma(N_i) = N_i \text{ for } i = 1,2}.
		\end{align*}
		\item 
		The inductive product on $\SG_{n_1} \times \SG_{n_2}$ is injective.
	\end{enumerate}
\end{lem}

\begin{proof}
	\begin{wideenumerate}
		\item Set $Y := \set{\sigma \in \SG_n \mid \sigma(N_i) = N_i \text{ for } i = 1,2}$.
		
		We first show $\SG_{n_1}\iprod\SG_{n_2} \subseteq Y$.
		Let $\sigma \in \SG_{n_1}\iprod\SG_{n_2}$.
		Then there are $\sigma_i \in \SG_{n_i}$ for $i = 1,2$ such that $\sigma = \sigma_1\iprod \sigma_2$. 
		By \Cref{thm:inductive_product_two_domains} we have $\sigma(N_i) = \sigma^{\varphi_i}(N_i) = N_i$ for $i = 1,2$. Hence, $\sigma \in Y$.
		
		Now we show $Y \subseteq \SG_{n_1}\iprod\SG_{n_2}$. Let $\sigma \in Y$.
		For $i=1,2$ set $\tilde \sigma_i = \sigma |_{N_i}$ (the restriction to $N_i$).
		Consider $i \in \set{1,2}$.
		Since $\sigma\in Y$, $\tilde \sigma_i(N_i) = N_i$ and thus $\tilde \sigma_i \in \SG(N_i)$.
		Therefore, $\sigma_i := \varphi\inv_i \circ \tilde\sigma_i \circ \varphi_i$ is an element of $\SG_{n_i}$.
		Moreover, $\sigma_i^{\varphi_i}$ considered as an element of $\SG_n$ leaves each element of $N_{2-i+1}$ fixed. 
		Hence, we have
		\begin{align*}
		 (\sigma_1 \iprod \sigma_2)|_{N_i} = \sigma^{\varphi_1}_1\sigma^{\varphi_2}_2|_{N_i} = \sigma^{\varphi_i}_i|_{N_i} = \tilde \sigma_i |_{N_i} = \sigma |_{N_i}.
		\end{align*}
		Consequently, $\sigma = \sigma_1 \iprod \sigma_2$.
		
		\item Since $|N_i| = n_i$ for $i = 1,2$, the cardinality of~$Y$ is $n_1!n_2!$. This is also the cardinality of $\SG_{n_1} \times \SG_{n_2}$. As the image of $\SG_{n_1} \times \SG_{n_2}$ under $\iprod$ is~$Y$, it follows that $\iprod$ is injective. \qedhere
	\end{wideenumerate}
\end{proof}

Recall that for $\alpha \vDash_e n$, each element of $\Sigma_\alpha$ has the property that its length is maximal in its conjugacy class.
We want to use this property to prove our main result.

Consider $\sigma = \sigma_1 \iprod \sigma_2$ such that $\sigma_1$ has type $(n_1)$.
We seek a formula for $\ell(\sigma)$ depending on $\sigma_1$ and $\sigma_2$.
We are particularly interested in the case where the $n_1$-cycle $\sigma_1$ is oscillating.

Given $\sigma \in \SG_n$ let $\Inv(\sigma) := \set{(i,j) \mid 1\leq i < j \leq n, \sigma(i) > \sigma(j) } $ be the \emph{set of inversions} of $\sigma$.
Then $\ell(\sigma) = |\Inv(\sigma)|$ by \cite[Proposition 1.5.2]{Bjorner2006}.

\begin{lem} 
	\label{thm:inductive_product_and_length}
	Let $\sigma_1\in \SG_{n_1}$ be an $n_1$-cycle, $\sigma_2\in \SG_{n_2}$, $\sigma := \sigma_1 \iprod \sigma_2$,
	\begin{align*}
		P &:= \set{i\in [k] \mid \sigma_1(i) > k}, \\
		Q &:= \set{i\in [k+1,n_1] \mid \sigma_1(i) \leq k},
	\end{align*}
	$p := |P|$ and $q := |Q|$.
Then we have
\begin{align*}
	\ell(\sigma) = \ell(\sigma_1) + \ell(\sigma_2) + (p+q)n_2.
\end{align*}
Moreover,
\begin{enumerate}
	\item $p,q \leq \floor{\frac{n_1}2}$,
	\item if $\sigma_1$ is oscillating, then $p = q = \floor*{\frac{n_1}2}$.
\end{enumerate}
\end{lem}

\begin{proof}
	Let $i,j\in [n]$ and $m:= \floor*{\frac{n_1}{2}}$.
	We distinguish three types of pairs $(i,j)$ and count the number of inversions of $\sigma$ type by type.
	\begin{enumerate}[label=\bf{Type \arabic*.}, nosep, wide]
		\item 
		There is an $r\in \set{1,2}$ such that $i,j\in N_r$.
		In this case let $t\in \set{i,j}$ and set $t' := \varphi_r\inv(t)$.
		Then $t' \in [n_r]$. 
		From \Cref{thm:inductive_product_two_domains} we obtain 
		\begin{align*}
		\sigma(t)= \varphi_r(\sigma_r(t')).
		\end{align*}
		In addition, we have 
		\begin{align*}
			\varphi_r(\sigma_r(i')) > \varphi_r(\sigma_r(j')) \iff \sigma_r(i') > \sigma_r(j')
		\end{align*}
		since $\varphi_r$ is a stricly increasing function.
		As $\varphi_r\inv$ is stricly increasing as well, we also have that
		\begin{align*}
		 i<j \iff i' < j'.
		 \end{align*}
		Hence,
		\begin{align*}
		(i,j) \in \Inv(\sigma)
		&\iff \text{$i < j$ and $\sigma(i) > \sigma(j)$} \\
		&\iff \text{$i' < j'$ and $\varphi_r(\sigma_r(i')) > \varphi_r(\sigma_r(j'))$} \\
		&\iff \text{$i' < j'$ and $\sigma_r(i') > \sigma_r(j')$} \\
		&\iff (i',j') \in \Inv(\sigma_r).
		\end{align*}
		Thus, the number of inversions of Type~1 is
		\begin{align*}
		|\Inv(\sigma_1)| + |\Inv(\sigma_2)| = \ell(\sigma_1) + \ell(\sigma_2).
		\end{align*}

		\item We have $i \in N_1$, $j\in N_2$ and $i<j$.
		Assume that $(i,j)$ is of this type and 
		recall that $N_1 = [k] \cup [k + n_2 + 1, n]$ and $N_2 = [k+1, k+n_2]$ where $k = \ceil{\frac{n}2}$.
		Since $i<j$, we have  $i\leq k$ which in particular means that $\varphi_1^{-1}(i) = i$.
		As $\sigma(j) \in N_2$,  $k+1 \leq \sigma(j) \leq k + n_2$.
		Moreover, $\sigma(i) = \sigma_1^{\varphi_1}(i)$ by \Cref{thm:inductive_product_two_domains}.
		Consequently,
		\begin{align*}
			\sigma(i) = \sigma_1^{\varphi_1}(i) = \varphi_1(\sigma_1(i)) = 
			\begin{cases}
				\sigma_1(i) < \sigma(j) &\myif \sigma_1(i)  \leq k \\
				\sigma_1(i) + n_2 > \sigma(j) &\myif \sigma_1(i) > k.
			\end{cases}
		\end{align*}
		Therefore,
		\begin{align*}
			 (i,j) \in \Inv(\sigma) \iff \sigma_1(i) > k.
		\end{align*}
		Hence, the number of inversions of Type~2 is the cardinality of the set $P \times N_2$.
		Thus, we have $pn_2$ inversions of Type~2.

		\item We have $i\in N_2$, $j \in N_1$ and $i<j$. Let $(i,j)$ be of Type~$3$.
		Then from $i<j$ we obtain $j\geq k+n_2+1$.
		In particular this type can only occur if $n_1 > 1$ because otherwise $n =  1 + n_2< j$.
		
		Since $i \in N_2$, also $\sigma(i)\in N_2$.
		That is, $k+1 \leq \sigma(i) \leq k+n_2$.	
		Moreover, from $i < j$ and $i\in N_2$ it follows that $j \geq k + n_2+1$.
		Thus, 
		\begin{align*}
			j':=\varphi_1\inv(j) = j - n_2
		\end{align*}
		and  $j' \in [k+1,n_1]$.
		Hence,
		\begin{align*}
			\sigma(j) = \sigma_1^{\varphi_1}(j) = \varphi_1(\sigma_1(j')) =
				\begin{cases}
			\sigma_1(j') < \sigma(i) &\myif \sigma_1(j')  \leq k \\
			\sigma_1(j') + n_2 > \sigma(i) &\myif \sigma_1(j') > k.
			\end{cases}
		\end{align*}
		That is, 
		\begin{align*}
					 (i,j) \in \Inv(\sigma) \iff \sigma_1(j') \leq k
					 \iff j' \in Q
					 \iff j \in \varphi_1(Q)
		\end{align*}
		where we use that $j'\in [k+1, n_1]$ for the second equivalence.
		Consequently, the set of inversion of Type~3 is the set $N_2 \times \varphi_1(Q)$.
		Since $\varphi_1$ is a bijection, it follows that there are exactly $qn_2$ inversions of this type.
	\end{enumerate} 

	Summing up the number of inversions of each type, we obtain the formula for the length of $\sigma$.
	
	We now prove $(1)$ and $(2)$.
	\begin{wideenumerate}
		\item 
		By definition, $\sigma_1(P) \subseteq [k+1,n_1]$ and $Q\subseteq [k+1,n_1]$.
		The cardinality of  $[k+1,n_1]$ is $\floor*{\frac{n_1}{2}}$.
		Therefore, $p,q\leq \floor*{\frac{n_1}{2}}$.
		 \item
		 Assume that $\sigma_1$ is oscillating.
		 Suppose first that $n$ is even.
		 Then $k = \frac{n_1}2$.
		 Because $\sigma_1$ is oscillating, we obtain that
		 \begin{align*}
		 	\sigma_1([k]) = [k+1,n_1] \quad \text{and} \quad \sigma_1([k+1,n_1]) = [k]
		 \end{align*}
		 from \Cref{def:oscillating_n-cycle} and  \Cref{thm:oscillating_n_cycle_complements}.
		 Hence, $p = q = k = \floor{\frac{n_1}2}$.
		 
		Suppose now that $n$ is odd.
		Then $k = \frac{n_1+1}2$.
		Since $\sigma_1$ is oscillating, \Cref{def:oscillating_n-cycle} and  \Cref{thm:oscillating_n_cycle_complements} yield that there is an $m\in \set{k-1,k}$ such that
		\begin{align*}
				\sigma_1([m])  = [n_1-m+1, n_1]   \quad\text{and}\quad   \sigma_1([m+1, n_1])  = [n_1-m].
		\end{align*}
		It is not hard to see that this implies $p = q = k-1 = \floor{\frac{n_1}2}$.	 
		 \qedhere		 
	\end{wideenumerate}
\end{proof}

We have seen in \Cref{exa:inductive_product} that the elements in stair form $\sigma_{(5,3)}$ and $\sigma_{(6,3)}$ can be decomposed as
\begin{align*}
\sigma_{(5,3)} = \sigma_{(5)} \iprod \sigma_{(3)}^{w_0} \quad \text{and} \quad \sigma_{(6,3)} = \sigma_{(6)} \iprod \sigma_{(3)}
\end{align*}
where $w_0$ is the longest element of $\SG_3$.
We want to show that these are special cases of a general rule for decomposing the element in stair form $\sigma_\alpha$.
Before we state the rule, we compare the sequences 
used to define the element in stair form in \Cref{thm:element_in_stair_form} for  compositions of~$n$, $n_1$ and $n_2$.

\begin{lem} 
		\label{thm:inductive_product_and_x_sequence}
	For $m \in \N_0$ let $x^{(m)}$ be the sequence $(x^{(m)}_1,\dots, x^{(m)}_m)$ given by $x^{(m)}_{2i-1} = i$ and $x^{(m)}_{2i} = m-i+1$.
Set $x:=x^{(n)}$, $y := x^{(n_1)}$ and  $z := x^{(n_2)}.$
	\begin{enumerate}
		\item We have $\varphi_1(y_i) = x_i$ for all $i\in [n_1]$.
		\item If $n_1$ is even then 	$\varphi_2(z_i) = x_{i+n_1}$ for all $i\in [n_2]$.
		\item If $n_1$ is odd then 	$\varphi_2(w_0(z_i)) = x_{i+n_1}$ for all $i\in [n_2]$ where $w_0$ is the longest element of $\SG_{n_2}$.
	\end{enumerate}
\end{lem}

\begin{proof}
	Recall that $k = \ceil*{\frac{n_1}2}$ and $(n_1,n_2)\vDash_0 n$ by \Cref{not:inductive_product}. 
	Let $i\in \N$.
	We mainly do straight forward calculations.
	\begin{wideenumerate}
		\item
			 Assume $2i-1 \in [n_1]$. Then $i\leq k$ and thus $\varphi_1(i) = i$. Consequently,
			\begin{align*}
				\varphi_1(y_{2i-1}) = \varphi_1(i) = i = x_{2i-1}.
			\end{align*}
			Now, assume $2i \in [n_1]$. Then
			\begin{align*}
				n_1-i+1 
				=\ceil*{n_1-i+1} 
				&\geq \ceil*{n_1-\frac{n_1}2+1} \\
				&= \ceil*{\frac{n_1}2+1} 
				=\ceil*{\frac{n_1}2} + 1 = k + 1,
			\end{align*} 
			\ie $\varphi_1(n_1-i+1) = n_1+n_2-i+1$.
			Therefore,
			\begin{align*}
			\varphi_1(y_{2i}) = \varphi_1(n_1-i+1) = n_1+n_2-i+1 = n-i+1 = x_{2i}.
			\end{align*}
		\item 
			Assume that $n_1$ is even. Then $n_1 = 2k$. If $2i-1\in [n_2]$ then we have 
			\begin{align*}
			2(k+i)-1 = n_1 + 2i-1\leq n_1+n_2 = n.
			\end{align*}
			Thus,
			\begin{align*}
				\varphi_2(z_{2i-1}) = \varphi_2(i) = k + i = x_{2(k+i)-1} = x_{2i -1+n_1}.
			\end{align*}
			Suppose $2i\in [n_2]$. Then $2(k+i) = n_1 + 2i \leq n$ and
			\begin{align*}
				\varphi_2(z_{2i}) = k+n_2-i+1 &= (n-2k-n_2) + k + n_2 -i+1 \\
				&= n-k-i+1 \\
				&= x_{2(k+i)} = x_{2i + n_1}.
			\end{align*}
		\item
		Assume that $n_1$ is odd.
		In this case $n_1 = 2k-1$.
		Let $w_0$ be the longest element of $\SG_{n_2}$.
		We have $w_0(j) = n_2-j+1$ for all $j\in [n_2]$.
		If $2i-1\in [n_2]$ then $2i-1 + n_1 \in [n]$ and
		\begin{align*}
			\varphi_2(w_0(z_{2i-1})) 
			&= \varphi_2(w_0(i)) \\
			&= \varphi_2(n_2-i+1) \\
			&= n_2 + k - i + 1 \\
			&= (n - 2k + 1 - n_2) + n_2 + k - i + 1 \\
			&= n - (k+i-1) + 1 \\
			&= x_{2(i+k-1)} \\
			&= x_{2i-1 + 2k-1} = x_{2i-1 + n_1}.
		\end{align*}
		If $2i \in [n_2]$ then $2i + n_1 \in [n]$ and
		\begin{align*}
			\varphi_2(w_0(z_{2i})) 
			= \varphi_2(w_0(n_2-i+1)) 
			= \varphi_2(i) 
			= i+k 
			= x_{2(i+k) - 1}
			= x_{2i + n_1}. &\qedhere
		\end{align*} 
	\end{wideenumerate}
\end{proof}

\begin{exa}
	Consider $n=9$, $n_1 = 6$ and $n_2 = 3$.
	Then $k = 3$.
	Using the notation from \Cref{thm:inductive_product_and_x_sequence} we obtain
	\begin{align*}
	 x &= (1,9,2,8,3,7,4,6,5), \\
	 y &= (1,6,2,5,3,4), \\
	 z &= (1,3,2).
	\end{align*}
	Then $x = (\varphi_1(y_1) ,\dots, \varphi_1(y_6), \varphi_2(z_1), \varphi_2(z_2), \varphi_2(z_3))$ as predicted by \Cref{thm:inductive_product_and_x_sequence}.
	Moreover, $x,y$ and~$z$ are the sequences used to define the elements in stair form $\sigma_{(6,3)}$, $\sigma_{(6)}$ and $\sigma_{(3)}$, respectively.
	Therefore,
	\begin{align*}
	 \sigma_{(6,3)} 
	 = (\varphi_1(y_1) ,\dots, \varphi_1(y_6))( \varphi_2(z_1), \varphi_2(z_2), \varphi_2(z_3))
	 = \sigma_{(6)}^{\varphi_1} \sigma_{(3)}^{\varphi_2}
	 =  \sigma_{(6)} \iprod \sigma_{(3)}.
	 \end{align*}
	 This also illustrates the idea of the proof of the next lemma.
\end{exa}

\begin{lem} 
	\label{thm:inductive_product_and_element_in_stair_from}
	Let $\alpha = \parts{\alpha}{l} \vDash_e n$ with $l\geq1$.
	Then we have the following.
	\begin{enumerate}
		\item If $\alpha_1$ is even then $\sigma_\alpha = \sigma_{(\alpha_1)} \iprod \sigma_{(\alpha_2,\dots, \alpha_l)}$.
		\item If $\alpha_1$ is odd then $\sigma_\alpha = \sigma_{(\alpha_1)} \iprod \left(\sigma_{(\alpha_2,\dots, \alpha_l)}\right)^{w_0}$ where $w_0$ is the longest element of $\SG_{\alpha_2 + \dots + \alpha_l}$.
	\end{enumerate}
\end{lem}
\begin{proof} 
	Set   $n_1 := \alpha_1$ and $n_2 := \alpha_2 + \dots + \alpha_l$. 
	As in \Cref{thm:inductive_product_and_x_sequence}, let $x^{(m)}$ be the sequence $(x^{(m)}_1,\dots, x^{(m)}_m)$ given by $x^{(m)}_{2i-1} = i$ and $x^{(m)}_{2i} = m-i+1$ for $m \in \N_0$ and set $x:=x^{(n)}$, $y := x^{(n_1)}$ and  $z := x^{(n_2)}.$
	We have that
	\begin{enumerate}
		\item $\sigma_{\alpha}$ has the cycles
		\begin{align*}
			\sigma_{\alpha_i} = \left( x_{\alpha_1 + \dots + \alpha_{i-1} + 1}, x_{\alpha_1 + \dots + \alpha_{i-1} + 2}, \dots , x_{\alpha_1 + \dots + \alpha_{i-1} + \alpha_i} \right )
		\end{align*}
		for $i  = 1,\dots, l$,
		\item
		$\sigma_{(\alpha_1)} = \left( y_{1}, y_2, \dots , y_{n_1} \right )$ and
		\item $\sigma_{(\alpha_2,\dots, \alpha_l)}$ has the cycles
		\begin{align*}
		\tilde\sigma_{\alpha_i} = \left( z_{\alpha_2 + \dots + \alpha_{i-1} + 1}, z_{\alpha_2 + \dots + \alpha_{i-1} + 2}, \dots , z_{\alpha_2 + \dots + \alpha_{i-1} + \alpha_i} \right )
		\end{align*}
		for $i  = 2,\dots, l$.
	\end{enumerate}
	Assume that $\alpha_1$ is even and set $\sigma := \sigma_{(\alpha_1)} \iprod \sigma_{(\alpha_2,\dots, \alpha_l)}$.
	From \Cref{thm:inductive_product_cycles} we obtain that $\sigma$ has the cycles
		$(\sigma_{(\alpha_1)})^{\varphi_1}$ and $(\tilde \sigma_{(\alpha_i)})^{\varphi_2}$ for $i = 2,\dots, l$.
	By \Cref{thm:inductive_product_and_x_sequence}, $\varphi_1(y_j) = x_j$ for $j\in [n_1]$ and $\varphi_2(z_j) = x_{\alpha_1 + j}$ for $j\in [n_2]$. As a consequence,
	\begin{align*}
		(\sigma_{(\alpha_1)})^{\varphi_1} 
		= (\varphi_1(y_1),\dots, \varphi_1(y_{\alpha_1})) 
		= (x_1,\dots, x_{\alpha_1}) 
		= \sigma_{\alpha_1}
	\end{align*}
	and 
	\begin{align*}
	(\tilde \sigma_{\alpha_i})^{\varphi_2 }
	&=   \left( \varphi_2(z_{\alpha_2 + \dots + \alpha_{i-1} + 1}), \dots , \varphi_2(z_{\alpha_2 + \dots + \alpha_{i-1} + \alpha_i}) \right )\\
	&=  \left( x_{\alpha_1 + \dots + \alpha_{i-1} + 1}, \dots , x_{\alpha_1 + \dots + \alpha_{i-1} + \alpha_i} \right )\\
	&= \sigma_{\alpha_i}
	\end{align*}
	for $i = 2,\dots, l$. Hence, $\sigma = \sigma_\alpha$.
	
	Now let $\alpha_1$ be odd.
	Set $\sigma := \sigma_{(\alpha_1)} \iprod \left(\sigma_{(\alpha_2,\dots, \alpha_l)}\right)^{w_0}$ where $w_0$ is the longest element of $\SG_{\alpha_2 + \dots + \alpha_l}$.
	 Then $\sigma$ has the cycles
	$(\sigma_{(\alpha_1)})^{\varphi_1}$ and $( (\tilde\sigma_{(\alpha_i)})^{w_0})^{\varphi_2}$ for $i = 2,\dots, l$.
	Moreover, from \Cref{thm:inductive_product_and_x_sequence} we have that $\varphi_2(w_0(z_j)) = x_{\alpha_1 + i}$ for $j\in [n_2]$. Thus,
	\begin{align*}
		(\tilde\sigma_{(\alpha_i)})^{w_0})^{\varphi_2}	
		&=   \left( \varphi_2(w_0(z_{\alpha_2 + \dots + \alpha_{i-1} + 1})), \dots , \varphi_2(w_0(z_{\alpha_2 + \dots + \alpha_{i-1} + \alpha_i}) \right ))\\
		&=  \left( x_{\alpha_1 + \dots + \alpha_{i-1} + 1}, \dots , x_{\alpha_1 + \dots + \alpha_{i-1} + \alpha_i} \right )\\
		&= \sigma_{\alpha_i}
	\end{align*}
	for $i = 2,\dots, l$. As we have already shown that	$(\sigma_{(\alpha_1)})^{\varphi_1}= \sigma_{\alpha_1}$,  it follows that $\sigma = \sigma_\alpha$.
\end{proof}

The upcoming \Cref{thm:inductive_product_structures_of_Sigma} is the main result of this  \namecref{sec:equivalence_classes:inductive_product}.
It enables us to decompose $\Sigma_{\alpha}$ if $\alpha_1$ is even.
Before we can state the result, we need to introduce some more notation.
For $\alpha\vDash_e n$ we define 
\begin{align*}
\Sigma_\alpha^\times := \set{\sigma \in \Sigma_\alpha \mid \partition(\sigma) = \partition(\sigma_\alpha)}.
\end{align*}
In  \Cref{thm:inductive_product_structures_of_Sigma} the set $ \left(\Sigma^\times_{\alpha}\right)^{w_0}$ appears where $w_0$ the longest element of $\SG_n$.
Let $\sigma \in \Sigma_\alpha$.
Then by \Cref{thm:conjugation_with_w0_and_equivalence_classes}, $\sigma^{w_0} \in \Sigma_\alpha$.
Since $\partition(\sigma^{w_0}) = w_0(\partition(\sigma))$, we have
\begin{align}
	\label{eq:Sigma_star_and_orbits}
	\sigma \in \left(\Sigma^\times_{\alpha}\right)^{w_0} 
	\iff \partition(\sigma^{w_0}) = \partition(\sigma_\alpha)
	\iff \partition(\sigma) = \partition(\sigma^{w_0}_\alpha).
\end{align}

\begin{thm}
	\label{thm:inductive_product_structures_of_Sigma}
	Let $\alpha = \parts{\alpha}{l} \vDash_e n$ with $l\geq 1$.
	\begin{enumerate}
		\item If $\alpha_1$ is even then $\Sigma_{\alpha} = \Sigma_{(\alpha_1)} \iprod \Sigma_{(\alpha_2,\dots, \alpha_l)}$.
		\item If $\alpha_1$ is odd then $\Sigma^\times_{\alpha} = \Sigma^\times_{(\alpha_1)} \iprod \left(\Sigma^\times_{(\alpha_2,\dots, \alpha_l)}\right)^{w_0}$ where $w_0$ is the longest element of $\SG_{\alpha_2 + \cdots  + \alpha_l}$.
	\end{enumerate}
\end{thm}

\begin{proof}
	Let $\alpha^{(1)} := (\alpha_1)$, $\alpha^{(2)} := (\alpha_2,\dots, \alpha_l)$, $n_1 := |\alpha^{(1)}|$, $n_2 := |\alpha^{(2)}|$ and $w_0$ be the longest element of $\SG_{n_2}$.
	We use the inductive product on $\SG_{n_1} \times \SG_{n_2}$ and the related notation.
	The proofs of~(1) and~(2) have a lot in common. Hence, we do them simultaneously as much as possible and 
	separate the cases $\alpha_1$ even and $\alpha_1$ odd only when necessary.
	
	If $l =1$ then $\alpha = \alpha^{(1)}$, $\alpha^{(2)} = \emptyset$ and thus 
	\begin{align*}
	\Sigma_{\alpha^{(1)}} \iprod \Sigma_{\alpha^{(2)}} = \Sigma_{\alpha} \iprod \SG_0 = \Sigma_\alpha.
	\end{align*}
	Moreover, $\Sigma^\times_{(\alpha_1)} = \Sigma_{(\alpha_1)}$ and $\left(\Sigma^\times_{\emptyset}\right)^{w_0} = \Sigma_{\emptyset}$. Thus we have (1) and (2) in this case.

	Now suppose $l\geq 2$.
	Let $\sigma := \sigma_\alpha$, $\sigma_1 := \sigma_{\alpha^{(1)}}$ and $\sigma_2 := \sigma_{\alpha^{(2)}}$ if $\alpha_1$ is even and $\sigma_2 = \sigma^{w_0}_{\alpha^{(2)}}$ if $\alpha_1$ is odd. 
	From \Cref{thm:inductive_product_and_element_in_stair_from} we have $\sigma = \sigma_1 \iprod \sigma_2$.
	By \Cref{thm:parametrizations_of_kim}, $\sigma_{\alpha^{(i)}} \in \Sigma_{\alpha^{(i)}}$ for $i=1,2$.
	In addition, \Cref{thm:conjugation_with_w0_and_equivalence_classes} then yields that $\sigma^{w_0}_{\alpha^{(2)}} \in \Sigma_{\alpha^{(2)}}$.
	Consequently, $\sigma_i\in \Sigma_{\alpha^{(i)}}$ for $i=1,2$.
	
	We begin with the inclusions  ``$\subseteq$''. Let $\tau \in \Sigma_{\alpha}$ with $\partition(\tau) = \partition(\sigma)$ if $\alpha_1$ is odd.
	First we show $\tau \in \SG_{n_1} \iprod \SG_{n_2}$.
	By \Cref{thm:inductive_product_image_and_injectivity}, we have to show $\tau(N_i) = N_i$ for $i=1,2$. Since $\set{N_1,N_2}$ is a set partition of $[n]$, it suffices to show $\tau(N_1) = N_1$.
 	As $\sigma_1\in \SG_{n_1}$ is an $n_1$-cycle, $\partition(\sigma_1) = \set{[n_1]}$.
	Moreover, \Cref{thm:inductive_product_orbits} yields $\partition(\sigma) = \varphi_1(\partition(\sigma_1)) \cup \varphi_2(\partition(\sigma_2))$.
	Thus,
	\begin{align*}
	N_1 = \varphi_1([n_1]) \in \varphi_1(\partition(\sigma_1)) \subseteq \partition(\sigma).
	\end{align*}
	If $\alpha_1$ is even then $N_1\in \partition_e(\sigma)$.
	Moreover, \Cref{thm:characterization_of_Sigma_using_length} yields $\partition_e(\tau) = \partition_e(\sigma)$. Thus,  $N_1\in \partition(\tau)$ which means that $\tau(N_1) = N_1$.
	If $\alpha_1$ is odd then $\partition(\tau) = \partition(\sigma)$ by assumption. Hence, $N_1\in P(\sigma) = P(\tau)$ and thus $\tau(N_1) = N_1$.
	
	Because $\tau \in \SG_{n_1} \iprod \SG_{n_2}$, there are $\tau_1 \in \SG_{n_1}$ and $\tau_2 \in \SG_{n_2}$ such that $\tau = \tau_1 \iprod \tau_2$. 
	Let $i\in \set{1,2}$.
	We want to show $\tau_i \in \Sigma_{\alpha^{(i)}}$.
	Recall that $\sigma_i  \in \Sigma_{\alpha^{(i)}}$.
	Thus, from \Cref{thm:characterization_of_Sigma_using_length} it follows that $\tau_i \in \Sigma_{\alpha^{(i)}}$ if and only if
	\begin{enumerate}[label = (\roman*)]
	\item $\sigma_i$ and $\tau_i$ are conjugate in $\SG_{n_i}$,
	\item $\ell(\sigma_i) = \ell(\tau_i)$ and
	\item $\partition_e(\sigma_i) = \partition_e(\tau_i)$.
	\end{enumerate}
	Therefore, we show that $\tau_i$ satisfies (i) -- (iii).
	Let~$i$ be arbitrary again.
	\begin{enumerate}[label = (\roman*), wide, nosep]
		\item 
		For a permutation $\xi$, let $\typems(\xi)$ be the multiset of cycle lengths of $\xi$. Assume $\xi = \xi_1 \iprod \xi_2$ for $\xi_i \in \SG_{n_i}$ and $i=1,2$. From \Cref{thm:inductive_product_cycles} it follows that
		\begin{align}
		\label{eq:cycle_type_multisets}
		\typems(\xi) = \typems(\xi_1) \cup \typems(\xi_2).
		\end{align}
		Since $\tau = \tau_1 \iprod \tau_2$,
		\Cref{thm:inductive_product_orbits} implies
		$\partition(\tau) = \varphi_1(\partition(\tau_1)) \cup \varphi_2(\partition(\tau_2))$.
		Therefore, from  $N_1\in \partition(\tau)$ it follows that $\partition(\tau_1) = \set{[n_1]}$. That is, $\tau_1$ is an $n_1$-cycle of $\SG_{n_1}$. 
		By definition, $\sigma_1$ is an $n_1$-cycle of $\SG_{n_1}$ too.
		Thus, $\typems(\tau_1 ) = \typems(\sigma_1)$. Since $\tau \in  \Sigma_\alpha$, $\tau$ and $\sigma$ are conjugate so that $\typems(\tau) = \typems(\sigma)$. Because of \Cref{eq:cycle_type_multisets} and  $\typems(\tau_1 ) = \typems(\sigma_1)$, it follows that also $\typems{(\tau_2)} = \typems{(\sigma_2)}$.
		In other words, $\tau_i$ and $\sigma_i$ are conjugate for $i=1,2$.		
		\item 
		Let $m := \floor*{\frac{n_1}{2}}$.
		By \Cref{thm:inductive_product_and_length}, there are $p,q\leq m$ such that
		\begin{align*}
		\ell(\tau ) = \ell(\tau_1) + \ell(\tau_2) + (p+q)n_2.
		\end{align*}
		Moreover,  we have $\ell(\tau_i) \leq \ell(\sigma_i)$ for $i=1,2$ because $\tau_i$ and $\sigma_i$ are conjugate and $\sigma_i\in \Sigma_{\alpha^{(i)}}$.
		On the other hand,  $\sigma_1$ is oscillating by \Cref{thm:characterization_of_Sigma_(n)} and hence \Cref{thm:inductive_product_and_length} yields
		\begin{align*}
		\ell(\sigma ) = \ell(\sigma_1) + \ell(\sigma_2) + 2mn_2.
		\end{align*}
		Since $\tau \in \Sigma_\alpha$, we have $\ell(\tau) = \ell(\sigma)$.
		Therefore, we obtain from the equalities for $\ell(\tau)$ and $\ell(\sigma)$ and the inequalities for $\ell(\tau_1),\ell(\tau_2), p$ and~$q$ that $\ell(\tau_1) = \ell(\sigma_1)$ and $\ell(\tau_2) = \ell(\sigma_2)$.
		\item 
		\Cref{thm:inductive_product_orbits} states that
		\begin{align}
		\label{eq:inductive_product_orbits}
		\partition(\xi) = \varphi_1(\partition(\xi_1)) \cup \varphi(\partition(\xi_2))
		\end{align}
		for $\xi = \sigma,\tau$. 
		This equality remains valid if we replace $\partition$ by $\partition_e$.
		From $\tau \in \Sigma_\alpha$ and \Cref{thm:characterization_of_Sigma_using_length} it follows that $\partition_e(\tau) = \partition_e(\sigma)$. Hence,
		\begin{align*}
		\varphi_1(\partition_e(\tau_1)) \cup \varphi_2(\partition_e(\tau_2))=
		\varphi_1(\partition_e(\sigma_1)) \cup \varphi_2(\partition_e(\sigma_2)).
		\end{align*}
		Since $\varphi_1$ and $\varphi_2$ are bijections and the images of $\varphi_1$ and $\varphi_2$ are disjoint, it follows that $\partition_e(\tau_i) = \partition_e(\sigma_i)$ for $i=1,2$.
		This finishes the proof of $\tau \in \Sigma_{\alpha^{(1)}} \iprod \Sigma_{\alpha^{(2)}}$.
	\end{enumerate}

	It remains to show that $\tau_1 \in \Sigma^\times_{\alpha^{(1)}}$ and  $\tau_2 \in \left(\Sigma^\times_{\alpha^{(2)}}\right)^{w_0}$ if $\alpha_1$ is odd.
	Thus, assume that $\alpha_1$ is odd. 
	We have already seen that $\partition(\tau_1) = \partition(\sigma_1)$. Hence, $\tau_1 \in \Sigma^\times_{\alpha^{(1)}}$.
	Since $\alpha_1$ is odd, $\partition(\tau) = \partition(\sigma)$ by assumption and therefore we deduce from \Cref{eq:inductive_product_orbits} as above that $\partition(\tau_2) = \partition(\sigma_2)$. Now we can use that $\sigma_2 = \sigma_{\alpha^{(2)}}^{w_0}$ and obtain $\tau_2 \in \left(\Sigma^\times_{\alpha^{(2)}}\right)^{w_0}$ from \Cref{eq:Sigma_star_and_orbits}.

	We continue with the inclusions ``$\supseteq$''. 
	Let $\tau_i \in \Sigma_{\alpha^{(i)}}$ for $i = 1,2$ and $\tau := \tau_1 \iprod \tau_2$. If $\alpha_1$ is odd, assume that in addition $\tau_1 \in \Sigma^\times_{\alpha^{(1)}}$ and  $\tau_2 \in \left(\Sigma^\times_{\alpha^{(2)}}\right)^{w_0}$ which by \Cref{eq:Sigma_star_and_orbits} is equivalent to $\partition(\tau_i) = \partition(\sigma_i)$ for $i=1,2$.
	
	We want to show $\tau \in \Sigma_\alpha$ and again use \Cref{thm:characterization_of_Sigma_using_length} to do this. 
	That is, we show the properties (i) -- (iii) for $\tau$ and $\sigma$.
	\begin{enumerate}[label = (\roman*), wide, nosep]
		\item For $i\in \set{1,2}$ we have $\typems(\tau_i) =\typems(\sigma_i)$ since $\tau_i\in \Sigma_{\alpha^{(i)}}$. Hence, from \Cref{eq:cycle_type_multisets} it follows that  $\typems(\tau) = \typems(\sigma)$, \ie $\tau$ and $\sigma$ are conjugate.
		\item 
		Since $\tau_1,\sigma_1\in \Sigma_{\alpha^{(1)}}$, they are oscillating $n_1$-cycles by \Cref{thm:characterization_of_Sigma_(n)}.
		Therefore, \Cref{thm:inductive_product_and_length} yields
		\begin{align*}
		\ell(\xi ) = \ell(\xi_1) + \ell(\xi_2) + 2mn_2
		\end{align*}
		for $\xi = \sigma,\tau$ and $m = \floor{\frac{n_1}2}$. Moreover, as $\sigma_i,\tau_i \in \Sigma_{\alpha^{(i)}}$, $\ell(\tau_i) = \ell(\sigma_i)$ for $i = 1,2$. 
		Hence, $\ell(\tau) = \ell(\sigma)$.
		\item Since $\xi = \xi_1 \iprod \xi_2$ for $\xi = \sigma, \tau$, Equation \Cref{eq:inductive_product_orbits} holds. This equation remains true if we substitute $\partition$ by $\partition_e$.
		In addition, from \Cref{thm:characterization_of_Sigma_using_length} we obtain that $\partition_e(\tau_i) = \partition_e(\sigma_i)$ for $i=1,2$. Thus, $\partition_e(\tau) = \partition_e(\sigma)$.
	\end{enumerate}
	Because of (i) -- (iii) we can now apply \Cref{thm:characterization_of_Sigma_using_length}  and obtain that $\tau \in \Sigma_\alpha$.
	In the case where $\alpha_1$ is odd, it remains to show $\partition(\tau) = \partition(\sigma)$. But this is merely a consequence of $\partition(\tau_i) = \partition(\sigma_i)$ for $i=1,2$ and \Cref{eq:inductive_product_orbits}.
\end{proof}

We now infer from \Cref{thm:inductive_product_structures_of_Sigma} that the inductive product provides a bijection from $\Sigma_{(\alpha_1)} \times \Sigma_{(\alpha_2,\dots, \alpha_l)}$ to $\Sigma_\alpha$ for all $\alpha \vDash_e n$ with even $\alpha_1$.

\begin{cor}
	\label{thm:indcutive_product_bijection}
	Let $\alpha = \parts{\alpha}{l} \vDash_e n$ with $l\geq 1$.
	\begin{enumerate}
		\item 
		If $\alpha_1$ is even then the map $\Sigma_{(\alpha_1)} \times \Sigma_{(\alpha_2,\dots, \alpha_l)} \to \Sigma_{\alpha}$,  $(\sigma_1,\sigma_2)\mapsto \sigma_1 \iprod \sigma_2$ is a bijection.
		\item
		If $\alpha_1$ is odd then the map $\Sigma^\times_{(\alpha_1)}  \times \left(\Sigma^\times_{(\alpha_2,\dots, \alpha_l)}\right)^{w_0} \to \Sigma^\times_{\alpha}$, $(\sigma_1,\sigma_2)\mapsto \sigma_1 \iprod \sigma_2$ where $w_0$ is the longest element of $\SG_{\alpha_2 + \cdots  + \alpha_l}$ is a bijection.
	\end{enumerate}
\end{cor}

\begin{proof}
	By \Cref{thm:inductive_product_image_and_injectivity} the two maps in question are injective. 
	\Cref{thm:inductive_product_structures_of_Sigma} shows that they are also surjective.
\end{proof}

Recall that, given a maximal composition $\alpha = \parts{\alpha}l\vDash_e n$, there exists $0\leq j \leq l$ such that $\alpha_1,\dots, \alpha_j$ are even and $\alpha_{j+1} \geq \dots \geq \alpha_l$ are odd. 
Using Part~(1) of \Cref{thm:indcutive_product_bijection} iteratively, we obtain the following decomposition of the elements of $\Sigma_\alpha$.

\begin{cor}
	\label{thm:inductive_product_splitting_even_entries_off}
	Let $\alpha = \parts{\alpha}l\vDash_e n$,
	$\sigma\in \SG_n$ of type $\alpha$ 
	and
	$0\leq j \leq l$ be such that $\alpha' := (\alpha_{j+1}, \dots, \alpha_{l})$ are the odd parts of $\alpha$.
	Then $\sigma \in \Sigma_\alpha$ if and only if there are $\sigma_i\in \Sigma_{(\alpha_i)}$ for $i = 1,\dots, j$ and $\tau \in \Sigma_{\alpha'}$ such that
	\begin{align*}
		\sigma = \sigma_1 \iprod \sigma_2 \iprod  \cdots  \iprod \sigma_j \iprod \tau 
	\end{align*}
	where the product is evaluated from right to left.
\end{cor}

\begin{exa}
	Consider $\alpha = (2,4,3,1,1)\vDash_e 11$. From \Cref{tbl:Sigma_(n)} and \Cref{thm:odd_hook_example_(3_1_1)} we obtain
	\begin{align*}
	\Sigma_{(2)} &= \set{(1,2)}, \\
	\Sigma_{(4)} &= \set{(1,4,2,3), (1,3,2,4)}, \\
	\Sigma_{(3,1,1)} &= \set{(1,5,2),(1,2,5),(1,5,3),(1,3,5),(1,5,4),(1,4,5)}.
	\end{align*}
	By \Cref{thm:inductive_product_splitting_even_entries_off}, $\Sigma_\alpha$ consists of all elements $(1,2) \iprod \left( \sigma \iprod \tau \right)$ with $\sigma \in \Sigma_{(4)}$ and $\tau \in \Sigma_{(3,1,1)}$. Thus, $|\Sigma_\alpha|= 12$.
	For instance,
	\begin{align*}
	(1,2) \iprod \left( (1,3,2,4) \iprod (1,3,5)\right) 
	&= (1,2)\iprod (1,8,2,9)(3,5,7) \\
	&= (1,11)(2,9,3,10)(4,6,8)
	\end{align*}
	is an element of $\Sigma_{\alpha}$.
\end{exa}

\begin{rem}
	\label{thm:inductive_product_remark_reduction_to_odd_partitions}
	For compositions with one part $\alpha = (n)$, \Cref{thm:characterization_of_Sigma_(n)} provides a combinatorial characterization of $\Sigma_{(n)}$.
	Therefore, \Cref{thm:inductive_product_splitting_even_entries_off} reduces the problem of describing $\Sigma_\alpha$ for each maximal composition $\alpha$ to the case where $\alpha$ has only odd parts. 
	These $\alpha$ are the partitions consisting of odds parts.
	
	If $\alpha$ is an odd hook, then \Cref{thm:characterization_of_Sigma_for_odd_hook} yields that the hook properties characterize the elements of $\Sigma_\alpha$.
	That is, we have a description of $\Sigma_\alpha$ for all maximal compositions $\alpha$ whose odd parts form a hook.
\end{rem}

Let $\alpha \vDash_e n$ and $\alpha'$ be the composition formed by the odd parts of $\alpha$.
We infer from \Cref{thm:inductive_product_splitting_even_entries_off}  a formula that expresses $\card{\Sigma_\alpha}$ as a product of $\card{\Sigma_{\alpha'}}$ and a factor that only depends on the even parts of $\alpha$.
In the case where $\alpha'$ is an odd hook, we can determine $\card{\Sigma_{\alpha'}}$ explicitly and thus obtain a closed formula.

\begin{cor}
	\label{thm:inductive_product_cardinality_formula}
	Let $\alpha = \parts{\alpha}l\vDash_e n$,
	$0\leq j \leq l$ be such that $(\alpha_1,\dots, \alpha_j)$ are the even and  $\alpha' := (\alpha_{j+1}, \dots, \alpha_{l})$ are the odd parts of $\alpha$, 
	$n' := |\alpha'|$,
	$P := \set{i\in [j] \mid \alpha_i \geq 4}$, 
	$p := |P|$
	and
	$q := -2p + \frac{1}{2}\sum_{i\in P}\alpha_i$.
	Then
	\begin{align*}
	|\Sigma_\alpha| &= 2^p 3^{q}|\Sigma_{\alpha'}|. 
	\end{align*}
	Moreover, if $\alpha'$ is a hook $(r, 1^{n'-r})$ then
	\begin{align*}
		|\Sigma_\alpha| &= 
		\begin{cases}
			 2^p 3^{q} &\text{if $r \leq 1$} \\
		 	(n'-r+1)2^{p'} 3^{q'} &\text{if $r\geq3$}
		\end{cases}
	\end{align*}
	where $p' := p+1$ and $q' :=q+\frac{r-3}{2}$. 
\end{cor}

\begin{proof}
	Since $\alpha_1, \dots, \alpha_j$ are the even parts of $\alpha$,  \Cref{thm:inductive_product_splitting_even_entries_off} implies that
	\begin{align}
	\label{eq:inductive_product_cardinality_product_rule}
		|\Sigma_\alpha| = \left|\Sigma_{\alpha'}\right| \prod_{i=1}^j |\Sigma_{(\alpha_i)}|.
	\end{align}
	For the same reason, \Cref{thm:sizes_of_Sigma_n} yields
	\begin{align*}
			|\Sigma_{(\alpha_i)}| = 
		\begin{cases}
		1 & \myif \text{$n\leq2$} \\
		2 \cdot 3^{\frac{\alpha_i-4}2}  & \myif \text{$n\geq 4$}. \\
		\end{cases}
	\end{align*}
	for $i = 1,\dots, j$.
	Therefore, 
	\begin{align*}
	\prod_{i=1}^j |\Sigma_{(\alpha_i)}| 
	= \prod_{i\in P} 2 \cdot 3^{\frac{\alpha_i-4}2} 
	= 2^p 3^{-2p + \frac{1}{2}\sum_{i\in P}\alpha_i}
	= 2^p 3^q.
	\end{align*}
	and with \Cref{eq:inductive_product_cardinality_product_rule} we get the first statement.
	
	For the second part, assume that $\alpha'$ is a hook.
	Then, by the choice of~$j$, $\alpha'$ is an odd hook.
	It remains to compute $|\Sigma_{\alpha'}|$.
	If $\alpha' = \emptyset$ or $\alpha' = (1^{n'})$ we have $|\Sigma_\alpha'|=1$.
	If $\alpha' = (r,1^{n'-r})$ with $r\geq 3$ then \Cref{thm:odd_hook_cardinality_of_Sigma_alpha} provides the formula 
	\begin{equation*}
	|\Sigma_{\alpha'}| = 
	2 (n'-r+1) 3^{\frac{r-3}2}.
	\qedhere
	\end{equation*}
\end{proof}

\begin{exa}
	Consider $\alpha = (2,8,4,5,1,1,1) \vDash_e 22$.
	Then $\alpha' = (5,1,1,1) \vDash_e 8$ is a hook,
	$P = \set{2,3}$,
	$p' = 2 + 1$
	and
	$q' = - 2\cdot 2 + \frac{1}{2}(8+4)  + \frac{5-3}{2} = 3$.
	Thus, \Cref{thm:inductive_product_cardinality_formula} yields
	$|\Sigma_\alpha| = (8-5+1) 2^3 3^3 = 864$.
\end{exa}

Let $\alpha = (l, 1^{n-l})\vDash_e n$ be a hook.
From \Cref{thm:odd_hook_bijection} we know how to construct $\Sigma_\alpha$ from $\Sigma_{(l)}$ if~$k$ is odd.
If~$l$ is even, we obtain $\Sigma_\alpha$ in the following way.

\begin{cor}
	\label{thm:even_hook_bijection}
	Let $\alpha = (l,1^{n-l})\vDash_e n$ be an even hook and $\id \in \SG_{n-k}$.
	Then the map $\Sigma_{(l)} \to \Sigma_\alpha$, $\sigma \mapsto \sigma \iprod \id$ is a bijection. 
\end{cor}

\begin{proof}
	Recall $\Sigma_{(1^{n-l})} = \set{\id}$. 
	Then \Cref{thm:indcutive_product_bijection} yields that the map from the claim is a bijection.
\end{proof}

\begin{exa}
	Consider $\alpha = (4,1,1)$ and $\id \in \SG_2$. From \Cref{tbl:Sigma_(n)} we read
	\begin{align*}
	\Sigma_{(4)} &= \set{(1,4,2,3), (1,3,2,4)}
	\end{align*}
	Hence, \Cref{thm:even_hook_bijection} yields
	\begin{align*}
		\sigma_\alpha 
		=\set{\sigma \iprod \id \mid \sigma \in \Sigma_{(4)}}
		=\set{ (1,6,2,5), (1,5,2,6) }.	
	\end{align*}
\end{exa}

In \Cref{thm:characterization_of_Sigma_for_odd_hook} we showed that $\Sigma_\alpha$ is characterized by the hook properties if $\alpha$ is an odd hook.
In the remainder of the \namecref{sec:equivalence_classes} we want to prove that the same is true for even hooks.
We first show that $\iprod$ is compatible with the concepts of being oscillating and having connected intervals.

\begin{lem}
	\label{thm:inductive_product_osc_and_c.i}
	Let $\sigma_1\in \SG_{n_1}$, $\sigma_2\in \SG_{n_2}$ and $\sigma := \sigma_1 \iprod \sigma_2$. Then $\sigma$ is oscillating (has connected intervals) if and only if $\sigma_1$ and $\sigma_2$ are oscillating (have connected intervals).
\end{lem}

\begin{proof}
	Let  $\sigma_r = \sigma_{r,1}\sigma_{r,2}\cdots \sigma_{r,p_r}$  be a decomposition in disjoint cycles for $r= 1,2$.
	Fix an $r\in \set{1,2}$ and a cycle $(c_1,\dots c_t) = \sigma_{r,j}$ of $\sigma_r$.
	Then by \Cref{thm:inductive_product_cycles} we have that
	\begin{align*}
		\sigma_{r,j}^{\varphi_r} = (\varphi_r(c_1),\dots, \varphi_r(c_t)).
	\end{align*}
	As $\varphi_r$ is strictly increasing, it preserves the relative order of the cycle elements so  that
	\begin{align*}
		\cst( \sigma_{r,j}) = \cst(\sigma_{r,j}^{\varphi_r}).
	\end{align*}
	In addition, \Cref{thm:inductive_product_cycles} provides the cycle decomposition
	\begin{align*}
				\sigma = \sigma^{\varphi_1}_{1,1} \cdots \sigma^{\varphi_1}_{1,p_1} 
				\cdot
				\sigma^{\varphi_2}_{2,1} \cdots \sigma^{\varphi_2}_{2,p_2}. 
		\end{align*}
	of $\sigma$.
	Hence, $\sigma$ is oscillating if and only $\sigma_1$ and $\sigma_2$ are oscillating.
	For the same reason, $\sigma$ has connected intervals if and only if $\sigma_1$ and $\sigma_2$ have connected intervals.
\end{proof}

We now generalize \Cref{thm:characterization_of_Sigma_for_odd_hook} to all hooks. 
The hook properties can be looked up in \Cref{def:hook_properties}.

\begin{thm}
	\label{thm:characterization_of_Sigma_for_arbitrary_hook}
	Let $\alpha\hdash n$ be a hook and $\sigma\in \SG_n$ of type $\alpha$. Then 
	$\sigma \in \Sigma_\alpha$ if and only if $\sigma$ satisfies the hook properties.
\end{thm}
\begin{proof}
	Let $\alpha = (l, 1^{n-l})\hdash n$ and $\sigma \in \SG_n$ be of type $\alpha$.
	The case where~$l$ is odd was done in \Cref{thm:characterization_of_Sigma_for_odd_hook}.
	Therefore, assume that~$l$ is even.
	If $l = n$ then the third hook property is satisfied and therefore the $n$-cycle $\sigma \in \SG_n$ has the hook properties if and only if it is oscillating and has connected intervals. 
	By \Cref{thm:characterization_of_Sigma_(n)} this is equivalent to $\sigma \in \Sigma_{(n)}$.
	Therefore we now assume $l < n$.
	Write $\sigma = (d_1,\dots, d_l)$ omitting the trivial cycles. 
	We consider the inductive product on $\SG_l \times \SG_{n-l}$ and $\id \in \SG_{n-l}$.
	Following \Cref{not:inductive_product} we then have that
	\begin{align*}
		N_1  = \varphi_1([l]) = \left[\frac{l}2\right] \cup \left[n - \frac{l}2 +1, n\right].
	\end{align*}
	Note that $\sigma$ satisfies the third hook property if and only if $\set{d_1, \dots, d_l} = N_1$.
	
	We begin with the implication form left to right.
	Assume that $\sigma \in \Sigma_{(l,1^{n-l})}$. 
	By \Cref{thm:even_hook_bijection} there is $\tau \in \Sigma_{(l)}$ such that $\sigma = \tau \iprod \id$.
	Certainly $\id$ is oscillating and has connected intervals.
	Moreover, $\tau$ has these properties by \Cref{thm:characterization_of_Sigma_(n)}.
	Therefore, $\sigma$ is oscillating with connected intervals by \Cref{thm:inductive_product_osc_and_c.i}.
	Because $\sigma = \tau \iprod \id$, \Cref{thm:inductive_product_cycles} implies that we can write $\tau = (c_1,\dots, c_l)$ such that $d_i = \varphi_1(c_i)$ for $i = 1,\dots, l$.
	Therefore,
	\begin{align*}
		\set{d_1,\dots, d_l} 
		= \varphi_1(\set{c_1,\dots c_l})
		= \varphi_1([l]) 
		= N_1
	\end{align*}
	which means that $\sigma$ satisfies the third hook property.
	
	We now show the implication from right to left.
	Assume that $\sigma$ fulfills the hook properties. 
	Then the third hook property yields that
	$\set{d_1,\dots, d_l} = N_1$ which implies that $\sigma(N_1) = N_1$.
	Therefore, $\sigma \in \SG_l \iprod \SG_{n-l}$ by \Cref{thm:inductive_product_image_and_injectivity}, \ie there are $\sigma_1\in \SG_l$ and $\sigma_2 \in \SG_{n-l}$ such that $\sigma = \sigma_1 \iprod \sigma_2$.
	From \Cref{thm:inductive_product_two_domains} we obtain that $\sigma|_{N_1}= \sigma_1^{\varphi_1}$ so that we can write $\sigma_1$ as $\sigma_1 = (c_1 ,\dots, c_l)$ with $c_i = \varphi_1\inv(d_i)$ for $i = 1,\dots, l$.
	It follows that 
	$\sigma_1$ is an $l$-cycle of $\SG_l$.
	Since $\sigma$ fixes each element of $N_2$, it follows from \Cref{thm:inductive_product_two_domains} that $\sigma_2 = \id$.
	As $\sigma$ is oscillating with connected intervals, \Cref{thm:inductive_product_osc_and_c.i} implies that $\sigma_1$ has these properties as well.
	Thus, $\sigma_1 \in \Sigma_{(l)}$ by \Cref{thm:characterization_of_Sigma_(n)}. 
	Hence, we can apply \Cref{thm:even_hook_bijection} and obtain that
	\begin{displaymath}
		\sigma 
		= \sigma_1 \iprod \id \in \Sigma_{(l,1^{n-l})}.
		\qedhere
	\end{displaymath}
\end{proof}

\begin{rem}
	\label{thm:remark_on_generalisation_of_description_of_Sigma_alpha}
	In \Cref{thm:inductive_product_remark_reduction_to_odd_partitions} we  reduced the problem of describing $\Sigma_\alpha$ for all maximal compositions $\alpha$ to the partitions with only odd parts. 
	As we have such a description for odd hooks, it remains to find a combinatorial description of $\Sigma_\alpha$ in the case where $\alpha$ is a partition of odd parts which is not a hook.
	Then $\Sigma_\alpha$ consists of all permutations of type $\alpha$ of maximal length.
	Unfortunately, the situation is a lot more complex.
	One reason for this is the following.
	For any subset $\Sigma$ of $\SG_n$ define
	\begin{align*}
	\partition(\Sigma) := \set{\partition(\sigma) \mid \sigma \in \Sigma}.
	\end{align*}
	
	In general, $\partition(\sigma_\alpha)$ is not the only element of $\partition(\Sigma_\alpha)$ and there seems to be no obvious way to describe $P(\Sigma_\alpha)$.
	Moreover, the number of $\sigma \in \Sigma_\alpha$ whose orbits yield the same set partition of $[n]$ depends on this very set partition.
	For example, $\Sigma_{(3,3)}$ consists of the following elements where elements with the same orbit partition occur in the same row.
	\begin{align*}
	\begin{array}{llll}
	(1,6,2)(3,4,5) & (1,2,6)(3,4,5) & (1,6,2)(3,5,4) &
	(1,2,6)(3,5,4) \\
	(1,6,3)(2,4,5) & (1,6,3)(2,5,4) & (1,3,6)(2,4,5) &
	(1,3,6)(2,5,4) \\
	(1,4,5)(2,6,3) & (1,5,4)(2,3,6) & (1,5,4)(2,6,3) &
	(1,4,5)(2,3,6) \\
	(1,6,4)(2,3,5) & (1,4,6)(2,3,5) & (1,6,4)(2,5,3) &
	(1,4,6)(2,5,3) \\
	(1,6,5)(2,3,4) & (1,5,6)(2,3,4) & (1,5,6)(2,4,3) &
	(1,6,5)(2,4,3) \\
	(1,5,3)(2,4,6) & (1,3,5)(2,6,4) \\
	\end{array}
	\end{align*}	
\end{rem}

\bibliographystyle{amsalpha}
\bibliography{../../bibtex/dislit}

\end{document}